\documentclass[10pt]{article}
\usepackage[english]{babel}
\usepackage{amsfonts}
\usepackage{amsthm}
\usepackage{amsmath}    
\usepackage{bm}         
\usepackage[utf8]{inputenc}
\usepackage{libertine}
\usepackage{graphicx}   
\usepackage{verbatim}   
\usepackage{color}      
\usepackage{subfigure}  
\usepackage{epsfig}
\usepackage{tikz}       
\usepackage{stmaryrd}
\usepackage{hyperref}   
\numberwithin{equation}{section} 
\usepackage{cleveref} 
\usepackage{fancyhdr} 
\usepackage{rotating} 
\usepackage{mathrsfs}
\usepackage{thmtools}
\usepackage{caption}
\usepackage{pgfplots,xparse}
\pgfplotsset{compat=newest}

\declaretheoremstyle[notefont=\bfseries,notebraces={}{},%
    headpunct={},postheadspace=1em]{mystyle}
\declaretheorem[style=mystyle,numbered=no,name=Theorem]{thm-hand}

\declaretheoremstyle[notefont=\bfseries,notebraces={}{},%
    headpunct={},postheadspace=1em]{mystyle}
\declaretheorem[style=mystyle,numbered=no,name=Proposition]{prop-hand}

\declaretheoremstyle[notefont=\bfseries,notebraces={}{},%
    headpunct={},postheadspace=1em]{mystyle}
\declaretheorem[style=mystyle,numbered=no,name=Corollary]{Cor-hand} 

\declaretheoremstyle[notefont=\bfseries,notebraces={}{},%
    headpunct={},postheadspace=1em]{mystyle}
\declaretheorem[style=mystyle,numbered=no,name=Lemma]{lem-hand}

\declaretheoremstyle[notefont=\bfseries,notebraces={}{},%
    headpunct={},postheadspace=1em]{mystyle}
\declaretheorem[style=mystyle,numbered=no,name=Remark]{rem-hand}

\usepackage{longtable}
\usepackage{booktabs}

\makeatletter
\renewcommand*\env@matrix[1][\arraystretch]{%
  \edef\arraystretch{#1}%
  \hskip -\arraycolsep
  \let\@ifnextchar\new@ifnextchar
  \array{*\c@MaxMatrixCols c}}
\makeatother

\renewcommand{\thefigure}{\arabic{section}.\arabic{figure}}

\allowdisplaybreaks

\newcommand{\avg}[1]{\ensuremath{\{\!\!\{#1\}\!\!\}} }
\newcommand{\jump}[1]{\ensuremath{[\![#1]\!]} }

\newcommand {\R}{\mathbb{R}}
\newcommand {\N}{\mathbb{N}}

\newcommand{\matx}[1]{\mathcal{\mathcal{#1}}} 
\newcommand{\blockmatx}[1]{\boldsymbol{\mathcal{#1}}}
\newcommand{\blockvec}[1]{\overset{\leftrightarrow}{#1}}
\newcommand{\blockveccon}[1]{\overset{\leftrightarrow}{\tilde{#1}}}

\newcommand{\interpolation}[2]{\mathbb{I}^{#1}\!#2}
\newcommand{\Dprojection}[2]{\vec{\mathbb{D}}^{#1}\!#2}

\newcommand{\pderivative}[2]{\frac{\partial #1}{\partial #2}}

\begin{document}




\title{ \Large Entropy Stable Discontinuous Galerkin Schemes on Moving Meshes with Summation-by-Parts Property for Hyperbolic Conservation
Laws}

\setcounter{footnote}{0}
\font\myfont=cmr12 at 10pt 
\author{  {\myfont Gero Schnücke\footnote{\ Department for Mathematics and Computer Science, University of Cologne, Email: gschnuec@math.uni-koeln.de}, 
\ Nico Krais\footnote{\ Institute of Aerodynamics and Gas Dynamics (IAG),  University of Stuttgart, Email: krais@iag.uni-stuttgart.de},  
\ Thomas Bolemann\footnote{\ Department for Mathematics and Computer Science, University of Cologne 
\newline Institute of Aerodynamics and Gas Dynamics (IAG), University of Stuttgart, Email: bolemann@iag.uni-stuttgart.de} \ and     
Gregor J. Gassner \footnote{\ Department for Mathematics and Computer Science, University of Cologne, Email:   ggassner@math.uni-koeln.de 
\newline 
Center for Data and Simulation Science (CDS), University of Cologne, Url: www.cds.uni-koeln.de}}} 

\pagestyle{fancy}
\lhead{G. Schn\"ucke, N. Krais, T. Bolemann and G. J. Gassner}
\rhead{Entropy Stable DG Schemes on Moving Meshes}
\date {\myfont \today}
\maketitle
{\hrule height 1.5pt} 
\begin{abstract}
\noindent
This work is focused on the  entropy analysis of a semi-discrete nodal discontinuous Galerkin spectral element method (DGSEM) on moving meshes for hyperbolic conservation laws. The DGSEM is constructed with a local tensor-product Lagrange-polynomial basis computed from Legendre-Gauss-Lobatto (LGL) points. Furthermore, the collocation of interpolation and quadrature nodes is used in the spatial discretization. This approach leads to discrete derivative approximations in space that are summation-by-parts (SBP) operators. On a static mesh, the SBP property and suitable two-point flux functions, which satisfy the entropy condition from Tadmor, allow to mimic results from the continuous entropy analysis on the discrete level. In this paper, Tadmor's condition is extended to the moving mesh framework. Based on the moving mesh entropy condition, entropy conservative two-point flux functions for the homogeneous shallow water equations and the compressible Euler equations are constructed. Furthermore, it will be proven that the semi-discrete moving mesh DGSEM is an entropy conservative scheme when a two-point flux function, which satisfies the moving mesh entropy condition, is applied in the split form DG framework. This proof does not require any exactness of quadrature in the spatial integrals of the variational form. Nevertheless, entropy conservation is not sufficient to tame discontinuities in the numerical solution and thus the entropy conservative moving mesh DGSEM is modified by adding numerical dissipation matrices to the entropy conservative fluxes. Then, the method becomes entropy stable such that the discrete mathematical entropy is bounded at any time by its initial and boundary data when the boundary conditions are specified appropriately. 
 
Besides the entropy stability, the time discretization of the moving mesh DGSEM will be investigated and 
it will be proven that the moving mesh DGSEM  satisfies the free stream preservation property for an arbitrary $s$-stage Runge-Kutta method. 
 
The theoretical properties of the moving mesh DGSEM will be validated by numerical experiments for the compressible Euler equations. 
\vspace{1\baselineskip}
{\hrule height 1.5pt}
\end{abstract}

\newpage
\section{Introduction}\label{sec:Int}

A lot of applications in engineering and physics  
require the approximation of conservation laws on time-dependent domains, e.g., domains with moving boundaries.  
For instance, moving mesh discontinuous Galerkin (DG) methods have been investigated in~\cite{Boscheri2017, Fu2018, Lomtev1999JCP, Nguyen2010JFS, Persson2009}. 
In particular,  moving mesh discontinuous Galerkin spectral element methods (DGSEM) have been constructed and analyzed in~\cite{Kopriva2016, Minoli2011, Winters2014}. In the literature, there are also moving mesh methods with the capability to change the connectivity of the mesh, e.g. with finite volume (FV) methods \cite{Loubere2010,Springel2010} and with a DG method~\cite{Persson2015}. In general, moving mesh methods are well suited to preserve motion related properties like the Galilean-invariance. These properties are necessary to describe physical processes like the formation of disc galaxies \cite{Marinacci2013}. 

A common way to approximate conservation laws on time-dependent domains is to use the Arbitrary Lagrangian-Eulerian (ALE) approach \cite{Donea2004}. 
In this approach the conservation law is transformed from the time-dependent domain onto a time-independent reference domain. The motion of the mesh on the physical domain is part of the transformation. Thus, the grid velocity field appears as a new quantity in the equation on the reference domain. On the one hand the ALE transformation simplifies the discretization, since a static mesh can be used in the reference domain. On the other hand, the new quantities in the equation on the reference domain complicate the discrete stability analysis, even in the linear case \cite{Kopriva2016}. More precisely, the appearance of the grid velocity field leads to a non-conservative equation. On the continuous level, the geometric quantities and the grid velocity field satisfy an additional balance law. This balance law is often referred to as geometric conservation law (GCL) \cite{Farhat2001, Farhat2000, Lesoinne1996, Lombard1979, Mavriplis2006}. The GCL allows to write the transformed equation in conservation form. However, on the discrete level a numerical scheme does not necessarily satisfy a discrete GCL (D-GCL). Farhat et al. \cite{Farhat2001, Farhat2000, Lesoinne1996} proved that the absence of this property has a critical effect on the accuracy and stability of a moving mesh method. In particular, the preservation of constant states is no longer guaranteed. 

In this work, moving mesh DGSEM to solve non-linear, symmetrizable and hyperbolic systems of conservation laws are investigated. It is well known that symmetrizable systems are equipped with an entropy/entropy flux pair  \cite{Godunov1961,Mock1980}. For scalar conservation laws, entropy admissible criteria provide the unique physical relevant weak solution \cite{DiPerna1983,Kruzkov1970}. For systems of conservation laws, entropy admissible criteria are not enough to ensure well posedness \cite{Chiodaroli2015}. Nevertheless, the entropy is an essential quantity to analyze systems of conservation laws. In particular, for gas dynamics a possible mathematical entropy is the scaled negative thermodynamic entropy which shows that the mathematical model correctly captures the second law of thermodynamics \cite{Barth2018}. The entropy is conserved for smooth solutions of a conservation law and decays for discontinuous solutions \cite{Harten1983,Tadmor2003}. 

Therefore, it is desirable to construct numerical schemes for conservation laws which reflect the properties of the entropy on the discrete level. Tadmor \cite{Tadmor1987} developed a discrete entropy criterion to construct a specific class of two-point flux functions for low-order FV methods. FV methods with these class of two-point flux functions preserve entropy on the discrete level. Moreover, these FV methods can be modified by adding dissipation to the numerical fluxes such that the entropy is decreasing for all times. Therefore, two-point fluxes with Tadmor's discrete entropy condition are called entropy conservative fluxes. Fisher and Carpenter \cite{Fisher2013} proved that low-order entropy conservative fluxes can be used to construct high-order schemes when the derivative approximations in space are summation-by-parts (SBP) operators. A SBP operator provides a discrete analogue of the integration-by-parts formula \cite{Fernandez2014,Gassner2013,kreiss1}. It is worth to mention that the derivative matrix in the DGSEM provides a SBP operator, if the tensor-product Lagrange-polynomial basis is computed from Legendre-Gauss-Lobatto (LGL) points and interpolation and quadrature are collocated. Gassner et al. \cite{Gassner2016, Gassner2017} showed that split forms of the partial differential equations can be discretely recovered when specific choices of numerical volume fluxes in the flux form volume integral of Fisher and Carpenter are chosen.
Thus, the following restrictions need to be satisfied to construct entropy stable DGSEM:  
\begin{itemize}
\item[(1)] The derivative matrix satisfies the SBP property.
\item[(2)] There are two-point flux functions with Tadmor's discrete entropy condition that can be extended to high-order in a split form DG framework. 
\end{itemize}
This methodology has been used in the construction of high-order entropy stable DGSEM on quadrilateral/hexahedral elements, e.g.  \cite{Bohm2018,Gassner2016,Wintermeyer2017}, or on triangular/tetrahedral elements, e.g. \cite{chan2018,Chen2017,crean2018}. All these methods are provably entropy stable and the semi-discrete entropy analysis for them is based merely on the properties of the SPB operators and the assumption that the time integration is exact. The exactness of quadrature in the spatial integrals of the variational form is not necessary. 

Up to now, the methodology to construct entropy stable DGSEM has been developed for static meshes. The focus of this work is the construction of entropy stable moving mesh DGSEM. The remainder of the paper is organized as follows: In Section \ref{Sec:BasicConcepts}, the ALE transformation is presented for one dimensional systems of conservation laws. The continuous entropy analysis for the conservation law on the reference domain is given in the Subsection \ref{1DEntropyAnalysisALE}. Then, in Subsection \ref{Sec:EntropyMovingMeshFV}, a moving mesh FV scheme is constructed and a natural extension of Tadmor's discrete entropy condition for the ALE framework is presented. The capability of the moving mesh FV scheme to preserve constant states is investigated in the Subsection \ref{sec:FreestreamPreservationFV}. The ALE transformation and continuous entropy analysis is presented in the Subsections \ref{Sec:ComponentsALE2D}, \ref{Sec:TransCL} and \ref{EntropyAnalysis}. The framework for the spectral element discretization with the SPB operator is given in the Subsection \ref{sec:Spectral} and the DG split form framework is presented in the Subsection \ref{sec:SemiDiscreteDG}. The moving mesh DGSEM is finally presented in the Subsection \ref{sec:SemiDiscreteDG}. A discrete entropy analysis for the moving mesh DGSEM is given in the Subsections  \ref{sec:EntropyConservation} and \ref{sec:EntropyStability}. Furthermore, in Subsection \ref{sec:FreestreamPreservationDGSEM} it is proven that the moving mesh DGSEM satisfies the free stream preservation property. In Section 4,  numerical examples with the compressible Euler equations are presented to validate our theoretical findings.
   
\section{Basic concepts for 1D conservations laws in the ALE framework}\label{Sec:BasicConcepts}
In this section, we discuss the ALE transformation for the one dimensional symmetrizable system     
\begin{equation}\label{eq:consLaw1D}
\pderivative{\textbf{u}}{t} + \pderivative{\textbf{f}}{x}  = \textbf{0}.
\end{equation}
The system consists of $p$ equations, is defined on a time dependent interval $\left[a\left(t\right),b\left(t\right)\right]$ and is equipped with suitable boundary conditions. The time dependent interval has continuous differentiable boundary functions and can be mapped on the reference interval $\left[-1,1\right]$  by a time-dependent affine linear mapping $x=\chi\left(\xi,t\right)$. The mapping provides the quantities 
\begin{equation}\label{GeometricQuantities1D}
J\left(\xi,t\right):=\pderivative{}{\xi}\left(\chi\left(\xi,t\right)\right), 
\qquad \text{and} \qquad \nu\left(\xi,t\right):=\frac{d}{dt}\left(\chi\left(\xi,t\right)\right),
\end{equation}
where $\nu$ is the grid velocity. The chain rule and \eqref{GeometricQuantities1D} provide 
\begin{equation}\label{eq:ChainRule1D}
J\frac{d\mathbf{u}}{dt}=J\pderivative{\mathbf{u}}{t}+\nu\pderivative{\mathbf{u}}{\xi}\qquad\Leftrightarrow\qquad J\pderivative{\mathbf{u}}{t}=J\frac{d\mathbf{u}}{dt}-\nu\pderivative{\mathbf{u}}{\xi}.
\end{equation}
Hence, by applying the chain rule and rearranging terms the conservation law \eqref{eq:consLaw1D} becomes
\begin{equation}\label{eq:consLaw1Da}
J\frac{d\textbf{u}}{dt}+\pderivative{\textbf{f}}{\xi}=\nu\pderivative{\textbf{u}}{\xi}
\end{equation}
on the reference element. The product rule allows to write the equation \eqref{eq:consLaw1Da} as       
\begin{equation}\label{eq:consLaw1Db}
J\frac{d\textbf{u}}{dt}+\pderivative{\textbf{g}}{\xi}=-\pderivative{\nu}{\xi}\mathbf{u},
\end{equation}
where $\textbf{g}=\mathbf{f}-\nu\mathbf{u}$. The equations \eqref{eq:consLaw1Da} and \eqref{eq:consLaw1Db} are not desirable for a numerical discretization, since the equations are not in conservation form and a source term appears. The discretization of equations with source terms provides additional numerical problems, since the preservation of certain states is not guaranteed \cite{Chandrashekar2017, Gassner2016WB, Wintermeyer2017}. Nevertheless, Schwarz's theorem on the symmetry of second derivatives provides  
\begin{equation}\label{GCL1D}
\frac{dJ}{dt}=\frac{d}{dt}\pderivative{\chi}{\xi}=\pderivative{}{\xi}\frac{d\chi}{dt}=\pderivative{\nu}{\xi}.
\end{equation}
The equation \eqref{GCL1D} is typically referred to as the GCL. Thus, equation \eqref{eq:consLaw1D} can be written on the reference interval as the conservation law 
\begin{equation}\label{eq:consLaw1DRef}
\frac{d J\mathbf{u}}{dt}+\pderivative{\mathbf{g}}{\xi}=\mathbf{0}.
\end{equation} 

\subsection{Entropy analysis in one dimension}\label{1DEntropyAnalysisALE}
The well known results by Godunov \cite{Godunov1961} and Mock  \cite{Mock1980} state that a system of conservation laws is symmetrizable if and only if it is equipped with a strictly convex entropy function. Hence, it exists an entropy/entropy flux pair $\left(s,f^{s}\right)$ for the system \eqref{eq:consLaw1D}. Then the entropy variables are defined by  
\begin{equation}\label{eq:1DEntropyVariables}
\textbf{w} := \pderivative{s}{\textbf{u}},
\end{equation}
and the entropy flux satisfies 
\begin{equation}\label{eq:1DEntropyFlux}
\textbf{w}^{T}\left(\pderivative{\textbf{f}}{x}\right)= 
\left(\pderivative{f^{s}}{x}\right)
\qquad\Leftrightarrow\qquad
\textbf{w}^{T}\left(\pderivative{\textbf{f}}{\xi}\right)=\left(\pderivative{f^{s}}{\xi}\right).   
\end{equation}
Thus, the chain rule provides 
\begin{align}\label{Temporalpart1D}
\begin{split}
\mathbf{w}^{T}\left(\frac{d J\mathbf{u}}{dt}\right)
=& \quad J\mathbf{w}^{T}\left(\frac{d\mathbf{u}}{dt}\right)+\left(\frac{dJ}{dt}\right)\mathbf{w}^{T}\mathbf{u} \\ 
=& \quad \left(\frac{d Js}{dt}\right)+\left(\frac{dJ}{dt}\right)\left(\mathbf{w}^{T}\mathbf{u}-s\right) \\
=& \quad \left(\frac{d Js}{dt}\right)+\left(\pderivative{\nu}{\xi}\right)\left(\mathbf{w}^{T}\mathbf{u}-s\right),
\end{split}
\end{align}
where we used the GCL \eqref{GCL1D} in the last step. Likewise, the chain rule and equation \eqref{eq:1DEntropyFlux} give    
\begin{align}\label{Spatialpart1D}
\begin{split}
\mathbf{w}^{T}\left(\pderivative{\mathbf{g}}{\xi}\right) 
=& \quad \mathbf{w}^{T}\left(\pderivative{\mathbf{f}}{\xi}\right)-\nu\mathbf{w}^{T}\left(\pderivative{\mathbf{u}}{\xi}\right)-\left(\pderivative{\nu}{\xi}\right)\mathbf{w}^{T}\mathbf{u} \\
=& \quad \left(\pderivative{f^{s}}{\xi}\right)-\nu\left(\pderivative{s}{\xi}\right)-\left(\pderivative{\nu}{\xi}\right)\mathbf{w}^{T}\mathbf{u} \\
=& \quad \pderivative{}{\xi}\left(f^{s}-\nu s\right)-\left(\pderivative{\nu}{\xi}\right)\left(\mathbf{w}^{T}\mathbf{u}-s\right). 
\end{split}
\end{align} 
For smooth solutions of the conservation law \eqref{eq:consLaw1D}, the equations \eqref{Temporalpart1D},\eqref{Spatialpart1D}, provide the entropy balance law on the reference element   
\begin{equation}\label{1D:EntropyConditionREF}
0=\mathbf{w}^{T}\left(\frac{dJ\mathbf{u}}{dt}\right)+\mathbf{w}^{T}\left(\pderivative{\mathbf{g}}{\xi}\right)=\left(\frac{dJs}{dt}\right)+\pderivative{}{\xi}\left(f^{s}-\nu s\right).
\end{equation}
Next, we integrate the equation \eqref{1D:EntropyConditionREF} over the reference element $\left[-1,1\right]$, apply the fundamental theorem of calculus and rearrange terms. This provides the identity    
\begin{equation}\label{EQ:TotalEntropy1D}
\frac{d}{dt}\int_{-1}^{1}Js\,d\xi=-\left.\left(f^{s}-\nu s\right)\right|_{-1}^{\ 1}.
\end{equation}
The equation \eqref{EQ:TotalEntropy1D} describes the temporal evolution of the total entropy and holds merely for smooth solutions of the conservation law \eqref{eq:consLaw1D}. Nevertheless, for  discontinuous solutions, we obtain the inequality   
\begin{equation}\label{IQ:TotalEntropy1D}
\frac{d}{dt}\int_{-1}^{1}Js\,d\xi\leq-\left.\left(f^{s}-\nu s\right)\right|_{-1}^{\ 1}.
\end{equation}
A numerical method, which satisfies a discrete analogue of the equation \eqref{EQ:TotalEntropy1D}, is called entropy conservative and a method, which satisfies a discrete analogue of the equation \eqref{IQ:TotalEntropy1D}, is called entropy stable approximation. 

On the discrete level, it is not easy to mimic the chain rule. We focus on SBP operators that allow to recover the integration-by-parts formula on the discrete level \cite{Fernandez2014,Gassner2013,kreiss1}. Thus, as an alternative, we re-write \eqref{Spatialpart1D} in such a way that the product rule can be used:  
\begin{equation}\label{1D:EntropyConditionREF1}
\left(\pderivative{\mathbf{w}}{\xi}\right)^{T}\mathbf{g}
=\left(\pderivative{\psi}{\xi}\right)-\nu\left(\pderivative{\phi}{\xi}\right),
\end{equation}
where $\phi:=\mathbf{w}^{T}\mathbf{u}-s$ and $\psi:=\mathbf{w}^{T}\mathbf{f}-f^{s}$ are the entropy potential and the entropy flux potential functions \cite{Tadmor2003}. 

\subsection{Discrete entropy conservation for moving mesh FV methods}\label{Sec:EntropyMovingMeshFV}
Next, we consider a FV discretization of the problem \eqref{eq:consLaw1D}. We divide the domain $\left[a\left(t\right),b\left(t\right)\right]$ in $K$ non-overlapping time-dependent elements $I_{k}\left(t\right):=\left[x_{k}\left(t\right),x_{k+1}\left(t\right)\right]$, $k=1,\dots,K$. Each interval $I_{k}\left(t\right)$ is mapped onto the reference Interval $\left[-1,1\right]$ by an affine linear mapping. The approximation for the element mean value of the exact solution $\textbf{u}$ in the element $I_{k}\left(t\right)$ is denoted by $\textbf{u}_{k}\left(t\right)$ and the approximation for $J$ in the element $I_{k}\left(t\right)$ is denoted by $\mathrm{J}_{k}$. The GCL \eqref{GCL1D} needs to be satisfied on the discrete level. Thus, the equation \eqref{GCL1D} is also discretized. Then, the semi-discrete moving mesh FV scheme is given by     
\begin{subequations}\label{FVScheme}
\begin{align}
\frac{d\mathrm{J}_{k}}{dt}=&\quad \left(\nu_{k+\frac{1}{2}}^{*}-\nu_{k-\frac{1}{2}}^{*}\right), \label{FVGCL} \\
\frac{d\left(\mathrm{J}_{k}\textbf{u}_{k}\right)}{dt}=&-\left(\mathbf{g}_{k+\frac{1}{2}}^{*}-\mathbf{g}_{k-\frac{1}{2}}^{*}\right). \label{FVCL} 
\end{align}
\end{subequations} 
The numerical flux functions $\nu_{k+\frac{1}{2}}^{*}$ and $\mathbf{g}_{k+\frac{1}{2}}^{*}$ in \eqref{FVScheme} satisfy    
\begin{equation}
\nu_{k+\frac{1}{2}}^{*}:=\frac{1}{2}\left(\nu_{k+\frac{1}{2}}^{L}+\nu_{k+\frac{1}{2}}^{R}\right),\qquad 
\text{and} \qquad 
\mathbf{g}_{k+\frac{1}{2}}^{*}:=\mathbf{g}^{*}\left(\nu_{k+\frac{1}{2}}^{L},\nu_{k+\frac{1}{2}}^{R},\mathbf{u}_{k},\mathbf{u}_{k+1}\right), 
\end{equation}  
where $\nu_{k+\frac{1}{2}}^{L}$ is the approximation of the grid velocity from the left side in the cell interface point $x_{k+1}$ and  $\nu_{k+\frac{1}{2}}^{R}$ is the approximation from the right side. Moreover, the flux function $\mathbf{g}^{*}$ needs to satisfy the consistence condition    
\begin{equation}\label{ConsistenceCondition}
\mathbf{g}^{*}\left(\nu^{L},\nu^{R},\mathbf{u},\mathbf{u}\right)=\mathbf{f}\left(\mathbf{u}\right)-\frac{1}{2}\left(\nu^{L}+\nu^{R}\right)\mathbf{u}.
\end{equation}
In the FV context, a discrete analogue of the equation \eqref{1D:EntropyConditionREF1} is given by   
\begin{equation}\label{FV:EntropyCondition}
\left(\mathbf{w}_{k+1}-\mathbf{w}_{k}\right)^{T}\mathbf{g}_{k+\frac{1}{2}}^{*}=\left(\psi_{k+1}-\psi_{k}\right)-\nu_{k+\frac{1}{2}}^{*}\left(\phi_{k+1}-\phi_{k}\right), 
\end{equation} 
\begin{equation}\label{FV:EntropyCondition1}
\phi_{k}:=\mathbf{w}_{k}^{T}\mathbf{u}_{k}-s\left(\mathbf{u}_{k}\right),\qquad\psi_{k}:=\mathbf{w}_{k}^{T}\mathbf{f}\left(\mathbf{u}_{k}\right)-f^{s}\left(\mathbf{u}_{k}\right),
\end{equation}
for all $k=1,\dots,K-1$. We note that for $\nu=0$ (static mesh special case) the equation \eqref{FV:EntropyCondition} becomes Tadmor's discrete entropy  condition for entropy conservative numerical flux functions \cite{Tadmor1987}. 
The following property of the FV method \eqref{FVScheme} is proven in Appendix \ref{FV:EC}.  
\begin{thm-hand}[2.1.]
Suppose the flux function $\mathbf{g}^{*}$ satisfies the condition \eqref{FV:EntropyCondition}.  
Then the semi-discrete moving mesh FV \eqref{FVScheme} scheme is entropy consistent such that 
\begin{equation}\label{FV:TotalEntropyConservation}
\frac{d}{dt}\sum_{k=1}^{K}\mathrm{J}_{k}s\left(\mathbf{u}_{k}\right)=-\left(\mathbf{w}_{K}^{T}\mathbf{g}_{K+\frac{1}{2}}^{*}-\psi_{K}+\nu_{K+\frac{1}{2}}^{*}\psi_{K}\right)+\left(\mathbf{w}_{1}^{T}\mathbf{g}_{\frac{1}{2}}^{*}-\psi_{1}+\nu_{\frac{1}{2}}^{*}\psi_{1}\right).
\end{equation}
\end{thm-hand}
The equation \eqref{FV:TotalEntropyConservation} is a discrete analogue of the identity \eqref{EQ:TotalEntropy1D}. Therefore, in the context of the ALE framework, the condition \eqref{FV:EntropyCondition} can be seen as a natural extension  of Tadmor's discrete entropy condition. Finally, we remark that a flux with the property \eqref{FV:EntropyCondition} can be constructed from an entropy consistent flux for FV methods on static meshes. 
\begin{rem-hand}[2.2.]
If a two-point flux function $\mathbf{f}^{*}$ for static meshes with the Tadmor condition
\begin{equation}\label{FV:EntropyConditionTadmor}
\left(\mathbf{w}_{k+1}-\mathbf{w}_{k}\right)^{T}\mathbf{f}_{k+\frac{1}{2}}^{*}=\left(\psi_{k+1}-\psi_{k}\right),\qquad k=1,\dots,K-1,
\end{equation}  
and a state function $\mathbf{u}^{\#}$ with the property 
\begin{equation}\label{FV:EntropyConditionState}
\left(\mathbf{w}_{k+1}-\mathbf{w}_{k}\right)^{T}\mathbf{u}_{k+\frac{1}{2}}^{\#}=\left(\phi_{k+1}-\phi_{k}\right),\qquad k=1,\dots,K-1,
\end{equation}
are available, a flux function $\mathbf{g}^{*}$ with the entropy condition \eqref{FV:EntropyCondition} can be constructed by 
\begin{equation}\label{FV_EC_Flux}
\mathbf{g}_{k+\frac{1}{2}}^{*}=\mathbf{f}_{k+\frac{1}{2}}^{*}-\nu_{k+\frac{1}{2}}^{*}\mathbf{u}_{k+\frac{1}{2}}^{\#}.
\end{equation}
State functions $\mathbf{u}^{\#}$ with the property \eqref{FV:EntropyConditionState} have been constructed by Friedrich et al. \cite{Friedrich2018} for the shallow water, compressible Euler and ideal magnetohydrodynamic (MHD) equations and are needed for the construction of fully discrete entropy stable space-time DGSEM methods. 
\end{rem-hand}

\subsection{Free stream preservation for moving mesh FV methods} \label{sec:FreestreamPreservationFV}
The geometric and metric terms are time-dependent in a moving mesh approximation. Thus, it is a priori not clear, if the common used time discretization strategies like Runge-Kutta (RK) methods, e.g \cite{Butcher1987}, provide a meaningful approximation for these terms. In particular, it is not clear, if the semi-discrete GCL \eqref{FVGCL} stays true after the time discretization.       A simple test to check, if the discretization of the geometric and metric terms is reasonable, is the free stream preservation test. The capability of the numerical scheme to preserve solutions with constant states is investigated in this test. 

In the following, the time discretization of the semi-discrete moving mesh FV scheme \eqref{FVScheme} is analyzed. We divide the interval $\left[0,T\right]$ in discrete time levels $t^{n}$. The step size of the time discretization is $\Delta t$. The FV solutions $\textbf{u}_{k}$, $\mathrm{J}_{k}$ and the numerical fluxes $\nu_{k+\frac{1}{2}}^{*}$, $\mathbf{g}_{k+\frac{1}{2}}^{*}$ are approximated in the time levels $t^{n}$, e.g. $\mathbf{u}_{k}\left(t^{n}\right)\approx\mathbf{u}_{k}^{n}$. Thus, the moving mesh forward Euler FV discretization is given by  
\begin{subequations}\label{FVScheme1}
\begin{align}
\mathrm{J}_{k}^{n+1}=&\quad\mathrm{J}_{k}^{n}+\Delta t\left(\nu_{k+\frac{1}{2}}^{*,n}-\nu_{k-\frac{1}{2}}^{*,n}\right), \label{Time:FVGCL} \\
\left(\mathrm{J}_{k}\mathbf{u}_{k}\right)^{n+1}=&\quad 
 \mathrm{J}_{k}^{n}\mathbf{u}_{k}^{n}-\Delta t\left(\mathbf{g}_{k+\frac{1}{2}}^{*,n}-\mathbf{g}_{k-\frac{1}{2}}^{*,n}\right), \label{Time:FVCL} \\
\mathbf{u}_{k}^{n+1}=&\quad\left(\mathrm{J}_{k}\mathbf{u}_{k}\right)^{n+1}/\mathrm{J}_{k}^{n+1}. \label{Time:FVUpdate}
\end{align}
\end{subequations} 
The last step is necessary to decouple the solution $\mathbf{u}_{k}$ from the metric term $\mathrm{J}_{k}$. However, it is also possible to write \eqref{FVScheme1} in an equivalent form without the step \eqref{Time:FVUpdate}. Therefore, we apply the D-GCL \eqref{Time:FVGCL} and realize   
\begin{equation}\label{Time:FVSpecialEquation1}
\mathrm{J}_{k}^{n}\mathbf{u}_{k}^{n}=\left(\mathrm{J}_{k}^{n+1}-\Delta t\left(\nu_{k+\frac{1}{2}}^{*,n}-\nu_{k-\frac{1}{2}}^{*,n}\right)\right)\mathbf{u}_{k}^{n}. 
\end{equation}
Then, we plug \eqref{Time:FVSpecialEquation1} in \eqref{Time:FVCL} and obtain 
\begin{subequations}\label{FVScheme2}
\begin{align}
\mathrm{J}_{k}^{n+1}=&\quad\mathrm{J}_{k}^{n}+\Delta t\left(\nu_{k+\frac{1}{2}}^{*,n}-\nu_{k-\frac{1}{2}}^{*,n}\right), \label{Time:FVGCL1} \\
\mathbf{u}_{k}^{n+1}=&\quad\mathbf{u}_{k}^{n}-\frac{\Delta t}{\mathrm{J}_{k}^{n+1}}\left[\left(\mathbf{g}_{k+\frac{1}{2}}^{*,n}+\nu_{k+\frac{1}{2}}^{*,n}\mathbf{u}_{k}^{n}\right)-\left(\mathbf{g}_{k-\frac{1}{2}}^{*,n}+\nu_{k-\frac{1}{2}}^{*,n}\mathbf{u}_{k}^{n}\right)\right]. \label{Time:FVCL1} 
\end{align}
\end{subequations} 
We note that the equation \eqref{Time:FVCL1} matches the structure of the one dimensional moving mesh FV method from Fazio and LeVeque \cite{Fazio2003}. Next, we analyze the forward Euler moving mesh FV scheme with respect to the free stream preservation property. For a solution $\mathbf{u}_{k}^{n}=\mathbf{c}=\left[c_{1},\dots,c_{p}\right]^T$ with constant components $c_{i}\in\R^{p}$, $i=1,\dots,p$, it follows by \eqref{ConsistenceCondition}
\begin{align}\label{FreestreamPreservationFV1}
\begin{split}
\mathbf{u}_{k}^{n+1}  
=& \quad 
\mathbf{c}-\frac{\Delta t}{\mathrm{J}_{k}^{n+1}}\left[\left(\mathbf{g}_{k+\frac{1}{2}}^{*,n}+\nu_{k+\frac{1}{2}}^{*,n}\mathbf{c}\right)-\left(\mathbf{g}_{k-\frac{1}{2}}^{*,n}+\nu_{k+\frac{1}{2}}^{*,n}\mathbf{c}\right)\right] \\
=& \quad \mathbf{c}-\frac{\Delta t}{\mathrm{J}_{k}^{n+1}}\left[\left(\left(\mathbf{f}\left(\mathbf{c}\right)-\nu_{k+\frac{1}{2}}^{*,n}\mathbf{c}\right)+\nu_{k+\frac{1}{2}}^{*,n}\mathbf{c}\right)-\left(\left(\mathbf{f}\left(\mathbf{c}\right)-\nu_{k-\frac{1}{2}}^{*,n}\mathbf{c}\right)+\nu_{k-\frac{1}{2}}^{*,n}\mathbf{c}\right)\right]=\mathbf{c}. 
\end{split}
\end{align}
Thus, the forward Euler moving mesh FV method satisfies the free stream preservation property. 

In the same way, the semi-discrete moving mesh FV scheme \eqref{FVScheme} can be discretized by an explicit $s$-stage RK method with the characteristic coefficients $\left\{ a_{\tau \sigma}\right\} _{\tau,\sigma=1}^{s}$, $\left\{ b_{\sigma}\right\} _{\sigma=1}^{s}$, $\left\{ c_{\sigma}\right\} _{\sigma=1}^{s}$. This provides, the fully-discrete moving mesh RK-FV method: 
\begin{subequations}\label{RK_FVScheme}
\begin{align}
& \text{for $\tau=1,\dots,s$:} \nonumber \\
& \mathrm{J}_{k}^{\left(\tau\right)}=\mathrm{J}_{k}^{n}+\Delta t\sum_{\sigma=1}^{\tau-1}a_{\tau \sigma }\left(\nu_{k+\frac{1}{2}}^{*,n+\sigma}-\nu_{k-\frac{1}{2}}^{*,n+\sigma}\right), \label{RK_FVScheme:GCL} \\
& \mathbf{u}_{k}^{\left(\tau\right)}=\mathbf{u}_{k}^{n}-\frac{\Delta t}{\mathrm{J}_{k}^{\left(\tau\right)}}\sum_{\sigma=1}^{\tau-1}a_{\tau \sigma}\left[\left(\mathbf{g}_{k+\frac{1}{2}}^{*,n+\sigma}+\nu_{k+\frac{1}{2}}^{*,n+\sigma}\mathbf{u}_{k}^{\left(\sigma\right)}\right)-\left(\mathbf{g}_{k-\frac{1}{2}}^{*,n+\sigma}+\nu_{k-\frac{1}{2}}^{*,n+\sigma}\mathbf{u}_{k}^{\left(\sigma\right)}\right)\right], \label{RK_FVScheme:CL} \\
& \nonumber \\ 
& \mathrm{J}_{k}^{n+1}=\mathrm{J}_{k}^{n}+\Delta t\sum_{j=1}^{s}b_{\sigma}\left(\nu_{k+\frac{1}{2}}^{*,n+\sigma}-\nu_{k-\frac{1}{2}}^{*,n+\sigma}\right), \label{RK_FVScheme:GCL1} \\
& \mathbf{u}_{k}^{n+1}=\mathbf{u}_{k}^{n}-\frac{\Delta t}{\mathrm{J}_{k}^{n+1}}\sum_{\sigma=1}^{s}b_{\sigma}\left[\left(\mathbf{g}_{k+\frac{1}{2}}^{*,n+\sigma}+\nu_{k+\frac{1}{2}}^{*,n+\sigma}\mathbf{u}_{k}^{\left(\sigma\right)}\right)-\left(\mathbf{g}_{k-\frac{1}{2}}^{*,n+\sigma}+\nu_{k-\frac{1}{2}}^{*,n+\sigma}\mathbf{u}_{k}^{\left(\sigma\right)}\right)\right], \label{RK_FVScheme:CL1} 
\end{align}
\end{subequations} 
where  $\nu_{k+\frac{1}{2}}^{L,n+\sigma}:=\nu_{k+\frac{1}{2}}^{L}\left(t^{n}+c_{\sigma}\Delta t\right)$, $\nu_{k+\frac{1}{2}}^{R,n+\sigma}:=\nu_{k+\frac{1}{2}}^{R}\left(t^{n}+c_{\sigma}\Delta t\right)$ and
\begin{align}
\nu_{k+\frac{1}{2}}^{*,n+\sigma}
=&\quad\frac{1}{2}\left(\nu_{k+\frac{1}{2}}^{L,n+\sigma}+\nu_{k+\frac{1}{2}}^{R,n+\sigma}\right), \\
\mathbf{g}_{k+\frac{1}{2}}^{*,n+\sigma}
=&\quad\mathbf{g}^{*}\left(\nu_{k+\frac{1}{2}}^{L,n+\sigma},\nu_{k+\frac{1}{2}}^{R,n+\sigma},\mathbf{u}_{k}^{\left(\sigma\right)},\mathbf{u}_{k}^{\left(\sigma\right)}\right).
\end{align}

The free stream preservation property for the moving mesh RK-FV method \eqref{RK_FVScheme} can be proven by the same calculation as in \eqref{FreestreamPreservationFV1}. In particular, the following statement can be proven by exactly the same arguments as for the forward Euler step. 
\begin{thm-hand}[2.3.]
Let $\mathbf{u}_{k}^{n}=\mathbf{c}:=\left(c_{1},\dots,c_{p}\right)^{T}\in\R^{p}$ be the solution of the fully-discrete FV \eqref{RK_FVScheme} scheme at time level $t^{n}$. Furthermore, the numerical flux $\mathbf{g}^{*}$ satisfies  \eqref{ConsistenceCondition}. Then, the constant states $c_{i}$, $i=1,\dots,p$, are preserved in each Runge-Kutta stage \eqref{RK_FVScheme:CL}. In particular, the solution of the fully-discrete FV method at time level $t^{n+1}$ is $\mathbf{u}_{k}^{n+1}=\mathbf{c}$.  
\end{thm-hand}

\section{Entropy stable DGSEM on moving meshes}\label{Sec:EntropyStableDGSEM}
In Section \ref{Sec:BasicConcepts} an entropy conservative first order moving mesh FV method has been introduced. The main goal of this work is the construction of an entropy stable high order moving mesh DGSEM. On static meshes, it is possible to construct high-order entropy stable DGSEM, if the derivative matrix is an SBP operator and entropy conservative two-point flux functions are available. This methodology has been used in the construction of high-order entropy stable DGSEM on quadrilateral/hexahedral elements, e.g.  \cite{Bohm2018,Gassner2016,Wintermeyer2017}. In this section, it will be shown that similar ideas can be used to construct high-order entropy stable moving mesh DGSEM, when there are two-point flux functions with the property \eqref{FV:EntropyCondition} available for the low-order FV method \eqref{FVScheme}. The construction of the entropy stable moving mesh DGSEM will be presented for an arbitrary symmetrizable and hyperbolic system of conservation laws  
\begin{equation}\label{eq:consLaw}
\pderivative{\textbf{u}}{t}+\sum_{i=1}^{3}\pderivative{\textbf{f}_{i}}{x_{i}}=\textbf{0},
\end{equation}
on a time-dependent domain $\Omega\left(t\right)\subseteq \R^3$. The vector of conserved variables is $\textbf{u}$ and $\textbf{f}_i,\,i=1,2,3$, are the physical flux vectors. The state vectors are of size $p$ depending on the number of equations in the system under consideration and the conservation law is subjected to appropriate initial and boundary conditions.

The block vector nomenclature in \cite{Gassner2017} simplifies the analysis of the system \eqref{eq:consLaw} on curved elements. Thus, we translate the conservation law \eqref{eq:consLaw} in block vector notation. A block vector is highlighted by the double arrow    
\begin{equation}\label{BlockFlux}
\blockvec{\textbf{f}}:=\begin{bmatrix}
\mathbf{f}_{1} \\
\mathbf{f}_{2} \\
\mathbf{f}_{3}
\end{bmatrix}.
\end{equation}
The dot product of two block vectors is given by   
\begin{equation}\label{ProductBlock}
\blockvec{\textbf{f}}\cdot\blockvec{\textbf{g}}:=\sum_{i=1}^{3}\textbf{f}_{i}^{T}\textbf{g}_{i}.
\end{equation}
Furthermore, the dot product of a vector $\vec{v}$ in the three dimensional space and a block vector is defined by    
\begin{equation}\label{ProductVecBlock}
\vec{v}\cdot\blockvec{\textbf{f}}:=\sum_{i=1}^{3}v_{i}\textbf{f}_{i}.
\end{equation}
We note that the dot product \eqref{ProductBlock} is a scalar quantity and the dot product \eqref{ProductVecBlock} is a vector in a $p$ dimensional space, where the number $p$ corresponds to the number of conserved variables in the conservation law \eqref{eq:consLaw}. The interaction between a vector $\vec{v}$ and the conserved variables is defined as the block vector   
\begin{equation}\label{BlockVector}
\vec{v}\,\textbf{u}:=\begin{bmatrix}v_{1}\mathbf{u}\\
v_{2}\mathbf{u}\\
v_{3}\mathbf{u}
\end{bmatrix}.
\end{equation}
Thus, in particular, the spatial gradient of the conserved variables is defined by     
\begin{equation}\label{BlockGradient}
\vec{\nabla}_{x}\mathbf{u}:=\begin{bmatrix}\pderivative{\mathbf{u}}{x_{1}}\\
\pderivative{\mathbf{u}}{x_{2}}\\
\pderivative{\mathbf{u}}{x_{3}}
\end{bmatrix}.
\end{equation}
The gradient of a vector valued function $\vec{g}=[g_{1},g_{2},g_{3}]^T$ is a second order tensor, written in matrix form as
\begin{equation}\label{JacobianMatrix}
\vec{\nabla}_{x}\vec{g}=\begin{bmatrix}\pderivative{g_{1}}{x_{1}} & \pderivative{g_{1}}{x_{2}} & \pderivative{g_{1}}{x_{3}}\\
\pderivative{g_{2}}{x_{1}} & \pderivative{g_{2}}{x_{2}} & \pderivative{g_{2}}{x_{3}}\\
\pderivative{g_{3}}{x_{1}} & \pderivative{g_{3}}{x_{2}} & \pderivative{g_{3}}{x_{3}}
\end{bmatrix}.
\end{equation}
The dot product \eqref{ProductBlock} and the spatial gradient \eqref{BlockGradient} are used to define the divergence of a block vector flux as  
\begin{equation}\label{BlockDivergence}
\vec{\nabla}_{x}\cdot\blockvec{\textbf{f}}:=\sum_{i=1}^{3}\pderivative{\textbf{f}_{i}}{x_{i}}.
\end{equation}
Moreover, for a vector valued function $\vec{g}$ and the conserved variables,   we have the product rule     
\begin{equation}\label{ProductRule}
\vec{\nabla}_{x}\cdot\left(\vec{g}\,\mathbf{u}\right)=\left(\vec{\nabla}_{x}\cdot\vec{g}\right)\mathbf{u}+\vec{g}\cdot\left(\vec{\nabla}_{x}\mathbf{u}\right)
\end{equation}
with respect to the dot products \eqref{ProductBlock} and \eqref{ProductVecBlock}. These notations allow to write the conservation law \eqref{eq:consLaw} in the compact form    
\begin{equation}\label{eq:consLawBlock}
\pderivative{\textbf{u}}{t} + \vec{\nabla}_x\cdot\blockvec{\textbf{f}} = \textbf{0}.
\end{equation}
 
\subsection{Building blocks of the ALE transformation for hexahedral curved meshes}\label{Sec:ComponentsALE2D}
In order to set up the moving mesh DGSEM in the Section \ref{sec:SemiDiscreteDG}, we make for all $t\in\left[0,T\right]$ the assumptions:  
\begin{itemize}
\item[(A1)] For a fixed number $K\in\N$ the physical domain $\Omega\left(t\right)$ can be subdivided into $K$ time-dependent, non-overlapping and conforming hexahedral elements, $e_{\kappa}(t)$, $\kappa=1,\dots,K$. These elements can have curved faces. 
\item[(A2)] The time-dependent elements $e_{\kappa}(t)$ are mapped into the spatial computational domain $E=\left[-1,1\right]^{3}$ with a bijective isoparametric transfinite mapping. Winters constructed in his PHD thesis \cite{Winters2014_A} a mapping for this set up. Like in \cite{Winters2014_A}, it is assumed that the curved faces satisfy for all $t\in \left[0,T\right]$     
\begin{align}\label{BoundaryPoints}
\begin{tabular}{lll}
$\vec{\Gamma}_{1}\left(-1,\xi^{3},t\right)=\vec{\Gamma}_{6}\left(-1,\xi^{3},t\right)$, & $\vec{\Gamma}_{2}\left(-1,\xi^{3},t\right)=\vec{\Gamma}_{6}\left(1,\xi^{3},t\right)$, & $\vec{\Gamma}_{3}\left(-1,\xi^{2},t\right)=\vec{\Gamma}_{6}\left(\xi^{2},-1,t\right)$,\tabularnewline
$\vec{\Gamma}_{1}\left(1,\xi^{3},t\right)=\vec{\Gamma}_{4}\left(-1,\xi^{3},t\right)$, & $\vec{\Gamma}_{2}\left(1,\xi^{3},t\right)=\vec{\Gamma}_{4}\left(1,\xi^{3},t\right)$, & $\vec{\Gamma}_{3}\left(1,\xi^{2},t\right)=\vec{\Gamma}_{4}\left(\xi^{2},-1,t\right)$, \tabularnewline
$\vec{\Gamma}_{1}\left(\xi^{1},-1,t\right)=\vec{\Gamma}_{3}\left(\xi^{1},-1,t\right)$, & $\vec{\Gamma}_{2}\left(\xi^{1},-1,t\right)=\vec{\Gamma}_{3}\left(\xi^{1},1,t\right)$, & $\vec{\Gamma}_{5}\left(-1,\xi^{2},t\right)=\vec{\Gamma}_{6}\left(\xi^{2},1,t\right)$, \tabularnewline
$\vec{\Gamma}_{1}\left(\xi^{1},1,t\right)=\vec{\Gamma}_{5}\left(\xi^{1},-1,t\right)$, & $\vec{\Gamma}_{2}\left(\xi^{1},1,t\right)=\vec{\Gamma}_{5}\left(\xi^{1},1,t\right)$, & 
$\vec{\Gamma}_{5}\left(1,\xi^{2},t\right)=\vec{\Gamma}_{4}\left(\xi^{2},1,t\right)$.
\tabularnewline 
\end{tabular}
\end{align}
The location of the curved faces is sketched in Figure \ref{ReferencePhysicalCell3D}. 
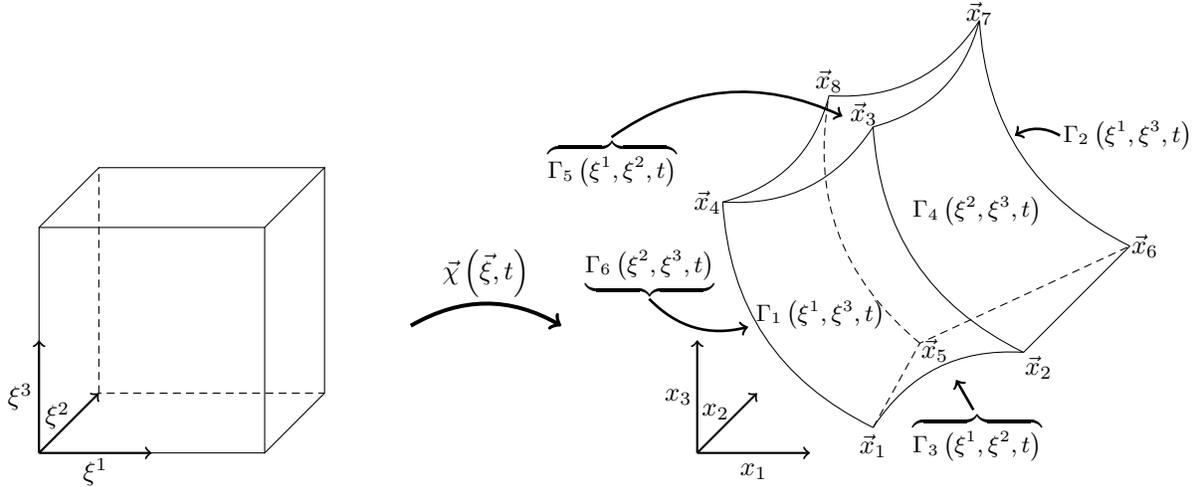
\begin{figure}[h]
\begin{center}
\begin{tikzpicture}[%
x={(1cm,0cm)},
y={(0cm,1cm)},
z={({0.5*cos(45)},{0.5*sin(45)})},scale=1]

\coordinate (A1) at (3,1,1); 
\coordinate (B1) at (5,2,1) ;
\coordinate (C1) at (3,5,1); 
\coordinate (D1) at (1,4,1); 
\coordinate (E1) at (1.5,0,7); 
\coordinate (F1) at (5,2,5); 
\coordinate (G1) at (3,5,5); 
\coordinate (H1) at (1,4,5);

\node at (3,0.75,1) {$\vec{x}_{1}$};
\node at (5.2,1.8,1) {$\vec{x}_{2}$};
\node at (2.7,5,1.5) {$\vec{x}_{3}$};
\node at (0.8,4,1) {$\vec{x}_{4}$};
\node at (1.7,-0.1,7) {$\vec{x}_{5}$};
\node at (5.2,2,5) {$\vec{x}_{6}$};
\node at (3,5.1,5) {$\vec{x}_{7}$};
\node at (1,4.2,5) {$\vec{x}_{8}$};

\node at (2.3,2.5,1)  {\small{$\Gamma_{1}\left(\xi^{1},\xi^{3},t\right)$}};
\node at (5.5,4,3.5)  {\small{$\Gamma_{2}\left(\xi^{1},\xi^{3},t\right)$}}; 
\draw[line width=1pt,bend right,->] (4.6,4,3.5) to (4,4,3.5);
\node at (3.5,0,3.5)  {\small{$\overbrace{\Gamma_{3}\left(\xi^{1},\xi^{2},t\right)}$}};
\draw[line width=1pt,,->] (3.45,0.35,3.5) to (3.22,0.75,3.5);
\node at (-1,4,2.5)  {\small{$\overbrace{\Gamma_{5}\left(\xi^{1},\xi^{2},t\right)}$}};
\draw[line width=1pt,bend left,->] (-0.69,4.65,1.6) to (2.15,4.7,2.25);
\node at (3.5,3,3.5)  {\small{$\Gamma_{4}\left(\xi^{2},\xi^{3},t\right)$}};
\node at (-0.5,2.5,2.5)  {\small{$\underbrace{\Gamma_{6}\left(\xi^{2},\xi^{3},t\right)}$}};
\draw[line width=1pt,bend right,->] (-0.19,2.5,1.6) to (1.25,2.25,1.25);

\draw[bend left,-] (A1) to   (B1);
\draw[bend left,-] (B1) to (C1); 
\draw[bend left,-]  (C1) to  (D1); 
\draw[bend right,-] (D1) to (A1); 
\draw[-] (B1) to (F1);
\draw[bend left,-] (F1) to (G1); 
\draw[bend left,-] (G1) to (C1);
\draw[bend left,-] (G1) to (H1); 
\draw[bend left,-] (H1) to (D1);
\draw[densely dashed] (A1) -- (E1) -- (F1);
\draw[bend left,densely dashed] (E1) to (H1);

\draw[thick,->] (0.75,0.75,0.75) -- (2.25,0.75,0.75); 
\node at (1.5,0.5,0.75)  {$x_{1}$};
\draw[thick,->] (0.75,0.75,0.75) -- (0.75,0.75,3); 
\node at (0.5,0.8,2.15)  {$x_{2}$};
\draw[thick,->] (0.75,0.75,0.75) -- (0.75,2.25,0.75);
\node at (0.5,1.5,0.75)  {$x_{3}$};

\coordinate (A) at (-8,0.75,0.75); 
\coordinate (B) at (-5,0.75,0.75) ;
\coordinate (C) at (-5,3.75,0.75); 
\coordinate (D) at (-8,3.75,0.75); 
\coordinate (E) at (-8,0.75,3); 
\coordinate (F) at (-5,0.75,3); 
\coordinate (G) at (-5,3.75,3); 
\coordinate (H) at (-8,3.75,3);

\draw[] (A) -- (B)  -- (C)  -- (D)  -- (A);
\draw[] (B) to (F);
\draw[] (F) -- (G) -- (C);
\draw[] (G) -- (H) -- (D);
\draw[densely dashed] (A) -- (E) -- (F);
\draw[densely dashed] (E) to (H);

\draw[thick,->] (-8,0.75,0.75) -- (-6.5,0.75,0.75); 
\node at (-7.25,0.5,0.75)  {$\xi^{1}$};
\draw[thick,->] (-8,0.75,0.75) -- (-8,0.75,3); 
\node at (-8.25,0.75,2.15)  {$\xi^{2}$};
\draw[thick,->] (-8,0.75,0.75) -- (-8,2.25,0.75);
\node at (-8.25,1.5,0.75)  {$\xi^{3}$};

\draw[line width=1.5pt,bend left,->] (-3.5, 2,2) to (-1.5, 2,2);
\node at (-2.5,2.7,2)  {$\vec{\chi}\left(\vec{\xi},t\right)$};
\end{tikzpicture}
\caption{\label{ReferencePhysicalCell3D} Left the reference element $E=[-1,1]^3$ and on the right a general hexahedral element $e_{\kappa}(t)$ with the curved faces $\vec{\Gamma}_{1}\left(\xi^{1},\xi^{3},t\right)$, $\vec{\Gamma}_{2}\left(\xi^{1},\xi^{3},t\right)$, $\vec{\Gamma}_{3}\left(\xi^{1},\xi^{2},t\right)$, $\vec{\Gamma}_{4}\left(\xi^{2},\xi^{3},t\right)$, $\vec{\Gamma}_{5}\left(\xi^{1},\xi^{2},t\right)$, and $\vec{\Gamma}_{6}\left(\xi^{2},\xi^{3},t\right)$. The mapping $\vec{\chi}\left(\vec{\xi},t\right)$ connects $E$ and $e_{\kappa}(t)$.}
\end{center}
\end{figure}
The curved faces of an element $e_{\kappa}(t)$ are approximated as interpolation polynomials up to degree $N$ such that     
\begin{equation}
\interpolation{N}{\left(\vec{\Gamma}_{i}\right)}\left(\eta,\zeta,t\right):=\sum_{j,k=0}^{N}\vec{\Gamma}_{i}\left(\eta_{j},\zeta_{k},t\right)\ell_{j}\left(\eta\right)\ell_{k}\left(\zeta\right),\quad i=1,2,3,4,5,6,  
\end{equation}
where $\left\{ \ell_{j}\right\} _{j=0}^{N}$, $\left\{ \ell_{k}\right\} _{k=0}^{N}$ are the Lagrange polynomials associated with the interpolation points $\left\{ \eta_{j}\right\} _{j=0}^{N}$ and $\left\{ \zeta_{k}\right\} _{k=0}^{N}$.  
\item[(A3)] The determinant $J$ of the Jacobian matrix $\vec{\nabla}_{\vec{\xi}}\vec{\chi}$ satisfies     
\begin{equation}\label{DeterminantJabobian}
J:=\text{det}\left(\vec{\nabla}_{\vec{\xi}}\vec{\chi}\right)>0,\qquad\forall t\in\left[0,T\right].
\end{equation}       
\end{itemize}
Mesh curving techniques are discussed by Hindenlang et al. \cite{Hindenlang2015} and methodologies to construct a moving mesh with the properties (A1)-(A3) are given in the literature e.g. the book of Huang and Russell \cite[Chapter 6, Chapter 7]{Huang2010}. The mapping $\vec{\chi}=\left[\chi_{1},\chi_{2},\chi_{3}\right]^{T}$ provides the grid velocity field 
  \begin{equation}\label{GridVelocity}
\vec{\nu}=\left[\nu_{1},\nu_{2},\nu_{3}\right]^{T}:=\left[\frac{d\chi_{1}}{dt},\frac{d\chi_{2}}{dt},\frac{d\chi_{3}}{dt}\right]^{T}=\frac{d\vec{\chi}}{dt}.
\end{equation}
It is desirable that the grid velocity is continuous, since the mesh should be conforming and watertight at each time level. The next statement provides conditions on the element boundaries to guarantee that the grid velocity becomes continuous. 
\begin{lem-hand}[3.1] 
Let $e_{1}(t)$ and $e_{2}(t)$ be two neighboring elements which share one of the faces 
\begin{equation}
\vec{\Gamma}_{1}^{1}=\vec{\Gamma}_{2}^{2},\quad\vec{\Gamma}_{3}^{1}=\vec{\Gamma}_{5}^{2},\quad\vec{\Gamma}_{4}^{1}=\vec{\Gamma}_{6}^{2}, \quad
\vec{\Gamma}_{1}^{2}=\vec{\Gamma}_{2}^{1},\quad\vec{\Gamma}_{3}^{2}=\vec{\Gamma}_{5}^{1},\quad\vec{\Gamma}_{4}^{2}=\vec{\Gamma}_{6}^{1},
\end{equation}
where $\vec{\Gamma}_{i}^{l}$, $l=1,2$,  and $i=1,2,3,4,5,6$, are the faces of the element $e_{l}(t)$. Furthermore, suppose that the faces $\vec{\Gamma}_{i}^{l}\left(\cdot,\cdot,t\right)$ are continuously differentiable in the time interval $\left[0,T\right]$. Then the grid velocity field is continuous in the 
points which belong to the face that the elements share. 
\end{lem-hand}
In Appendix \ref{sec:ProofGridVelocity} the Lemma 3.1 is proven in two dimensions. The three dimensional proof can be done by the same argumentation.   

\subsection{Transformation of the conservation law onto a reference element}\label{Sec:TransCL}
In the following, we show that the system \eqref{eq:consLawBlock} can be transformed from a time-dependent element $e_{\kappa}(t)$ onto the reference element $E$. The mapping $\vec{\chi}=\left[\chi_{1},\chi_{2},\chi_{3}\right]^{T}$ provides the covariant basis vectors   
\begin{equation}
\vec{a}_{i}:=\pderivative{\vec{\chi}}{\xi^{i}},\qquad i=1,2,3,
\end{equation}
and the volume weighted contravariant vectors 
\begin{equation}\label{ContravariantVectors}
Ja^{i}=a_{j}\times a_{k}, \quad (i, j, k) \ \text{cyclic}.
\end{equation}
Thereby, is $\vec{\xi}=\left(\xi^{1},\xi^{2},\xi^{3}\right)^{T}$ a vector in the reference element $E=[-1,1]^3$. The covariant and the volume weighted contravariant vectors represent the Jacobian matrix $\vec{\nabla}_{\vec{\xi}}\vec{\chi}$ and its adjoint matrix 
\begin{equation}\label{JacobiMatrixAdjoint}
\vec{\nabla}_{\vec{\xi}}\vec{\chi}=\begin{bmatrix}\vec{a}_{1} & \vec{a}_{2} & \vec{a}_{3}
\end{bmatrix},\qquad\text{adj}\left(\vec{\nabla}_{\vec{\xi}}\vec{\chi}\right)=\begin{bmatrix}\left(J\vec{a}^{1}\right)^{T}\\
\left(J\vec{a}^{2}\right)^{T} \\ 
\left(J\vec{a}^{3} \right)^{T}
\end{bmatrix}.
\end{equation}
Furthermore, the contravariant vectors satisfy the metric identities    
\begin{equation}\label{MetricIdentities}
\sum_{i=1}^{3}\pderivative{J\vec{a}^{i}}{\xi^{i}}=0.
\end{equation}
In particular, the covariant and the contravariant vectors allow to transform differential operators on the time-independent reference element $E$. On the reference element the gradient of a function $f$ is given by 
\begin{equation}\label{Trans:Gradient}
\vec{\nabla}_{x}f=\frac{1}{J}\left(\sum_{i=1}^{3}J\vec{a}^{i}\pderivative{f}{\xi^{i}}\right)=\frac{1}{J}\left[\text{adj}\left(\vec{\nabla}_{\vec{\xi}}\vec{\chi}\right)\right]^{T}\,\vec{\nabla}_{\xi}f
\end{equation}
and the divergence of a vector valued function $\vec{g}$ is given by    
\begin{equation}\label{Trans:Divergence}
\vec{\nabla}_{x}\cdot\vec{g}=\frac{1}{J}\sum_{i=1}^{3}\pderivative{}{\xi^{i}}\left(J\vec{a}^{i}\cdot\vec{g}\right)=\frac{1}{J}\vec{\nabla}_{\xi}\cdot\vec{\tilde{g}},
\end{equation}
where we used the contravariant flux 
\begin{equation}\label{ContravariantFunction}
\vec{\tilde{g}}:=\begin{bmatrix}J\vec{a}^{1}\cdot\vec{g}\\
J\vec{a}^{2}\cdot\vec{g}\\
J\vec{a}^{3}\cdot\vec{g}
\end{bmatrix}=\text{adj}\left(\vec{\nabla}_{\vec{\xi}}\vec{\chi}\right)\,\vec{g}.
\end{equation}
In \cite{Gassner2017}, the following block matrix has been introduced to combine the transformations \eqref{Trans:Gradient} and \eqref{Trans:Divergence} with the block vector notation    
\begin{equation}\label{TransformationBlock}
\blockmatx{M}=\left[\text{adj}\left(\vec{\nabla}_{\vec{\xi}}\vec{\chi}\right)\right]^{T}\otimes\matx{I_{\textit{p}}}=\left[\begin{matrix}Ja_{1}^{1}\matx{I_{\textit{p}}} & Ja_{1}^{2}\matx{I_{\textit{p}}} & Ja_{1}^{3}\matx{I_{\textit{p}}}\\[0.1cm]
Ja_{2}^{1}\matx{I_{\textit{p}}} & Ja_{2}^{2}\matx{I_{\textit{p}}} & Ja_{2}^{3}\matx{I_{\textit{p}}}\\[0.1cm]
Ja_{3}^{1}\matx{I_{\textit{p}}} & Ja_{3}^{2}\matx{I_{\textit{p}}} & Ja_{3}^{3}\matx{I_{\textit{p}}}
\end{matrix}\right],
\end{equation}
where $\otimes$ is the Kronecker product, the matrix $\matx{I_{\textit{p}}}$ is the $p\times p$ identity matrix and $Ja_{j}^{i}$ is the component of $J \vec{a}^{i}$ in the $j$-th Cartesian coordinate direction. The transformation of the gradient becomes     
\begin{equation}\label{BlockGradientReference}
\vec{\nabla}_{x}\mathbf{u}=\frac{1}{J}\blockmatx{M}\,\vec{\nabla}_{\xi}\mathbf{u}.
\end{equation}  
We note that for a vector valued function $\vec{g}$ the following identity holds  
\begin{equation}\label{BlockGradientReferenceAndVector}
\vec{g}\cdot\vec{\nabla}_{x}\mathbf{u}=\frac{1}{J}\vec{g}\cdot\blockmatx{M}\,\vec{\nabla}_{\xi}\mathbf{u}=\frac{1}{J}\vec{\tilde{g}}\cdot\vec{\nabla}_{\xi}\mathbf{u}.
\end{equation}
Moreover, by applying the metric identities \eqref{MetricIdentities}, the transformation of the divergence can be written as 
\begin{equation}\label{BlockDivergenceReference}
\vec{\nabla}_{x}\cdot\blockvec{\textbf{f}}=\frac{1}{J}\vec{\nabla}_{\xi}\cdot\blockmatx{M}^{T}\blockvec{\textbf{f}}.
\end{equation}  
Hence, the contravariant block vector flux is given by 
\begin{equation}\label{ContravariantBlockVector}
\blockveccon{\textbf{f}}
:= \begin{bmatrix}J\vec{a}^{1}\cdot\blockvec{\textbf{f}}\\
J\vec{a}^{2}\cdot\blockvec{\textbf{f}} \\
J\vec{a}^{3}\cdot\blockvec{\textbf{f}}
\end{bmatrix}
=\blockmatx{M}^{T}\blockvec{\textbf{f}}.
\end{equation} 
Since the elements $\left\{ e_{k}\left(t\right)\right\}_{k=1}^{K}$ are time-dependent, the time evolution of the quantity $J$ needs to be analyzed. Thus, we apply Jacobi's formula (cf. e.g. Bellman \cite{Bellman}) and obtain by  \eqref{GridVelocity}, \eqref{JacobiMatrixAdjoint} 
\begin{align}\label{GCL1}
\begin{split}
\frac{d J}{dt}= \text{tr}\left[\text{adj}\left(\vec{\nabla}_{\vec{\xi}}\vec{\chi}\right)\frac{d}{dt}\left(\vec{\nabla}_{\vec{\xi}}\vec{\chi}\right)\right]  =\sum_{i=1}^{3}J\vec{a}^{i}\cdot\left(\frac{d\vec{a}_{i}}{dt}\right) 
= \sum_{i=1}^{3}J\vec{a}^{i}\cdot\left(\pderivative{\vec{\nu}}{\xi^{i}}\right),
\end{split}
\end{align}
where $\text{tr}\left[\cdot\right]$ denotes the trace of a matrix. The metric identities \eqref{MetricIdentities} allow to write the equation \eqref{GCL1} in conservation form      
\begin{equation}\label{GCL2}
\frac{dJ}{dt}=\sum_{i=1}^{3}\pderivative{}{\xi_{i}}\left(J\vec{a}^{i}\cdot\vec{\nu}\right)=\vec{\nabla}_{\xi}\cdot\vec{\tilde{\nu}}.
\end{equation}
The equation \eqref{GCL2} is the GCL in three dimensions. The chain rule formula and the identity \eqref{BlockGradientReferenceAndVector} provide 
\begin{equation}\label{TimeDerivative1}
\frac{d\mathbf{u}}{dt} = \pderivative{\mathbf{u}}{t}+\frac{1}{J}\vec{\tilde{\nu}}\cdot\vec{\nabla}_{\xi}\mathbf{u}. 
\end{equation}
Next, we plug the GCL \eqref{GCL2} into equation \eqref{TimeDerivative1},  apply the product rule \eqref{ProductRule} and rearrange. This provides the equation 
\begin{equation}\label{TimeDerivative2}
J\pderivative{\mathbf{u}}{t}=\frac{d\left(J\mathbf{u}\right)}{dt}-\vec{\nabla}_{\xi}\cdot\left(\vec{\tilde{\nu}}\mathbf{u}\right).
\end{equation}
Finally, we combine the identities \eqref{BlockDivergenceReference} and \eqref{TimeDerivative2} to write the the conservation law \eqref{eq:consLawBlock} in the following form      
\begin{equation}\label{eq:consLawBlockRefE}
\frac{d\left(J\mathbf{u}\right)}{dt}+\vec{\nabla}_{\xi}\cdot\blockvec{\tilde{\mathbf{g}}}=\textbf{0},
\end{equation}
where
\begin{equation}\label{G-flux}
\blockvec{\mathbf{g}}=
\begin{bmatrix}\mathbf{g}_{1}\\
\mathbf{g}_{2}\\
\mathbf{g}_{3}\\
\end{bmatrix}
:=\begin{bmatrix}\mathbf{f}_{1}-\nu_{1}\mathbf{u}\\
\mathbf{f}_{2}-\nu_{2}\mathbf{u} \\
\mathbf{f}_{3}-\nu_{3}\mathbf{u}
\end{bmatrix}=\blockvec{\mathbf{f}}-\vec{\nu}\mathbf{u}.
\end{equation}
The formulation \eqref{eq:consLawBlockRefE} is the representation of the system \eqref{eq:consLawBlock} on the time-independent reference element $E$ for a time-dependent element $e_{\kappa}(t)$.   

\subsection{Entropy analysis in three dimensions}\label{EntropyAnalysis}
The system \eqref{eq:consLaw} is assumed to be symmetrizable. Thus, in particular, it is equipped with entropy/entropy flux pairs $\left(s,f_{i}^{s}\right)$, $i=1,2,3$, (cf. Godunov \cite{Godunov1961} and Mock  \cite{Mock1980}). The strictly convex function $s$ is the entropy function. The entropy function $s$ provides the entropy variables       
\begin{equation}\label{EntropyVariables}
\textbf{w}:=\pderivative{s}{\textbf{u}},
\end{equation}
and it follows by the chain rule  
\begin{equation}\label{PropertiesEntropy}
\pderivative{s}{t}=\mathbf{w}^{T}\pderivative{\mathbf{u}}{t},\qquad
\pderivative{s}{x_{i}}=\mathbf{w}^{T}\pderivative{\mathbf{u}}{x_{i}},\qquad i=1,2,3.
\end{equation}
The entropy flux functions and the flux functions in the conservation law are related and satisfy
\begin{equation}\label{PropertiesEntropyFluxFunctions}
\mathbf{w}^{T}\pderivative{\mathbf{f}_{i}}{x_{i}}=\pderivative{f_{i}^{s}}{x_{i}}
\qquad i=1,2,3.
\end{equation}  
The identities \eqref{BlockGradientReferenceAndVector} and \eqref{PropertiesEntropy} give    
\begin{equation}\label{PropertiesEntropyRefE1}
\mathbf{w}^{T}\left(\vec{\tilde{\nu}}\cdot\vec{\nabla}_{\xi}\mathbf{u}\right)=J\mathbf{w}^{T}\left(\vec{\nu}\cdot\vec{\nabla}_{x}\mathbf{u}\right)=J\left(\vec{\nu}\cdot\vec{\nabla}_{x}s\right)=\vec{\tilde{\nu}}\cdot\vec{\nabla}_{\xi}s.
\end{equation}
Hence, we obtain with the identity \eqref{TimeDerivative1} and the chain rule  
\begin{equation}\label{PropertiesEntropyRefE2}
J\mathbf{w}^{T}\frac{d\mathbf{u}}{dt}=J\pderivative{s}{t}+\vec{\tilde{\nu}}\cdot\vec{\nabla}_{\xi}s=J\frac{ds}{dt}.
\end{equation}
Therefore, the product rule provides the identity   
\begin{align}\label{PropertiesEntropyRefE3}
\begin{split}
\mathbf{w}^{T}\frac{d\left(J\mathbf{u}\right)}{dt} 
=& \quad J\frac{ds}{dt}+\left(\frac{dJ}{dt}\right)\mathbf{w}^{T}\mathbf{u} \\
=& \quad  \frac{d\left(Js\right)}{dt}+\left(\frac{dJ}{dt}\right)\left(\mathbf{w}^{T}\mathbf{u}-s\right)\\
=& \quad \frac{d\left(Js\right)}{dt}+\left(\vec{\nabla}_{\xi}\cdot\vec{\tilde{\nu}}\right)\left(\mathbf{w}^{T}\mathbf{u}-s\right),
\end{split}
\end{align}
where we used the GCL \eqref{GCL2} in the last step. Next, we apply the relation \eqref{PropertiesEntropyFluxFunctions} for the entropy flux functions and obtain 
\begin{equation}\label{EC:Condition1}
\mathbf{w}^{T}\pderivative{\mathbf{g}_{i}}{x_{i}}=\pderivative{}{x}\left(f_{i}^{s}-\nu_{i}s\right)-\left(\pderivative{v_{i}}{x_{i}}\right)\left(\mathbf{w}^{T}\mathbf{u}-s\right),\qquad i=1,2,3,
\end{equation}
by the same arguments as in the one dimensional computation \eqref{Spatialpart1D}. Next, we apply the vector notation $\vec{f}^{s}:=\left[f_{1}^{s},f_{2}^{s},f_{3}^{s}\right]^{T}$. 
Then \eqref{EC:Condition1} and the transformation formulas for the gradient and divergence in the Section \ref{Sec:TransCL} give  
\begin{equation}\label{EC:Condition3}
\mathbf{w}^{T}\left(\vec{\nabla}_{\xi}\cdot\blockvec{\tilde{\mathbf{g}}}\right)=\vec{\nabla}_{\xi}\cdot\left(\vec{\tilde{f}}^{s}-\vec{\tilde{\nu}}s\right)-\left(\vec{\nabla}_{\xi}\cdot\vec{\tilde{\nu}}\right)\left(\mathbf{w}^{T}\mathbf{u}-s\right).
\end{equation}
Finally, the identities \eqref{PropertiesEntropyRefE3} and \eqref{EC:Condition3} provide the balance law  
\begin{equation}\label{EntropyConditionReferenceElement}
0=\mathbf{w}^{T}\left(\frac{d\left(J\mathbf{u}\right)}{dt}+\vec{\nabla}_{\xi}\cdot\blockvec{\tilde{\mathbf{g}}}\right)=\frac{d\left(Js\right)}{dt}+\vec{\nabla}_{\xi}\cdot\left(\vec{\tilde{f}}^{s}-\vec{\tilde{\nu}}s\right).
\end{equation}
We integrate the equation \eqref{EntropyConditionReferenceElement} over the domain $E\times\left[0,T\right]$ and obtain
\begin{equation}\label{EQ:TotalEntropy}
\frac{d}{dt}\int_{E}Js\,d\vec{\xi}=-\int_{\partial E}\left(\vec{\tilde{f}}^{s}-\vec{\tilde{\nu}}s\right)^{T}\hat{n}\,dS.
\end{equation}
For discontinuous solutions the equation \eqref{EQ:TotalEntropy} is not satisfied, but it holds the inequality     
\begin{equation}\label{IQ:TotalEntropy}
\frac{d}{dt}\int_{E}Js\,d\vec{\xi}\leq-\int_{\partial E}\left(\vec{\tilde{f}}^{s}-\vec{\tilde{\nu}}s\right)^{T}\hat{n}\,dS.
\end{equation} 

\subsection{Buidling blocks for the spectral element approximation}\label{sec:Spectral}
The spectral element approximation based on a nodal approach with Lagrange basis functions given by    
\begin{equation}\label{LagrangeBasis}
\ell_{j}\left(\xi\right):=\prod_{i=0,\:j\neq i}^{N}\frac{\xi-\xi_{i}}{\xi_{j}-\xi_{i}},\quad j=0,\dots N,
\end{equation}     
where the nodal points $\left\{ \xi_{i}\right\} _{i=0}^{N}$ are the LGL points. We note that $\xi_{0}=-1$ and $\xi_{N}=1$. The Lagrange basis functions satisfy the cardinal property  
\begin{equation}\label{CardinalProperty}
\ell_{i}\left(\xi_{j}\right)=\delta_{ji},
\end{equation}  
where $\delta_{ji}$ is the Kronecker delta. On the reference element $E=\left[-1,1\right]^{3}$ the solution and fluxes of the system \eqref{eq:consLawBlockRefE} are approximated by tensor product Lagrange polynomials of degree $N$, e.g., 
\begin{equation}\label{SpatialPolynomial}
\mathbf{u}\left(\xi^{1},\xi^{2},\xi^{3},t\right)\approx\mathbf{U}\left(\xi^{1},\xi^{2},\xi^{3},t\right):=
\sum_{i,j,k=0}^{N}\mathbf{U}_{ijk}\left(t\right)\ell_{i}\left(\xi^{1}\right)\ell_{j}\left(\xi^{2}\right)\ell_{k}\left(\xi^{3}\right).
\end{equation}
In the following, polynomial approximations are highlighted  by capital letters, e.g., $\mathbf{U}$ is an approximation for the state vector $\textbf{u}$ and $\textbf{F}_{l}$, $l=1,2,3$, are approximations for the fluxes $\textbf{f}_{l}$, $l=1,2,3$. The determinant $J$ of the Jacobian matrix $\vec{\nabla}_{\vec{\xi}}\vec{\chi}$ is also approximated by tensor product Lagrange polynomials 
\begin{equation}\label{SpatialPolynomialMetricQuantity}
J\left(\xi^{1},\xi^{2},\xi^{3},t\right)\approx\mathrm{J}\left(\xi^{1},\xi^{2},\xi^{3},t\right):=\sum_{i,j,k=0}^{N}\mathrm{J}_{ijk}\left(t\right)\ell_{i}\left(\xi^{1}\right)\ell_{j}\left(\xi^{2}\right)\ell_{k}\left(\xi^{3}\right).
\end{equation}
In particular, the interpolation operator for a function $\mathbf{g}$ is given by 
\begin{equation}\label{InterpolationOperator}  
\interpolation{N}{\left(\mathbf{g}\right)}\left(\xi^{1},\xi^{2},\xi^{3}\right)=\sum_{i,j,k=0}^{N}\mathbf{g}_{ijk}\ell_{i}\left(\xi^{1}\right)\ell_{j}\left(\xi^{2}\right)\ell_{k}\left(\xi^{3}\right),
\end{equation}
where $\mathbf{g}_{ijk}:=\mathbf{g}\left(\xi_{i}^{1},\xi_{j}^{2},\xi_{k}^{3}\right)$ and $\left\{ \xi_{i}^{1}\right\} _{i=0}^{N}$, $\left\{\xi_{i}^{2}\right\} _{i=0}^{N}$, $\left\{\xi_{i}^{3}\right\} _{i=0}^{N}$ are sets of LGL points. Derivatives are approximated by exact differentiation of the polynomial interpolants. In general we have $ \left(\interpolation{N}{\left(g\right)}\right)'\neq \interpolation{N-1}{\left(g'\right)}$ (cf. e.g. \cite{CHQZ:2006,Kopriva2009}), as differentiation and interpolation only commute if there are no discretization errors. However, the contravariant coordinate vectors need to be discretized in such a way that the metric identities \eqref{MetricIdentities} are satisfied on the discrete level, too. Kopriva \cite{Kopriva2006} introduced the conservative curl form that computes
\begin{equation}\label{DiscreteContravariantVectors}
\mathrm{J}\vec{a}_{\beta}^{\alpha}:=-\hat{x}_{\alpha}\cdot\nabla_{\xi}\times\left(\interpolation{N}{\left(\chi_{\gamma}\nabla_{\delta}\chi_{m}\right)}\right),\quad\alpha=1,2,3,\quad\beta=1,2,3,\quad\left(\beta,\gamma,\delta\right)\ \text{cyclic},
\end{equation}
to approximate the metric terms. Here $\vec{\chi}=\left[\chi_{1},\chi_{2},\chi_{3}\right]^{T}$ represents the mapping from the element to the reference element and $\hat{x}_{i}$ is the unit vector in the $i$-th Cartesian coordinate direction. The representation \eqref{DiscreteContravariantVectors} ensures that   
\begin{equation}\label{DiscreteMetricIdentities1}
\sum_{\alpha=1}^{3}\pderivative{\interpolation{N}{\left(\mathrm{J}a_{\beta}^{\alpha}\right)}}{\xi^{\alpha}}=0,\qquad\beta=1,2,3.
\end{equation}
We note that in the LGL points $\xi_{i}^{1}$, $\xi_{j}^{2}$, $\xi_{k}^{3}$, $i,j,k=0,\dots,N$, the equation \eqref{DiscreteMetricIdentities1} gives 
\begin{equation}\label{DiscreteMetricIdentities}
\sum_{m=0}^{N}\left(\mathcal{D}_{im}\left(\mathrm{J}a_{\beta}^{1}\right)_{mjk}+\mathcal{D}_{jm}\left(\mathrm{J}a_{\beta}^{2}\right)_{imk}+\mathcal{D}_{km}\left(\mathrm{J}a_{\beta}^{3}\right)_{ijm}\right)=0,\qquad\beta=1,2,3.
\end{equation}
Integrals are approximated by a tensor product extension of a $2N-1$ accurate LGL quadrature formula. Hence, interpolation and quadrature nodes are collocated. In one spatial dimension the LGL quadrature formula is given by   
\begin{equation}\label{GLQ}
\int\limits_{-1}^{1}g\left(\xi\right)\,d\xi\approx\sum_{i=0}^{N}\omega_{i}g\left(\xi_{i}\right)=\sum_{i=0}^{N}\omega_{i}g_{i},  
\end{equation}
where $\omega_{i}$, $i=0,\dots,N$, are the quadrature weights and $\xi_{i}$, $i=0,\dots,N$, are the LGL quadrature points. The formula \eqref{GLQ} motivates the definition of the inner product notation 
\begin{equation}\label{Innerproduct}
\left\langle \mathbf{f},\mathbf{g}\right\rangle _{N}:=\sum_{i=0}^{N}\sum_{j=0}^{N}\sum_{k=0}^{N}\omega_{i}\omega_{j}\omega_{k}\mathbf{f}_{ijk}^{T}\mathbf{g}_{ijk}=\sum_{i,j,k=0}^{N}\omega_{ijk}\mathbf{f}_{ijk}^{T}\mathbf{g}_{ijk}
\end{equation}
for two functions $\mathbf{f}$ and $\mathbf{g}$. We note that the inner product \eqref{Innerproduct} satisfies        
\begin{equation}\label{MagicProperty}
\left\langle \interpolation{N}{\left(\mathbf{g}\right)},\boldsymbol{\varphi}\right\rangle _{N}=\left\langle \mathbf{g},\boldsymbol{\varphi}\right\rangle _{N}, \qquad \forall \boldsymbol{\varphi}\in\mathbb{P}^{N}\left(E,\R^{p}\right).
\end{equation}
Furthermore, for a block vector $\blockvec{\textbf{F}}$ and test functions $\boldsymbol{\varphi}\in\mathbb{P}^{N}\left(E,\R^{p}\right)$, we define the discrete surface integral   
\begin{align}\label{DiscreteSpatialSurface}
\begin{split}
\int\limits _{\partial E,N}\boldsymbol{\varphi}^T\left\{ \blockvec{\textbf{F}}\cdot\hat{n}\right\}\,dS 
:=&\quad\sum_{j,k=0}^{N}\omega_{j}\omega_{k}\left(\boldsymbol{\varphi}_{Njk}^{T}\left(\mathbf{F}_{1}\right)_{Njk}-\boldsymbol{\varphi}_{0jk}^{T}\left(\mathbf{F}_{1}\right)_{0jk}\right) \\
&+\sum_{i,k=0}^{N}\omega_{i}\omega_{k}\left(\boldsymbol{\varphi}_{iNk}^{T}\left(\mathbf{F}_{2}\right)_{iNk}-\boldsymbol{\varphi}_{i0k}^{T}\left(\mathbf{F}_{2}\right)_{i0k}\right) \\
& +\sum_{i,j=0}^{N}\omega_{i}\omega_{j}\left(\boldsymbol{\varphi}_{ijN}^{T}\left(\mathbf{F}_{3}\right)_{ijN}-\boldsymbol{\varphi}_{ij0}^{T}\left(\mathbf{F}_{3}\right)_{ij0}\right),
\end{split}
\end{align}
where $\hat{n}$ is the unit outward normal at the faces of the reference element $E$.

The spectral element approximation with LGL points for interpolation and quadrature  provides a SBP operator  $\matx{Q}=\matx{M}\,\matx{D}$ with the mass matrix $\matx{M}$ and the derivative matrix $\matx{D}$. 
The mass matrix and the derivative matrix are given by  
\begin{equation}\label{SBPcoefficient}
\mathcal{M}_{ij}=\omega_{i}\delta_{ij}, \qquad 
\mathcal{D}_{ij}=\ell_{j}'\left(\xi_{i}\right) \qquad 
i,j=0,\dots,N.
\end{equation}
A SBP operator satisfies the property     
\begin{equation}\label{SBP}
\matx{Q}+\matx{Q}^T=\matx{B},  
\end{equation}  
where $\matx{B}=\text{diag}\left(-1,0,\dots,0,1\right)$. A SBP operator provides a discrete analogue of the integration-by-parts formula \cite{Fernandez2014,Gassner2013,kreiss1}.   

\subsection{The semi-discrete discontinuous Galerkin method}\label{sec:SemiDiscreteDG}
Now, we apply the notation introduced in Section \ref{sec:Spectral} and construct a moving mesh DGSEM. Like for the moving mesh FV scheme \eqref{FVScheme}, we discretize the equations \eqref{GCL2} and \eqref{eq:consLawBlockRefE} simultaneous. In this way, it is ensured that the GCL \eqref{GCL2} is satisfied on the discrete level \cite{Kopriva2016, Minoli2011, Winters2014}. First, we replace the solution $\textbf{u}$ by \eqref{SpatialPolynomial}, the Jacobian 
$J$ by \eqref{SpatialPolynomialMetricQuantity} and approximate the fluxes by the interpolation operator \eqref{InterpolationOperator}. Next, we multiply the GCL \eqref{GCL2} with test functions $\varphi\in\mathbb{P}^{N}\left(E\right)$, the equation \eqref{eq:consLawBlockRefE} with 
$\boldsymbol{\varphi}\in\mathbb{P}^{N}\left(E,\R^{p}\right)$, integrate the resulting equations and use integration-by-parts to separate boundary and volume contributions. The volume integrals in the variational form are approximated with the LGL quadrature. Then, we insert numerical surface fluxes $\vec{\tilde{\nu}}^{*}$ and $\blockvec{\tilde{\mathbf{G}}}{}^{*}$ at the spatial element interfaces. Afterwards, we use the SBP property \eqref{SBP} for the volume contribution to get the standard DGSEM in strong form: 
\begin{subequations}\label{StandartDG}
\begin{align}
\left\langle \frac{d\mathrm{J}}{dt},\varphi\right\rangle_{N} 
=&\quad \left\langle \vec{\nabla}_{\xi}\cdot\interpolation{N}{\left(\vec{\tilde{\nu}}\right)},\varphi\right\rangle _{N}+\int\limits _{\partial E,N}\varphi\left(\tilde{\nu}_{\hat{n}}^{*}-\tilde{\nu}_{\hat{n}}\right)\,dS, \qquad \forall \varphi\in\mathbb{P}^{N}\left(E\right),  \\ 
\left\langle \frac{d\left(\mathrm{J}\mathbf{U}\right)}{dt},\boldsymbol{\varphi}\right\rangle_{N} 
=& -\left\langle \vec{\nabla}_{\xi}\cdot\interpolation{N}{\left(\blockvec{\tilde{\mathbf{G}}}\right)},\boldsymbol{\varphi}\right\rangle _{N}-\int\limits _{\partial E,N}\boldsymbol{\varphi}^{T}\left(\tilde{\mathbf{G}}_{\hat{n}}^{*}-\tilde{\mathbf{G}}_{\hat{n}}\right)\,dS, \qquad \forall 
\boldsymbol{\varphi}\in\mathbb{P}^{N}\left(E,\R^{p}\right).
\end{align}
\end{subequations}
where we used the inner product notation \eqref{Innerproduct} and the notation \eqref{DiscreteSpatialSurface} for the discrete surface integral. 

The approximation of $\vec{\tilde{\nu}}$ and the nonlinear flux $\blockvec{\tilde{\mathbf{g}}}$ by the interpolation operator \eqref{InterpolationOperator} causes aliasing errors in the standard strong form. The aliasing errors cannot be bounded and the errors are independent of the choice of the numerical surface flux. In Gassner \cite{Gassner2013} a detailed explanation and analysis of the aliasing problem is given. Furthermore, a specific reformulation of the volume integrals by using the skew-symmetry strategy has been developed to fix the aliasing problem. This approach has been enhanced and generalized by Gassner et al. in \cite{Gassner2017, Gassner2016} with a technique developed for high-order finite difference (FD) schemes (cf. Fisher and Carpenter \cite{Fisher2013}). The generalized approach is called split form DG framework. Here, we proceed similar as in \cite{Gassner2017} and replace the interpolation operators in the discrete volume integrals by derivative projection operators. The interpolation operator in the discrete equation for the GCL \eqref{GCL2} is replaced by 
\begin{align}\label{SpatialDerivativeProjectionOperator1}
\begin{split}
\Dprojection{N}\cdot\vec{\tilde{\nu}}_{ijk}:=\sum_{m=0}^{N}\quad\,2\,\mathcal{D}_{im}\avg{\vec{\nu}}_{\left(i,m\right)jk}\cdot\avg{\mathrm{J}\vec{a}^{1}}_{\left(i,m\right)jk} \\
\quad\ \ \quad+\,2\,\mathcal{D}_{jm}\avg{\vec{\nu}}_{i\left(j,m\right)k}\cdot\avg{\mathrm{J}\vec{a}^{2}}_{i\left(j,m\right)k} \\
\quad\ \ \quad+\,2\,\mathcal{D}_{km}\avg{\vec{\nu}}_{ij\left(k,m\right)}\cdot\avg{\mathrm{J}\vec{a}^{3}}_{ij\left(k,m\right)}
\end{split}
\end{align}
with the volume averages of the metric terms, e.g.  
\begin{equation}\label{SpatialVolumeAverages}
\avg{\cdot}_{(i,m)jk}:=\frac{1}{2}\left[\left(\cdot\right)_{ ijk}+\left(\cdot\right)_{ mjk}\right].
\end{equation}
The derivative projection operator in the discrete equation for \eqref{eq:consLawBlockRefE} is computed as in \cite{Gassner2017}. Thus, the operator is given by
\begin{align}\label{SpatialDerivativeProjectionOperator2}
\begin{split}
\Dprojection{N}\cdot\blockvec{\tilde{\mathbf{G}}}{}_{ijk}^{\text{EC}}:=&\sum_{m=0}^{N}\quad\,2\,\mathcal{D}_{im}\left(\blockvec{\mathbf{G}}{}^{\text{EC}}\left(\vec{\nu}_{ijk},\vec{\nu}_{mjk},\mathbf{U}_{ijk},\mathbf{U}_{mjk}\right)\cdot\avg{\mathrm{J}\vec{a}^{1}}_{\left(i,m\right)jk}\right)\\
& \quad \ \ \quad+2\,\mathcal{D}_{jm}\left(\blockvec{\mathbf{G}}{}^{\text{EC}}\left(\vec{\nu}_{ijk},\vec{\nu}_{imk},\mathbf{U}_{ijk},\mathbf{U}_{imk}\right)\cdot\avg{\mathrm{J}\vec{a}^{2}}_{i\left(j,m\right)k}\right) \\
& \quad \ \ \quad+2\,\mathcal{D}_{km}\left(\blockvec{\mathbf{G}}{}^{\text{EC}}\left(\vec{\nu}_{ijk},\vec{\nu}_{ijm},\mathbf{U}_{ijk},\mathbf{U}_{ijm}\right)\cdot\avg{\mathrm{J}\vec{a}^{3}}_{ij\left(k,m\right)}\right).
\end{split}
\end{align}
The flux $\blockvec{\mathbf{G}}{}^{\text{EC}}$ is consistent and symmetric such that, e.g.
\begin{align}\label{ConsistentDGECFlux}
\begin{split}
\blockvec{\mathbf{G}}{}^{\text{EC}}\left(\vec{\nu}_{ijk},\vec{\nu}_{mjk},\mathbf{U},\mathbf{U}\right)=\blockvec{\mathbf{F}}\left(\mathbf{U}\right)-\avg{\vec{v}}_{\left(i,m\right)jk}\mathbf{U},
\end{split}
\end{align}
and
\begin{equation}\label{FluxSymmetric}
\blockvec{\mathbf{G}}{}^{\text{EC}}\left(\vec{\nu}_{ijk},\vec{\nu}_{mjk},\mathbf{U}_{ijk},\mathbf{U}_{mjk}\right)=\blockvec{\mathbf{G}}{}^{\text{EC}}\left(\vec{\nu}_{mjk},\vec{\nu}_{ijk},\mathbf{U}_{mjk},\mathbf{U}_{ijk}\right),
\end{equation}   
for $i,j,k,m=0,\dots,N$. Furthermore, the flux functions $\mathbf{G}_{l}^{\text{EC}}$, $l=1,2,3$, satisfy in the  LGL points the moving mesh FV entropy condition \eqref{FV:EntropyCondition}. More precisely for $i,j,k,m=0,\dots,N$, the following discrete entropy conditions are satisfied 
\begin{align}\label{EntropyConditionDGECFlux}
\begin{split}
& \jump{\mathbf{W}}_{(i,m)jk}^{T}\mathbf{G}_{l}^{\text{EC}}\left(\vec{\nu}_{ijk},\vec{\nu}_{mjk},\mathbf{U}_{ijk},\mathbf{U}_{mjk}\right)=\jump{\Psi_{l}}_{(i,m)jk}-\avg{\nu_{l}}_{\left(i,m\right)jk}\jump{\Phi}_{(i,m)jk}, \\
& \jump{\mathbf{W}}_{i(j,m)k}^{T}\mathbf{G}_{l}^{\text{EC}}\left(\vec{\nu}_{ijk},\vec{\nu}_{imk},\mathbf{U}_{ijk},\mathbf{U}_{imk}\right)=\jump{\Psi_{l}}_{i(j,m)k}-\avg{\nu_{l}}_{i\left(j,m\right)k}\jump{\Phi}_{i(j,m)k}, \\
& \jump{\mathbf{W}}_{ij(k,m)}^{T}\mathbf{G}_{l}^{\text{EC}}\left(\vec{\nu}_{ijk},\vec{\nu}_{ijm},\mathbf{U}_{ijk},\mathbf{U}_{ijm}\right)=\jump{\Psi_{l}}_{ij(k,m)}-\avg{\nu_{l}}_{ij\left(k,m\right)}\jump{\Phi}_{ij(k,m)}.
\end{split}
\end{align}
The quantities $\Phi$ and $\Psi_{l}$  are polynomial approximations which satisfy in the LGL points      
\begin{equation}\label{PolynomialApproximation1}
\Phi_{ijk}=\left[\mathbf{w}\left(\mathbf{U}\right)\right]_{ijk}^{T}\mathbf{U}_{ijk}-S_{ijk}, \qquad 
\left(\Psi_{l}\right)_{ijk}:=\left[\mathbf{w}\left(\mathbf{U}\right)\right]_{ijk}^{T}\left(\mathbf{F}_{l}\right)_{ijk}-\left(F_{l}^{s}\right)_{ijk}, \qquad l=1,2,3, 
\end{equation}
where $\mathbf{W}_{ijk}$, $S_{ijk}$ and $\left(F_{l}^{s}\right)_{ijk}$ are the nodal values of the polynomials      
\begin{equation}\label{PolynomialApproximation2}
\mathbf{W}:=\interpolation{N}{\left(\mathbf{w}\left(\mathbf{U}\right)\right)}, \qquad
S:=\interpolation{N}{\left(s\left(\mathbf{U}\right)\right)}, \qquad 
F_{l}^{s}:=\interpolation{N}{\left(f_{l}^{s}\left(\mathbf{U}\right)\right)}, \qquad l=1,2,3. 
\end{equation}
Here, $s$ represents an entropy for the system \eqref{eq:consLaw} with the corresponding entropy flux functions $f_{l}^{s}$, $l=1,2,3$, and entropy variables $\textbf{w}$. Furthermore, the volume jumps in \eqref{EntropyConditionDGECFlux} are, e.g. 
\begin{equation}\label{SpatialVolumeJumps}
\jump{\cdot}_{(i,m)jk}:=\left(\cdot\right)_{ijk}-\left(\cdot\right)_{mjk}.
\end{equation}
In the Appendices \ref{sec:App Shallow Water Flux} and \ref{sec:App Euler}, flux functions with these properties are presented for the shallow water and  Euler equations.   

Finally, for each element $e_{\kappa}(t)$ the semi-discrete moving mesh DGSEM can be written in the following form:      
\begin{subequations}\label{MovingMeshDGSEM}
\begin{align}
\left\langle \frac{d\mathrm{J}}{dt},\varphi\right\rangle_{N} 
=&\quad \left\langle \Dprojection{N}\cdot\vec{\tilde{\nu}},\varphi\right\rangle _{N}+\int\limits _{\partial E,N}\varphi\left(\tilde{\nu}_{\hat{n}}^{*}-\tilde{\nu}_{\hat{n}}\right)\,dS, \qquad \forall \varphi\in\mathbb{P}^{N}\left(E\right), \label{DGSEM:D-GLC} \\ 
\left\langle \frac{d\left(\mathrm{J}\mathbf{U}\right)}{dt},\boldsymbol{\varphi}\right\rangle_{N} 
=& -\left\langle \Dprojection{N}\cdot\blockvec{\tilde{\mathbf{G}}}{}^{\text{EC}},\boldsymbol{\varphi}\right\rangle _{N}-\int\limits _{\partial E,N}\boldsymbol{\varphi}^{T}\left(\tilde{\mathbf{G}}_{\hat{n}}^{*}-\tilde{\mathbf{G}}_{\hat{n}}\right)\,dS, \qquad \forall 
\boldsymbol{\varphi}\in\mathbb{P}^{N}\left(E,\R^{p}\right). \label{DGSEM:CL}
\end{align}
\end{subequations}
We note that the equation \eqref{DGSEM:D-GLC} is the D-GCL for the moving mesh spectral element DG discretization.  

The unit outward facing normal vector and surface element on the element side are constructed from the element metrics by
\begin{equation}\label{ComputationNormal}
\vec{n}:=\frac{1}{\hat{s}}\sum_{l=1}^{3}\left(\mathrm{J}\vec{a}^{l}\right)\hat{n}^{l},\qquad\hat{s}:=\left|\sum_{l=1}^{3}\left(\mathrm{J}\vec{a}^{l}\right)\hat{n}^{l}\right|.
\end{equation}
Thus, the quantity $\tilde{\nu}_{\hat{n}}$ in \eqref{DGSEM:D-GLC} and the flux $\tilde{\mathbf{G}}_{\hat{n}}$ in \eqref{DGSEM:CL} are defined by 
\begin{align}
\tilde{\nu}_{\hat{n}}=&\left(\hat{s}\vec{n}\right)\cdot\vec{\nu}=\sum_{l=1}^{3}\hat{n}^{l}\left(\mathrm{J}a_{1}^{l}\nu_{1}+\mathrm{J}a_{2}^{l}\nu_{2}+\mathrm{J}a_{3}^{l}\nu_{3}\right), \label{GridvelocityNormal} \\
\tilde{\mathbf{G}}_{\hat{n}}=&\left(\hat{s}\vec{n}\right)\cdot\blockvec{\mathbf{G}}=\sum_{l=1}^{3}\hat{n}^{l}\left(\mathrm{J}a_{1}^{l}\mathbf{G}_{1}+\mathrm{J}a_{2}^{l}\mathbf{G}_{2}+\mathrm{J}a_{3}^{l}\mathbf{G}_{3}\right)=\left\{ \blockmatx{M}\,\blockvec{\mathbf{G}}\right\} \cdot\hat{n}. \label{BlockvectorNormal}
\end{align}
To define the numerical surface fluxes in \eqref{DGSEM:D-GLC} and \eqref{DGSEM:CL}, we introduce notation for states at the LGL nodes along an interface between two spatial elements to be a primary ``$-$'' and complement the notation with a secondary ``$+$'' to denote the value at the LGL nodes on the opposite side.
Then the orientated jump and
the arithmetic mean at the interfaces are defined by  
\begin{equation}\label{SurfaceJumpMean}
\jump{\cdot}:=\left(\cdot\right)^{+}-\left(\cdot\right)^{-},
\quad \text{and} \quad 
\avg{\cdot}:=\frac{1}{2}\left[\left(\cdot\right)^{+}+\left(\cdot\right)^{-}\right].  
\end{equation}
When applied to vectors, the average and jump operators are evaluated separately for each vector component. Then the normal vector $\vec{n}$ is defined unique to point from the ``$-$'' to the ``$+$'' side. This notation allows to compute the  contravariant surface numerical fluxes in \eqref{DGSEM:D-GLC} as
\begin{equation}
\tilde{\nu}_{\hat{n}}^{*}=\hat{s}\left(n_{1}\avg{v_{1}}+n_{2}\avg{v_{2}}+n_{3}\avg{v_{3}}\right).
\end{equation}
The contravariant surface numerical fluxes in \eqref{DGSEM:CL} are given by  
\begin{equation}\label{ContravariantSurfaceFlux}
\tilde{\mathbf{G}}_{\hat{n}}^{*}=\hat{s}\left(n_{1}\mathbf{G}_{1}^{\text{EC}}+n_{2}\mathbf{G}_{2}^{\text{EC}}+n_{3}\mathbf{G}_{3}^{\text{EC}}\right),
\end{equation}
where the Cartesian fluxes  $\mathbf{G}_{l}^{\text{EC}}$, $l=1,2,3$, satisfy \eqref{ConsistentDGECFlux}, \eqref{FluxSymmetric}, \eqref{EntropyConditionDGECFlux}. 
We note that these fluxes are the baseline choices without interface dissipation, to get a baseline scheme that is entropy conservative.  
\begin{rem-hand}[3.2.]
The discrete volume weighted contravariant vectors  $\mathrm{J}\vec{a}^{\alpha}$, $\alpha=1,2,3$, do not dependent on the solution $\mathrm{J}$ of \eqref{DGSEM:D-GLC}, since these vectors are computed by the conservative curl form \eqref{DiscreteContravariantVectors}. Thus, the discrete metric identities \eqref{DiscreteMetricIdentities1} are satisfied and the normal computation \eqref{ComputationNormal} is watertight. This means the normal vector and the surface element are continuous across element interfaces.
\end{rem-hand}

\subsection{Discrete entropy conservation for the semi-discrete method}\label{EntropyConservation}\label{sec:EntropyConservation}
The spatial integral of the entropy is bounded in time on the continuous level. Thus, it is desirable that a numerical method is stable in the sense that a discrete version of this integral is bounded in time, too. In the context of the moving mesh semi-discrete DGSEM \eqref{MovingMeshDGSEM}, we are interested to find an upper bound for the quantity 
\begin{equation}\label{DGSEM:DisreteEntropyIntegral}
\bar{S}:=\sum_{k=1}^{K}\left\langle s\left(\mathbf{U}\right),\mathrm{J}\right\rangle _{N}=\sum_{k=1}^{K}\left\langle \interpolation{N}{\left(s\left(\mathbf{U}\right)\right)},\mathrm{J}\right\rangle _{N}=\sum_{k=1}^{K}\left\langle S,\mathrm{J}\right\rangle _{N},
\end{equation}
where we used the property \eqref{MagicProperty} of the inner product \eqref{Innerproduct} and the notation \eqref{PolynomialApproximation2}. Next, we prove the following statement for the semi-discrete moving mesh DGSEM.     
\begin{thm-hand}[3.3.]
Suppose the flux functions $\blockvec{\mathbf{G}}{}^{\text{EC}}$ in the derivative projection operator \eqref{SpatialDerivativeProjectionOperator2}  and the numerical surface fluxes $\tilde{\mathbf{G}}_{\hat{n}}^{*}$ are computed by Cartesian fluxes  $\mathbf{G}_{l}^{\text{EC}}$, $l=1,2,3$, with the properties \eqref{ConsistentDGECFlux}, \eqref{FluxSymmetric}, \eqref{EntropyConditionDGECFlux}. Then the semi-discrete moving mesh DGSEM \eqref{MovingMeshDGSEM} satisfies the discrete entropy equation
\begin{equation}\label{MovingMeshEC}
\frac{d\bar{S}}{dt}=
\underset{\text{faces}}
{\sum_{\text{Boundary}}}\ \int\limits _{\partial E,N}\left(\tilde{\Psi}_{\hat{n}}-\tilde{\nu}_{\hat{n}}^{*}\Phi-\mathbf{W}^{T}\tilde{\mathbf{G}}_{\hat{n}}^{*}\right)\,dS,
\end{equation}
where $\bar{S}$ is given by \eqref{DGSEM:DisreteEntropyIntegral}, $\Phi$ as well as $\Psi_{l}$, $l=1,2,3$, are polynomials with nodal values \eqref{PolynomialApproximation1} and $\tilde{\Psi}_{\hat{n}}:=\left(\hat{s}\vec{n}\right)\cdot\vec{\Psi}$ with $\vec{\Psi}=\left[\Psi_{1},\Psi_{2},\Psi_{3}\right]^{T}$.  
\end{thm-hand}
\begin{proof}
We proceed similar as in the continuous entropy analysis, use the polynomial approximation $\boldsymbol{\varphi}=\interpolation{N}{\left(\mathbf{w}\left(\mathbf{U}\right)\right)}=\mathbf{W}$ as test function in the equation \eqref{DGSEM:CL} and obtain  
\begin{equation}\label{DGSEM:EntropyPreservation1}
\left\langle \frac{d\left(\mathrm{J}\mathbf{U}\right)}{dt},\mathbf{W}\right\rangle _{N}=-\left\langle \Dprojection{N}\cdot\blockvec{\tilde{\mathbf{G}}}{}^{\text{EC}},\mathbf{W}\right\rangle _{N}-\int\limits _{\partial E,N}\mathbf{W}^{T}\left(\tilde{\mathbf{G}}_{\hat{n}}^{*}-\tilde{\mathbf{G}}_{\hat{n}}\right)\,dS.
\end{equation}
First, we consider the left hand side in the equation \eqref{DGSEM:EntropyPreservation1}. We obtain by the chain rule   
\begin{align}\label{DGSEM:EntropyPreservation2}
\begin{split}
\left\langle \mathrm{J}\frac{d\mathbf{U}}{dt},\mathbf{W}\right\rangle _{N}
=&\left\langle \mathbf{w}\left(\mathbf{U}\right)^{T}\pderivative{\mathbf{U}}{t},\mathrm{J}\right\rangle _{N}+\left\langle \vec{\tilde{\nu}}\cdot\left(\mathbf{w}\left(\mathbf{U}\right)^{T}\vec{\nabla}_{\xi}\mathbf{U}\right),1\right\rangle _{N} \\
=&\left\langle \pderivative{}{t}s\left(\mathbf{U}\right),\mathrm{J}\right\rangle _{N}+\left\langle \vec{\tilde{\nu}}\cdot\vec{\nabla}_{\xi}s\left(\mathbf{U}\right),1\right\rangle _{N}=\left\langle \frac{d}{dt}s\left(\mathbf{U}\right),\mathrm{J}\right\rangle _{N},
\end{split}
\end{align}
where we used the identity $\interpolation{N}{\left(\mathbf{w}\left(\mathbf{U}\right)\right)}=\mathbf{W}$ and the property \eqref{MagicProperty} of the inner product \eqref{Innerproduct}. Since we assume time continuity for our semi-discrete analysis, we apply the product rule in time and obtain by \eqref{DGSEM:EntropyPreservation2} 
\begin{align}\label{DGSEM:EntropyPreservation3}
\begin{split}
\left\langle \frac{d\left(\mathrm{J}\mathbf{U}\right)}{dt},\mathbf{W}\right\rangle _{N}=&\quad\left\langle \frac{ds\left(\mathbf{U}\right)}{dt},\mathrm{J}\right\rangle _{N}+\left\langle \frac{dJ}{dt},\mathbf{w}\left(\mathbf{U}\right)^{T}\mathbf{U}\right\rangle _{N} \\
=&\quad\frac{d}{dt}\left\langle s\left(\mathbf{U}\right),J\right\rangle _{N}+\left\langle \frac{dJ}{dt},\mathbf{w}\left(\mathbf{U}\right)^{T}\mathbf{U}-s\left(\mathbf{U}\right)\right\rangle _{N}.
\end{split}
\end{align}
Next, the notations in \eqref{PolynomialApproximation1}, \eqref{PolynomialApproximation2} and the property \eqref{MagicProperty} of the inner product \eqref{Innerproduct} give 
\begin{align}\label{DGSEM:EntropyPreservation4}
\begin{split}
\left\langle \frac{d\left(\mathrm{J}\mathbf{U}\right)}{dt},\mathbf{W}\right\rangle _{N} 
=& \quad\frac{d}{dt}\left\langle s\left(\mathbf{U}\right),J\right\rangle _{N}+\left\langle \frac{dJ}{dt},\mathbf{w}\left(\mathbf{U}\right)^{T}\mathbf{U}-s\left(\mathbf{U}\right)\right\rangle _{N}=\frac{d}{dt}\left\langle S,J\right\rangle _{N}+\left\langle \frac{dJ}{dt},\Phi\right\rangle _{N} \\
=&\quad\frac{d}{dt}\left\langle S,J\right\rangle _{N}+\left\langle \Dprojection{N}\cdot\vec{\tilde{\nu}},\Phi\right\rangle _{N}+\int\limits _{\partial E,N}\left(\tilde{\nu}_{\hat{n}}^{*}-\tilde{\nu}_{\hat{n}}\right)\Phi\,dS,
\end{split}
\end{align}
where we used in the last step the D-GCL \eqref{DGSEM:D-GLC} with the test function $\varphi=\Phi$. We note that the quantity $\Phi$ is defined as a polynomial with the nodal values \eqref{PolynomialApproximation1}. In the Appendix \ref{DGSEM_EC}, the following equation is proven 
\begin{equation}\label{DGSEM:EntropyPreservation5}
\left\langle \Dprojection{N}\cdot\blockvec{\tilde{\mathbf{G}}}{}^{\text{EC}},\mathbf{W}\right\rangle _{N}=\int\limits _{\partial E,N}\left(\tilde{F}_{\hat{n}}^{s}-\tilde{\nu}_{\hat{n}}S\right)\,dS-\left\langle \Dprojection{N}\cdot\vec{\tilde{\nu}},\Phi\right\rangle _{N},
\end{equation}
where $\tilde{F}_{\hat{n}}^{s}=\left(\hat{s}\vec{n}\right)\cdot\vec{F}^{s}$ with $\vec{F}^{s}=\left[F_{1}^{s},F_{2}^{s},F_{3}^{s}\right]^{T}$. Here the polynomials $F_{l}^{s}$, $l=1,2,3$, are given by \eqref{PolynomialApproximation2}. Moreover, we obtain by \eqref{GridvelocityNormal} and \eqref{BlockvectorNormal} 
\begin{align}\label{DGSEM:EntropyPreservation6}
\begin{split}
& -\mathbf{W}^{T}\left(\tilde{\mathbf{G}}_{\hat{n}}^{*}-\tilde{\mathbf{G}}_{\hat{n}}\right)-\left(\tilde{F}_{\hat{n}}^{s}-\tilde{\nu}_{\hat{n}}S\right) \\
=&\quad\sum_{l=1}^{3}\left\{ \hat{s}n_{l}\left(\mathbf{W}^{T}\mathbf{F}_{l}-F_{l}^{s}\right)-\hat{s}n_{l}\left(\mathbf{W}^{T}\mathbf{U}-S\right)\right\} -\mathbf{W}^{T}\tilde{\mathbf{G}}_{\hat{n}}^{*} \\
=&\quad\tilde{\Psi}_{\hat{n}}-\tilde{\nu}_{\hat{n}}\Phi-\mathbf{W}^{T}\tilde{\mathbf{G}}_{\hat{n}}^{*}.
\end{split}
\end{align}  
Next, we plug the equations \eqref{DGSEM:EntropyPreservation4}, \eqref{DGSEM:EntropyPreservation5}, \eqref{DGSEM:EntropyPreservation6} in 
\eqref{DGSEM:EntropyPreservation1} and rearrange. This results in the equation 
\begin{align}\label{DGSEM:EntropyPreservation7}
\begin{split}
\frac{d}{dt}\left\langle S,J\right\rangle _{N}
=&-\int\limits _{\partial E,N}\left\{ \mathbf{W}^{T}\left(\tilde{\mathbf{G}}_{\hat{n}}^{*}-\tilde{\mathbf{G}}_{\hat{n}}\right)+\left(\tilde{F}_{\hat{n}}^{s}-\tilde{\nu}_{\hat{n}}S\right)-\left(\tilde{\nu}_{\hat{n}}^{*}-\tilde{\nu}_{\hat{n}}\right)\Phi\right\} \,dS \\
=&\quad \int\limits _{\partial E,N}\left(\tilde{\Psi}_{\hat{n}}-\tilde{\nu}_{\hat{n}}^{*}\Phi-\mathbf{W}^{T}\tilde{\mathbf{G}}_{\hat{n}}^{*}\right)\,dS.
\end{split}
\end{align} 
Then, we sum the equation \eqref{DGSEM:EntropyPreservation7} over all elements and use that the normal computation \eqref{ComputationNormal} is watertight. This provides the equation 
\begin{align}\label{DGSEM:EntropyPreservation8}
\begin{split}
\frac{d\bar{S} }{dt}
=& \quad \underset{\text{faces}}{\sum_{\text{Boundary}}}\ \int\limits _{\partial E,N}\left(\tilde{\Psi}_{\hat{n}}-{\nu}_{\hat{n}}^{*}\Phi\tilde-\mathbf{W}^{T}\tilde{\mathbf{G}}_{\hat{n}}^{*}\right)\,dS \\
& -\underset{\text{faces}}{\sum_{\text{Interior}}}\int\limits _{\partial E,N}\left(\jump{\tilde{\Psi}_{\hat{n}}}-\avg{\tilde{\nu}_{\hat{n}}}\jump{\Phi}-\jump{\mathbf{W}}^{T}\tilde{\mathbf{G}}_{\hat{n}}^{*}\right)\,dS.
\end{split}
\end{align}
Since the numerical surface fluxes $\tilde{\mathbf{G}}_{\hat{n}}^{*}$ are computed by Cartesian fluxes  $\mathbf{G}_{l}^{\text{EC}}$, $l=1,2,3$, with the properties \eqref{EntropyConditionDGECFlux}, it follows  
\begin{equation}
\jump{\tilde{\Psi}_{\hat{n}}}-\avg{\tilde{\nu}_{\hat{n}}}\jump{\Phi}-\jump{\mathbf{W}}^{T}\tilde{\mathbf{G}}_{\hat{n}}^{*}
=\sum_{l=1}^{3}\hat{s}n_{l}\left(\jump{\Psi_{l}}-\avg{\nu_{l}}\jump{\Phi}-\jump{\mathbf{W}}^{T}\mathbf{G}_{l}^{*}\right)=0.
\end{equation}
Hence, we obtain the equation 
\begin{equation}\label{DGSEM:EntropyPreservation9}
\frac{d\bar{S}}{dt} 
= \underset{\text{faces}}{\sum_{\text{Boundary}}}\ \int\limits _{\partial E,N}\left(\tilde{\Psi}_{\hat{n}}-\tilde{\nu}_{\hat{n}}^{*}\Phi-\mathbf{W}^{T}\tilde{\mathbf{G}}_{\hat{n}}^{*}\right)\,dS.
\end{equation}
This completes the proof of Theorem 3.3.
\end{proof}
\begin{rem-hand}[3.4.]
The boundary contribution in \eqref{MovingMeshEC} becomes zero in the case of periodic boundary conditions. In this case the entropy is conserved. 
\end{rem-hand}

\subsection{Discrete entropy stability for the semi-discrete method}\label{EntropyStability}\label{sec:EntropyStability}
Entropy conservation can be merely expected when a reversible process is described by a system of PDEs. In general, conservation laws are describing irreversible processes with discontinuous solutions. Hence, it cannot be expected that the entropy conservative moving mesh DGSEM provides a physical meaningful discretization for the system \eqref{eq:consLaw}. However, the entropy conservative flux can be augmented by an artificial dissipation term. 

In the literature, there are different strategies to add dissipation to an entropy conservative flux. Here, dissipation is added via a matrix operator. This approach, for instance, has been used in the context of gas dynamics by Chandrashekar \cite{Chandrashekar2013} or Winters et al. \cite{Winters2017}. 

The conservative variables $\textbf{u}$ can be written in dependence of the entropy variables $\textbf{w}$. Differentiation of the conservative variables $\textbf{u}=\textbf{u}(\textbf{w})$ provides the symmetric positive definite matrix $\pderivative{\textbf{u}}{\textbf{w}}$, since the system \eqref{eq:consLaw} is assumed to be symmetrizable (cf. e.g. \cite{Harten1983}). Thus, it follows by a Taylor expansion up to first order   
\begin{equation}\label{JumpEntropy}
\jump{\textbf{u}}=\left(\pderivative{\textbf{u}}{\textbf{w}}\right)\jump{\textbf{w}}+\mathcal{O}\left(\left|\jump{\textbf{w}}\right|^{2}\right),
\end{equation}
where the jump operator is defined by \eqref{SurfaceJumpMean} at the interfaces. Furthermore, the system \eqref{eq:consLaw} is hyperbolic. Thus, the flux Jacobian matrices $\pderivative{\textbf{f}_{l}}{\textbf{u}}$, $l=1,2,3$, are diagonalizable and have real eigenvalues $\left\{ \lambda_{i}^{l}\left(\textbf{u}\right)\right\} _{i=1}^{p}\subseteq\R$. The corresponding right eigenvector matrices are $\matx{R}_{l}$. According to the eigenvector scaling theorem (cf. Barth \cite[Theorem 4]{Barth2018}), there are symmetric block diagonal scaling matrices $\matx{T}_{l}$ with 
\begin{equation}\label{ScaledEigenvectorMatrix}
\pderivative{\mathbf{f}_{l}}{\mathbf{u}}=\tilde{\matx{R}}_{l}\,\Lambda_{l}\left(\mathbf{u}\right)\,\tilde{\matx{R}}_{l}^{-1},\qquad\pderivative{\mathbf{u}}{\mathbf{w}}=\tilde{\matx{R}}_{l}\,\tilde{\matx{R}}_{l}^{T},\qquad\tilde{\matx{R}}_{l}=\matx{R}_{l}\,\matx{T}_{l},\qquad l=1,2,3,
\end{equation}
where $\Lambda_{l}\left(\mathbf{u}\right):=\text{diag}\left(\lambda_{1}^{l}\left(\mathbf{u}\right),\dots,\lambda_{p}^{l}\left(\mathbf{u}\right)\right)$.
The flux Jacobian matrices $\pderivative{\textbf{g}_{l}}{\textbf{u}}=\pderivative{\textbf{f}_{l}}{\textbf{u}}-\nu_{l}\matx{I}_{p}$  have the real eigenvalues $\left\{ \lambda_{i}^{l}\left(\mathbf{u}\right)-\nu_{l}\right\} _{i=1}^{p}$ and the same right eigenvectors as the flux Jacobian $\pderivative{\textbf{f}_{l}}{\textbf{u}}$. We note that $\matx{I}_{p}$ is the $p\times p$ identity matrix. Hence, it follows  
\begin{equation}\label{DiagonalJacobianG}
\pderivative{\mathbf{g}_{l}}{\mathbf{u}}=\tilde{\matx{R}_{l}}\,\Lambda_{l}\left(\nu,\mathbf{u}\right)\,\tilde{\matx{R}}_{l}^{-1},\qquad\Lambda_{l}\left(\nu,\mathbf{u}\right):=\text{diag}\left(\lambda_{1}^{l}\left(\mathbf{u}\right)-\nu_{l},\dots,\lambda_{p}^{l}\left(\mathbf{u}\right)-\nu_{l}\right), \qquad l=1,2,3.
\end{equation}
Furthermore, we obtain by \eqref{ScaledEigenvectorMatrix}
 \begin{equation}\label{DecompositionJacobianG}
\left(\pderivative{\mathbf{g}_{l}}{\mathbf{w}}\right)=\left(\pderivative{\mathbf{g}_{l}}{\mathbf{u}}\right)\left(\pderivative{\mathbf{u}}{\mathbf{w}}\right)=\tilde{\matx{R}}_{l}\,\Lambda_{l}\left(\nu,\mathbf{u}\right)\,\tilde{\matx{R}}_{l}^{T},\qquad l=1,2,3. 
\end{equation}      

The equation \eqref{DecompositionJacobianG} motivates the definition of the following matrix dissipation operators      
\begin{equation}\label{EntropyBasedDissipation}
\matx{H}_{l}=\hat{\matx{R}}_{l}\,\left|\Lambda_{l}\right|\,\hat{\matx{R}}_{l}^{T},
\qquad \hat{\matx{R}}_{l}=\matx{R}_{l}^{\star}\matx{T}_{l}^{\star}
,\qquad l=1,2,3.
\end{equation}
where the matrices $\matx{R}_{l}^{\star}$, $\matx{T}_{l}^{\star}$, depend on the averaged values of the states $\mathbf{U}^{-}$, $\mathbf{U}^{+}$ and they are consistent with the right eigenvector matrix $\matx{R}_{l}$ and the scaling matrix $\matx{T}_{l}$. The matrix $\left|\Lambda_{l}\right|$ depends on the values $\left\{ \lambda_{i}^{l}\left(\mathbf{U}^{-}\right)-\nu_{l}^{-}\right\} _{i=1}^{p}$ and $\left\{ \lambda_{i}^{l}\left(\mathbf{U}^{+}\right)-\nu_{l}^{+}\right\} _{i=1}^{p}$. The matrix $\matx{H}_{l}$ needs to be a symmetric positive definite matrix. Therefore, the matrix $\left|\Lambda_{l}\right|$ has to be choosen carefully. Examples of matrices which lead to a symmetric positive definite dissipation matrix are Ismail-Roe type dissipation~\cite{Ismail2009}    
\begin{equation}
\left|\Lambda_{l}\right|^{\text{Roe}}:=\text{diag}\left(\left|\avg{\lambda_{1}^{l}\left(\mathbf{U}\right)-\nu_{l}}\right|,\dots,\left|\avg{\lambda_{p}^{l}\left(\mathbf{U}\right)-\nu_{l}}\right|\right),  \qquad l=1,2,3,
\end{equation}
or Rusanov type dissipation     
\begin{equation}
\left|\Lambda_{l}\right|^{\text{Rus}}:=\max\left(\underset{1\leq i\leq p}{\max}\left|\lambda_{i}^{l}\left(\mathbf{U}^{-}\right)-\nu_{l}^{-}\right|,\underset{1\leq i\leq p}{\max}\left|\lambda_{i}^{l}\left(\mathbf{U}^{+}\right)-\nu_{l}^{+}\right|\right)\matx{I}_{p}, \qquad l=1,2,3.
\end{equation}
The Rusanov type dissipation is known to be quite dissipative and leads to smearing of shocks and contact discontinuities. The Roe type dissipation is more accurate in the reproduction of the characteristic profile of shocks or contact discontinuities, but in the regime of strong shock fronts the Roe type dissipation fails in the  stabilization of the numerical method. Therefore, a convex combination of Roe and Rusanov type dissipation is a natural choice \cite{Chandrashekar2013}. 

The dissipation operator \eqref{EntropyBasedDissipation} is used to modify the Cartesian numerical surface flux as follows 
\begin{equation}\label{EntropyStableSFFlux}
\mathbf{G}_{l}^{\text{ES}}:=\mathbf{G}_{l}^{\text{EC}}-\frac{1}{2}\,\matx{H}_{l}\,\jump{\mathbf{W}}, \qquad l=1,2,3,
\end{equation}
where the Cartesian fluxes  $\mathbf{G}_{l}^{\text{EC}}$, $l=1,2,3$, satisfy \eqref{ConsistentDGECFlux}, \eqref{FluxSymmetric}, \eqref{EntropyConditionDGECFlux}. The contravariant surface numerical fluxes $\tilde{\mathbf{G}}_{\hat{n}}^{\text{ES}}$ are computed by \eqref{ContravariantSurfaceFlux}. 

The numerical fluxes $\tilde{\mathbf{G}}_{\hat{n}}^{\text{ES}}$ do not generate an entropy conservative scheme, but the result in Theorem 3.3 can be used to prove that the moving mesh DGSEM becomes entropy stable, such that the discrete mathematical entropy is bounded at any time by its initial data, when the numerical fluxes $\tilde{\mathbf{G}}_{\hat{n}}^{\text{ES}}$ are used. 
\begin{Cor-hand}[3.5.]
Suppose the flux functions $\blockvec{\mathbf{G}}{}^{\text{EC}}$ in the derivative projection operator \eqref{SpatialDerivativeProjectionOperator2} are computed by Cartesian fluxes  $\mathbf{G}_{l}^{\text{EC}}$, $l=1,2,3$, with the properties \eqref{ConsistentDGECFlux}, \eqref{FluxSymmetric}, \eqref{EntropyConditionDGECFlux} and the numerical surface fluxes $\tilde{\mathbf{G}}_{\hat{n}}^{*}=\tilde{\mathbf{G}}_{\hat{n}}^{\text{ES}}$ are computed by the Cartesian fluxes  $\mathbf{G}_{l}^{\text{ES}}$, $l=1,2,3$, given by \eqref{EntropyStableSFFlux}. Then the semi-discrete moving mesh DGSEM \eqref{MovingMeshDGSEM} satisfies the discrete entropy inequality 
\begin{equation}\label{DGSEM:EntropyStability0}
\frac{d\bar{S}}{dt}\leq\underset{\text{faces}}{\sum_{\text{Boundary}}}\ \int\limits _{\partial E,N}\left(\tilde{\Psi}_{\hat{n}}-\tilde{\nu}_{\hat{n}}^{*}\Phi-\mathbf{W}^{T}\tilde{\mathbf{G}}_{\hat{n}}^{\text{ES}}\right)\,dS.
\end{equation}
\end{Cor-hand}
\begin{proof}
We proceed as in the proof of Theorem 3.3 and obtain the equation   
\begin{align}\label{DGSEM:EntropyStability1}
\begin{split}
\frac{d\bar{S}}{dt}
=& \quad \underset{\text{faces}}{\sum_{\text{Boundary}}}\ \int\limits _{\partial E,N}\left(\tilde{\Psi}_{\hat{n}}-\tilde{\nu}_{\hat{n}}^{*}\Phi-\mathbf{W}^{T}\tilde{\mathbf{G}}_{\hat{n}}^{\text{ES}}\right)\,dS \\
& -\underset{\text{faces}}{\sum_{\text{Interior}}}\int\limits _{\partial E,N}\left(\jump{\tilde{\Psi}_{\hat{n}}}-\avg{\tilde{\nu}_{\hat{n}}}\jump{\Phi}-\jump{\mathbf{W}}^{T}\tilde{\mathbf{G}}_{\hat{n}}^{\text{ES}}\right)\,dS.
\end{split}
\end{align}
Since the numerical surface fluxes $\tilde{\mathbf{G}}_{\hat{n}}^{\text{ES}}$ are computed by the Cartesian fluxes \eqref{EntropyStableSFFlux} and the fluxes $\mathbf{G}_{l}^{\text{EC}}$, $l=1,2,3$, satisfy \eqref{EntropyConditionDGECFlux}, it follows  
\begin{align}\label{DGSEM:EntropyStability2}
\begin{split}
&\quad\jump{\tilde{\Psi}_{\hat{n}}}-\jump{\Phi}\avg{\tilde{\nu}_{\hat{n}}}-\jump{\mathbf{W}}^{T}\tilde{\mathbf{G}}_{\hat{n}}^{\text{ES}} \\
=&\quad\sum_{l=1}^{3}\hat{s}n_{l}\left(\jump{\Psi_{l}}-\jump{\Phi}\avg{\nu_{l}}-\jump{\mathbf{W}}^{T}\mathbf{G}_{l}^{\text{EC}}+\frac{1}{2}\jump{\mathbf{W}}^{T}\matx{H}_{l}\jump{\mathbf{W}}\right) \\
=&\quad\frac{1}{2}\sum_{l=1}^{3}\hat{s}n_{l}\jump{\mathbf{W}}^{T}\matx{H}_{l}\jump{\mathbf{W}}.
\end{split}
\end{align}
Since the matrices $\matx{H}_{l}$, $l=1,2,3$, are symmetric positive definite and the outward normal vectors of the curved elements are positive oriented, the equation \eqref{DGSEM:EntropyStability2} provides 
\begin{equation}\label{DGSEM:EntropyStability3}
-\underset{\text{faces}}{\sum_{\text{Interior}}}\int\limits _{\partial E,N}\left(\jump{\tilde{\Psi}_{\hat{n}}}-\avg{\tilde{\nu}_{\hat{n}}}\jump{\Phi}-\jump{\mathbf{W}}^{T}\tilde{\mathbf{G}}_{\hat{n}}^{\text{ES}}\right)\,dS\leq0.
\end{equation}
Hence, we obtain the inequality  
\begin{equation}\label{DGSEM:EntropyStability4}
\frac{d\bar{S}}{dt}\leq\underset{\text{faces}}{\sum_{\text{Boundary}}}\ \int\limits _{\partial E,N}\left(\tilde{\Psi}_{\hat{n}}-\tilde{\nu}_{\hat{n}}^{*}\Phi-\mathbf{W}^{T}\tilde{\mathbf{G}}_{\hat{n}}^{\text{ES}}\right)\,dS.
\end{equation}
\end{proof}



\subsection{Free stream preservation of the moving mesh DGSEM} \label{sec:FreestreamPreservationDGSEM}
In this section, we check the discretization of the geometric and metric terms in time. Since DG methods with the forward Euler discretization are unstable \cite{Chavent1989, DiPietro2011}, we investigate directly the discretization by an explicit RK method with $s\geq2$ stages and the characteristic coefficients $\left\{ a_{\tau \sigma}\right\} _{\tau,\sigma=1}^{s}$, $\left\{ b_{\sigma}\right\} _{\sigma=1}^{s}$, $\left\{ c_{\sigma}\right\} _{\sigma=1}^{s}$. It is worth to mention that a Courant Lewy Friedrichs (CFL) restriction is necessary when an explicit $s$-stage RK method is used in the DG framework. 
In order to present the RK discretization of the semi-discrete DGSEM \eqref{MovingMeshDGSEM}, it is beneficial to write the method in the equivalent nodal representation. This representation is for all $i,j,k=0,\dots,N$, given by 
\begin{subequations}\label{MovingMeshDGSEM_Nodal}
\begin{align}
\frac{d\mathrm{J}_{ijk}}{dt}\quad &=\quad\mathsf{V}\left(\left(\vec{\nu}\right)_{ijk}\right), \label{DGSEM_Nodal:D-GLC}\\
\frac{d\left(\mathrm{J_{ijk}}\mathbf{U}_{ijk}\right)}{dt}\quad &=\quad\mathsf{\textbf{G}}\left(\left(\vec{\nu}\right)_{ijk},\mathbf{U}_{ijk}\right), \label{DGSEM_Nodal:CL}
\end{align}
\end{subequations}
where the right hand sides are given by 
\begin{equation}\label{SpatialNodalPart1}
\mathsf{V}\left(\left(\vec{\nu}\right)_{ijk}\right):=\Dprojection{N}\cdot\vec{\tilde{\nu}}_{ijk}+\frac{1}{\omega_{i}\omega_{j}\omega_{k}}\int\limits _{\partial E,N}\ell_{i}\ell_{j}\ell_{k}\left(\tilde{\nu}_{\hat{n}}^{*}-\tilde{\nu}_{\hat{n}}\right)\,dS,
\end{equation}
\begin{equation}\label{SpatialNodalPart2}
\mathsf{\textbf{G}}\left(\left(\vec{\nu}\right)_{ijk},\mathbf{U}_{ijk}\right):=-\Dprojection{N}\cdot\blockvec{\tilde{\mathbf{G}}}{}_{ijk}^{\text{EC}}-\frac{1}{\omega_{i}\omega_{j}\omega_{k}}\int\limits _{\partial E,N}\ell_{i}\ell_{j}\ell_{k}\left(\tilde{\mathbf{G}}_{\hat{n}}^{*}-\tilde{\mathbf{G}}_{\hat{n}}\right)\,dS
\end{equation}
with the tensorial Lagrange polynomials $\ell_{i}\ell_{j}\ell_{k}$ given by \eqref{LagrangeBasis}. 

Next, as in the Section \ref{sec:FreestreamPreservationFV}, the interval $\left[0,T\right]$ is divided in time levels $t^{n}$. The step size of the time discretization is $\Delta t$. The DGSEM solutions, the fluxes and the grid velocity field are approximated in the time levels $t^{n}$, e.g. $\mathbf{U}\left(t^{n}\right)\approx\mathbf{U}^{n}$. Then, the RK discretization of the semi-discrete DGSEM is similar to the RK discretization of the semi-discrete moving mesh FV scheme \eqref{FVScheme}. Thus, the moving mesh split form RK-DGSEM is given by    
\begin{subequations}\label{RK_MovingMeshDGSEM_Nodal}
\begin{align}
& \text{for $\tau=1,\dots,s$:} \nonumber \\ \label{RK_DGSEM_Nodal:D-GLC_A}
& \mathrm{J}_{ijk}^{\left(\tau\right)}=\mathrm{J}_{ijk}^{n}+\Delta t\sum_{\sigma=1}^{\tau-1}a_{\tau\sigma}\mathsf{V}\left(\left(\vec{\nu}\right)_{ijk}^{n+\sigma}\right),  \\   \label{DGSEM_Nodal:D-CL_A}  
& \mathbf{U}_{ijk}^{\left(\tau\right)}=\mathbf{U}_{ijk}^{n}+\frac{\Delta t}{\mathrm{J}_{ijk}^{\left(\tau\right)}}\sum_{\sigma=1}^{\tau-1}a_{\tau\sigma}\,\left(\mathsf{\textbf{G}}\left(\left(\vec{\nu}\right)_{ijk}^{n+\sigma},\mathbf{U}_{ijk}^{(\sigma)}\right)-\mathsf{V}\left(\left(\vec{\nu}\right)_{ijk}^{n+\sigma}\right)\mathbf{U}_{ij}^{n}\right), \\
\nonumber \\
&\mathrm{J}_{ijk}^{n+1}=\mathrm{J}_{ijk}^{n}+\Delta t\sum_{\sigma=1}^{s}b_{\sigma}\mathsf{V}\left(\left(\vec{\nu}\right)_{ijk}^{n+\sigma}\right), \label{DGSEM_Nodal:D-GLC_B}  \\
&\mathbf{U}_{ijk}^{n+1}=\mathbf{U}_{ijk}^{n}+\frac{\Delta t}{\mathrm{J}_{ijk}^{\left(n+1\right)}}\sum_{\sigma=1}^{s}b_{\sigma}\,\left(\mathsf{\textbf{G}}\left(\left(\vec{\nu}\right)_{ijk}^{n+\sigma},\mathbf{U}_{ijk}^{(\sigma)}\right)-\mathsf{V}\left(\left(\vec{\nu}\right)_{ijk}^{n+\sigma}\right)\mathbf{U}_{ijk}^{n}\right), \label{RK_DGSEM_Nodal:D-CL_B} 
\end{align}
\end{subequations}
where $\left(\vec{\nu}\right)_{ijk}^{n+\sigma}:=\vec{\nu}\left(\xi_{i}^{1},\xi_{j}^{2},\xi_{k}^{3},t^{n}+c_{\sigma}\Delta t\right)$ and $\left\{ \xi_{i}^{1}\right\} _{i=0}^{N}$, $\left\{\xi_{i}^{2}\right\} _{i=0}^{N}$, $\left\{\xi_{i}^{3}\right\} _{i=0}^{N}$ are sets of LGL points. Next, we prove that the fully-discrete split form RK-DGSEM \eqref{RK_MovingMeshDGSEM_Nodal} satisfies the free stream preservation property.   
\begin{thm-hand}[3.6.]
Suppose the solution of the fully-discrete split form RK-DGSEM \eqref{RK_MovingMeshDGSEM_Nodal} is given by $\mathbf{U}_{ijk}^{n}=\mathbf{C}:=\left(c_{1},\dots,c_{p}\right)^{T}\in\R^{p}$ for all elements $e_{\kappa}(t^{n})$, $\kappa=1,\dots,K$, and the numerical fluxes satisfy \eqref{ConsistentDGECFlux}. Then, the constant states $c_{l}$, $l=1,\dots,p$, are preserved in each Runge-Kutta stage \eqref{DGSEM_Nodal:D-CL_A}. In particular, the solution of the fully-discrete DGSEM method at time level $t^{n+1}$ is $\mathbf{U}_{ijk}^{n+1}=\mathbf{C}$.
\end{thm-hand}
\begin{proof}
Let $\tau\in\left\{1,\dots,s\right\}$ be an arbitrary fixed index. We are interested to investigate the $\tau$-th RK stage. Hence, without loss of generality, we can assume that $\textbf{U}^{(\sigma)}=\textbf{C}$ for all $\sigma=0,\dots,\tau-1$. Then, since the flux $\blockvec{\mathbf{G}}{}^{\text{EC}}$ satisfies \eqref{ConsistentDGECFlux}, it follows  
\begin{align}\label{FreestreamPreDGSEM1}
\begin{split}
\Dprojection{N}\cdot\blockvec{\tilde{\mathbf{G}}}{}_{ijk}^{\text{EC}} =&
\quad 2\sum_{m=0}^{N}\mathcal{D}_{im}\avg{\mathrm{J}\vec{a}^{1}}_{\left(i,m\right)jk}\cdot\blockvec{\mathbf{F}}\left(\mathbf{C}\right) \\
&  + 2\sum_{m=0}^{N}\mathcal{D}_{jm}\avg{\mathrm{J}\vec{a}^{2}}_{i\left(j,m\right)k}\cdot\blockvec{\mathbf{F}}\left(\mathbf{C}\right) \\
&  +2\sum_{m=0}^{N}\mathcal{D}_{km}\avg{\mathrm{J}\vec{a}^{3}}_{ij\left(k,m\right)}\cdot\blockvec{\mathbf{F}}\left(\mathbf{C}\right)-\Dprojection{N}\cdot\left(\vec{\tilde{\nu}}\right)_{_{ijk}}^{n+\sigma}\mathbf{C}.
\end{split}
\end{align}
Furthermore, since the metric terms are computed by the conservative curl form      \eqref{DiscreteContravariantVectors}, we obtain 
\begin{align}\label{FreestreamPreDGSEM1a}
\begin{split}
& 2\sum_{m=0}^{N}\left(\mathcal{D}_{im}\avg{\mathrm{J}\vec{a}^{1}}_{\left(i,m\right)jk}+\mathcal{D}_{jm}\avg{\mathrm{J}\vec{a}^{2}}_{i\left(j,m\right)k}+\mathcal{D}_{km}\avg{\mathrm{J}\vec{a}^{3}}_{ij\left(k,m\right)}\right)  \\
=&\ \ \sum_{m=0}^{N}\left(\mathcal{D}_{im}\left(\mathrm{J}\vec{a}^{1}\right)_{mjk}+\mathcal{D}_{jm}\left(\mathrm{J}\vec{a}^{2}\right)_{imk}+\mathcal{D}_{km}\left(\mathrm{J}\vec{a}^{3}\right)_{ijm}\right)=0.
\end{split}
\end{align}
Here we used the split form Lemma from Gassner et al. \cite[Lemma 1]{Gassner2016} in the first step and in the second step we used the identity \eqref{DiscreteMetricIdentities} for the discrete metric identities. Thus, it follows that  
\begin{equation}\label{FreestreamPreDGSEM1b}
\Dprojection{N}\cdot\blockvec{\tilde{\mathbf{G}}}{}_{ijk}^{\text{EC}}
=-\Dprojection{N}\cdot\left(\vec{\tilde{\nu}}\right)_{_{ijk}}^{n+\sigma}\mathbf{C}.
\end{equation}
Similar, since the flux $\blockvec{\mathbf{G}}{}^{*}$ satisfies \eqref{ConsistentDGECFlux}, it follows that   
\begin{align}\label{FreestreamPreDGSEM2}
\begin{split}
\tilde{\mathbf{G}}_{\hat{n}}^{*}-\tilde{\mathbf{G}}_{\hat{n}} 
=& \quad
\sum_{l=1}^{3}\hat{s}n_{l}\left(\mathbf{F}_{l}\left(\mathbf{C}\right)-\avg{\nu_{l}^{n+\sigma}}\mathbf{C}\right)-\left(\hat{s}\vec{n}\right)\cdot\left(\blockvec{\mathbf{F}}\left(\mathbf{C}\right)-\left(\vec{\nu}\right)^{n+\sigma}\mathbf{C}\right) \\
=& -\left(\hat{s}\vec{n}\right)\cdot\left(\avg{\left(\vec{\nu}\right)^{n+\sigma}}-\left(\vec{\nu}\right)^{n+\sigma}\right)\mathbf{C} 
=-\left(\tilde{\nu}_{\hat{n}}^{*,n+\sigma}-\tilde{\nu}_{\hat{n}}^{n+\sigma}\right)\mathbf{C}.
\end{split}
\end{align}
Thus, the equations \eqref{FreestreamPreDGSEM1b} and \eqref{FreestreamPreDGSEM2} give
\begin{equation}
\mathsf{\textbf{G}}\left(\left(\vec{\nu}\right)_{ijk}^{n+\sigma},\mathbf{C}\right)=\left(\Dprojection{N}\cdot\left(\vec{\tilde{\nu}}\right)_{_{ijk}}^{n+\sigma}+\frac{1}{\omega_{i}\omega_{j}\omega_{k}}\int\limits _{\partial E,N}\ell_{i}\ell_{j}\ell_{k}\left(\tilde{\nu}_{\hat{n}}^{*,n+\sigma}-\tilde{\nu}_{\hat{n}}^{n+\sigma}\right)\,dS\right)\mathbf{C}=\mathsf{V}\left(\left(\vec{\nu}\right)_{ijk}^{n+\sigma}\right)\mathbf{C}.
\end{equation}
Hence, the solution of the RK stage \eqref{DGSEM_Nodal:D-CL_A} is given by  
\begin{equation}
\mathbf{U}_{ijk}^{\left(\tau\right)}=\mathbf{C}+\frac{\Delta t}{\mathrm{J}_{ijk}^{\left(\tau\right)}}\sum_{\sigma=1}^{\tau-1}a_{\tau\sigma}\,\left(\mathsf{\textbf{G}}\left(\left(\vec{\nu}\right)_{ijk}^{n+\sigma},\mathbf{C}\right)-\mathsf{V}\left(\left(\vec{\nu}\right)_{ijk}^{n+\sigma}\right)\mathbf{C}\right)=\mathbf{C}.
\end{equation}
Since the parameter $\tau$ was arbitrary chosen, it follows $\mathbf{U}_{ijk}^{\left(\tau\right)}$ for all $\tau=1,\dots,s$. In particular, it follows    
\begin{equation}
\mathbf{U}_{ijk}^{n+1}=\mathbf{C}+\frac{\Delta t}{\mathrm{J}_{ijk}^{n+1}}\sum_{\sigma=1}^{s}b_{\sigma}\,\left(\mathsf{\textbf{G}}\left(\left(\vec{\nu}\right)_{ijk}^{n+\sigma},\mathbf{C}\right)-\mathsf{V}\left(\left(\vec{\nu}\right)_{ijk}^{n+\sigma}\right)\mathbf{C}\right)=\mathbf{C}.
\end{equation}
This completes the proof of Theorem 3.6.
\end{proof}

\section{Numerical results}\label{sec:numResults}
The numerical computations in this section are performed with the open source code \texttt{FLEXI}\footnote{\url{www.flexi-project.org}} and the three-dimensional high-order meshes for the simulations are generated with the open source tool \texttt{HOPR}\footnote{\url{www.hopr-project.org}}.

\begin{figure}[h]
\begin{center}
\includegraphics[width=3.2in]{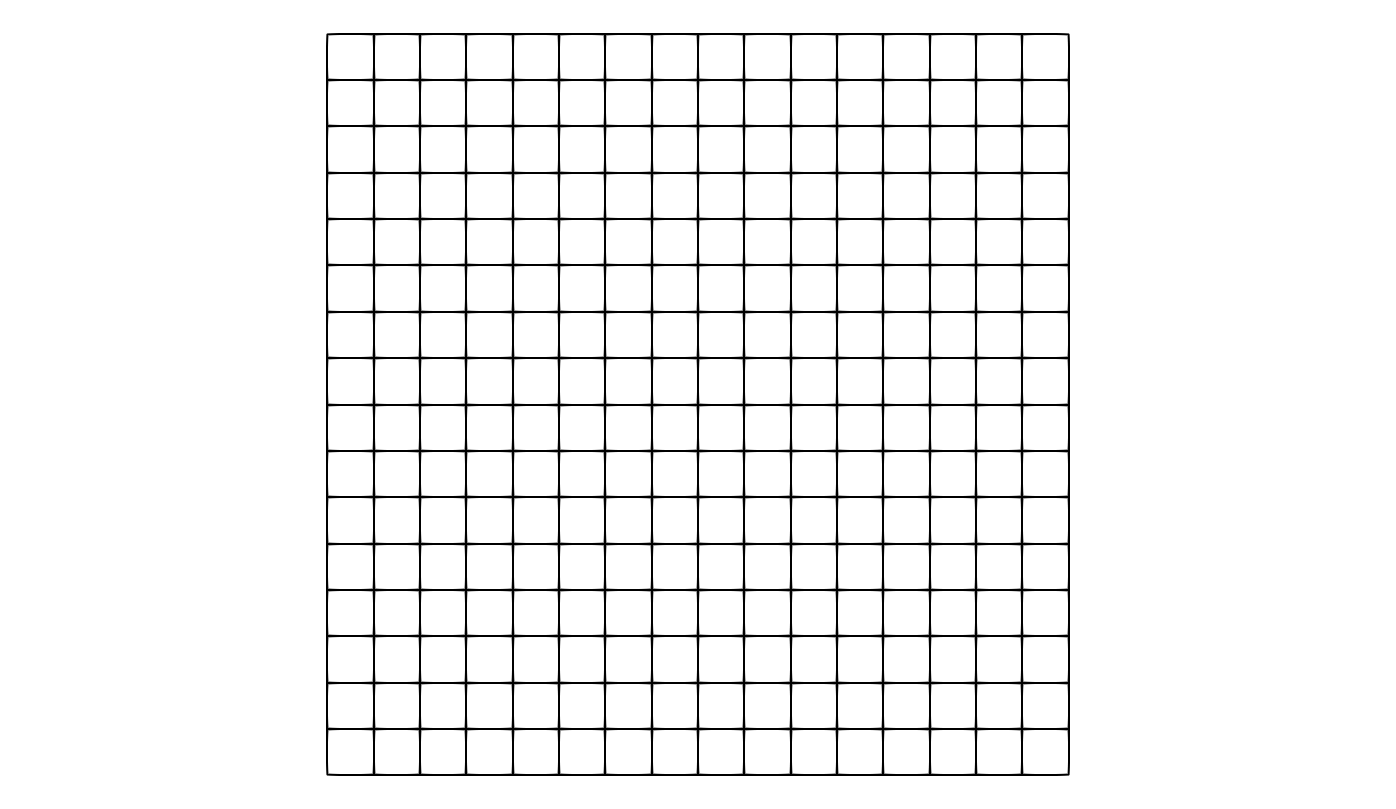}
\includegraphics[width=3.2in]{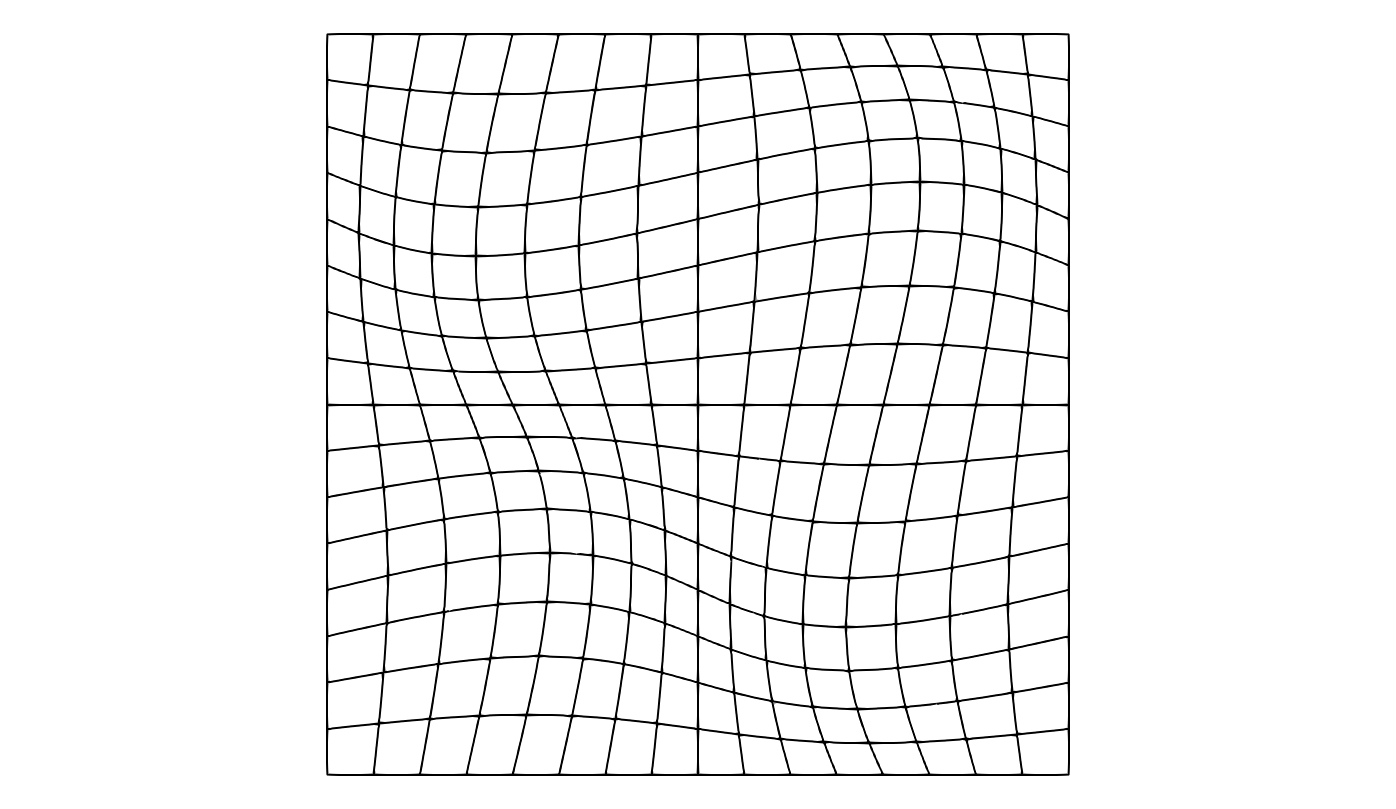}
\caption{A slice through a three dimensional mesh with $K=16^3$ elements at initial time (left) and at its maximal distortion (right).}
\label{figuremesh}
\end{center}
\end{figure}

We present tests on three dimensional moving hexahedral curved meshes for the compressible Euler equations. Based on these tests we will evaluate the theoretical findings of the previous sections. The three dimensional compressible Euler equations are given by  
\begin{equation}\label{eq:Euler}
\pderivative{\textbf{u}}{t} + \vec{\nabla} \cdot \blockvec{\textbf{f}}=\textbf{0}. 
\end{equation}
The state vector and the components of the block vector flux, $\blockvec{\textbf{f}}$, are given by
\begin{equation}
\mathbf{u}=\begin{bmatrix}\rho\\
\rho u_{1}\\
\rho u_{2}\\
\rho u_{3}\\
E
\end{bmatrix},\qquad\textbf{f}_{1}=\begin{bmatrix}\rho u_{1}\\
\rho u_{1}^{2}+p\\
\rho u_{1}u_{2}\\
\rho u_{1}u_{3}\\
\left(E+p\right)u_{1}
\end{bmatrix},\qquad\textbf{f}_{2}=\begin{bmatrix}\rho u_{2}\\
\rho u_{1}u_{2}\\
\rho u_{2}^{2}+p\\
\rho u_{2}u_{3}\\
\left(E+p\right)u_{2}
\end{bmatrix},\qquad\textbf{f}_{3}=\begin{bmatrix}\rho u_{3}\\
\rho u_{1}u_{3}\\
\rho u_{2}u_{3}\\
\rho u_{3}^{2}+p\\
\left(E+p\right)u_{3}
\end{bmatrix},
\end{equation}
where the conserved variables are the density $\rho$, the momentum $\rho\vec{u}=\left[\rho u_{1},\rho u_{2},\rho u_{3}\right]^{T}$ and the total energy $E$. In order to close the system, we assume an ideal gas such that the pressure is defined as
\begin{equation}
p = (\gamma-1)\left(E - \frac{\rho}{2}\left|\vec{u}\right|^2\right),
\end{equation}  
where $\gamma$ is the adiabatic exponent. We choose $\gamma =1.4$ in the following experiments. The system \eqref{eq:Euler} is investigated in the domain $\Omega=\left[x_{\text{min}},x_{\text{max}}\right]^3$. At initial time $t=0$ the domain is divided in $K$ non-overlapping and conforming cartesian hexahedral elements $e_{\kappa}(0)$, $\kappa=1,\dots,K$. For each element $e_{\kappa}(0)$, $\kappa = 1,\dots,K$, the temporal distribution of a grid point
\begin{equation}
\vec{x}_{\kappa}\left(0\right)=\left(x_{1}^{\kappa}\left(0\right),x_{2}^{\kappa}\left(0\right),x_{3}^{\kappa}\left(0\right)\right)^{T}\in e_{\kappa}(0)
\end{equation}  
is given by   
\begin{equation}\label{GridpointDistribution3D}
\vec{x}_{\kappa}\left(t\right)=\vec{x}_{\kappa}\left(0\right)+0.05 L \sin\left(2\pi t\right)\sin\left(\frac{2\pi}{L}x_{1}^{\kappa}\left(0\right)\right)\sin\left(\frac{2\pi}{L}x_{2}^{\kappa}\left(0\right)\right)\sin\left(\frac{2\pi}{L}x_{3}^{\kappa}\left(0\right)\right),
\end{equation}
where $L:=x_{\text{max}}-x_{\text{min}}$. In Figure~\ref{figuremesh}, we show a slice through a three dimensional mesh with $K=16^3$ elements at initial time and at its maximal distortion. 

Furthermore, the five stage fourth order low-storage explicit RK method from Kennedy, Carpenter and Lewis \cite{Kennedy2000} is used for the  time-integration in the numerical experiments. The CFL restriction is computed as in \cite{Cockburn2001} 
\begin{equation}\label{CFL1}
\frac{\Delta t}{\underset{1\leq\kappa\leq K}{\min}\left|h_{\kappa}\left(t^{n}\right)\right|}\leq\frac{C_{\text{CFL}}}{\left(2N+1\right)\lambda_{\text{max}}},
\end{equation}
where $h_{\kappa}\left(t^{n}\right)$ is the minimum element size of $e_{\kappa}(t^{n})$, $C_{\text{CFL}}\in\left(0,1\right]$ and $\lambda_{\text{max}}$ is the largest advective wave speed at the current time level traveling in either the $x_{1},x_{2},x_{3}$-direction.

\subsection{Experimental convergence rates}\label{sec:convergence3D}
In this section, we verify the high-order approximation of the moving DGSEM \eqref{MovingMeshDGSEM}. For this purpose, we investigate the domain $\Omega=\left[-1,1\right]^3$ and apply the method of manufactured solutions. Thus, we assume a solution of the form
\begin{equation}\label{MFS:Euler3D}
\mathbf{U}\left(\vec{x},t\right)=\begin{bmatrix}\rho\left(\vec{x},t\right)\\
\rho u_{1}\left(\vec{x},t\right)\\
\rho u_{2}\left(\vec{x},t\right)\\
\rho u_{3}\left(\vec{x},t\right)\\
E\left(\vec{x},t\right)
\end{bmatrix}=\begin{bmatrix}2+0.1\sin\left(\pi(x_1+x_2+x_3 - 2 \cdot 0.3t\right))\\
2+0.1\sin\left(\pi(x_1+x_2+x_3 - 2 \cdot 0.3t)\right)\\
2+0.1\sin\left(\pi(x_1+x_2+x_3 - 2 \cdot 0.3t)\right)\\
2+0.1\sin\left(\pi(x_1+x_2+x_3 - 2 \cdot 0.3t)\right)\\
\left[2+0.1\sin\left(\pi(x_1+x_2+x_3 - 2 \cdot 0.3t)\right)\right]^{2}
\end{bmatrix}.
\end{equation}

We plug solution \eqref{MFS:Euler3D} into the Euler system and compute the residual using a computer algebra system. This term is used as a source term in our convergence tests. We note that this term is handled and discretized as a solution independent part in the numerical computation.

We run the convergence test with periodic boundary conditions. Furthermore, the moving mesh DGSEM \eqref{MovingMeshDGSEM} is applied with the flux function in Appendix \ref{Chandrashekar_Flux} as volume and surface flux. In addition, the surface flux is stabilized by the dissipation operator in Appendix \ref{EulerDiss}. Besides using the grid point distribution given in \eqref{GridpointDistribution3D}, we also compute static reference solutions, by setting the grid velocity to zero. In this case, the moving mesh DGSEM \eqref{MovingMeshDGSEM} degenerates to the split form DGSEM for static meshes \cite{Gassner2016, Gassner2017}.

In Table \ref{tab:Euler3DEOCN=3}, we listed the experimental order of convergence (EOC) and $\mathrm{L}^{2}$ errors for the conserved variables that we obtain for polynomials with odd degree $N=3$ on a static mesh (top) and on a  moving mesh (bottom). The convergence rates on the moving mesh are not as good as on a static mesh, which can be justified by the high distortion in the mesh from the grid point distribution formula \eqref{GridpointDistribution3D}. However, with an increasing number of elements the same convergence rates as on a static mesh are almost reached. Moreover, the experimental order of convergence (EOC) and $\mathrm{L}^{2}$ errors for the conserved variables that we obtain for polynomials with even degree $N=4$ are listed in the Table \ref{tab:Euler3DEOCN=4}. We observe a similar behavior as for the odd degree $N=3$. This indicate the high-order approximation properties of the moving mesh DGSEM.   

\begin{center}
{ \scriptsize 
\begin{tabular}{c|cccccccccc}\toprule[1.5pt]
   $K$ &
  ${L}^{2}\left(\rho\right)$ &  $EOC(\rho)$ &  ${L}^{2}\left(\rho u_{1} \right)$  &   $EOC(\rho u_{1})$  & ${L}^{2}\left(\rho u_{2} \right)$  &   $EOC(\rho u_{2})$  & ${L}^{2}\left(\rho u_{3} \right)$  &   $EOC(\rho u_{3})$  & ${L}^{2}\left(E\right)$   & $EOC(E)$  \\\midrule
$2^3$   & 2.84E-02       & -     & 2.74E-02     & - &  2.74E-02   & - & 2.74E-02   &  -   & 5.47E-02  & -    \\
$4^3$   & 5.54e-03     & 2.36  	& 5.43E-03   & 2.34  	& 5.43E-03 	& 2.34  & 5.43E-03 & 2.34 & 1.03E-03 & 2.40    \\   
$8^3$  & 4.35E-05   &  6.99   & 4.28E-05      & 6.99  	& 4.28E-05 	 & 6.99  & 4.28E-05 & 6.99  & 1.06E-04  & 6.61      \\
$16^3$ & 2.10E-06    &  4.37  	& 2.07E-06   &  4.37 	& 2.08E-06 			& 4.37   & 2.07E-06 & 4.37  & 5.33E-06 & 4.31   \\
$32^3$  & 1.26E-07     & 4.06   & 1.24E-07 	& 4.06	  & 1.24E-07	& 4.06  & 1.24E-07 & 4.06  & 3.19E-07 & 4.06   \\ 
$64^3$ &  7.82E-09    & 4.01  & 7.67E-09 &	4.01   &  7.67E-09 	& 4.01 & 7.67E-09 & 4.01   & 1.97E-08 & 4.01   \\
\midrule
$2^3$   & 4.16E-02       & -     & 3.73E-02    & - &  3.73E-02   & - & 3.73E-02  &  -   & 5.61E-02  & -    \\
$4^3$   & 3.77E-03     & 3.46  	& 3.52E-03  & 3.41   	& 3.52E-03  	& 3.41   & 3.52E-03  & 3.41 & 6.06E-03 & 3.21     \\   
$8^3$  & 1.99E-04    & 4.25     & 1.75E-04       & 4.33   	& 1.75E-04  	 & 4.33  & 1.75E-04  &  4.33  & 3.24E-04   & 4.23       \\
$16^3$ &  5.37E-06   & 5.21    	& 4.91E-06    & 5.16  	& 4.91E-06 			& 5.16   & 4.91E-06  & 5.16  &  1.20E-05 & 4.75    \\
$32^3$  & 2.18E-07     & 4.62    & 2.07E-07	&  4.57 &  2.07E-07 	& 4.57  & 2.07E-07   & 4.57    &	5.83E-07 & 4.36   \\ 
$64^3$ & 1.45E-08     & 3.92  & 1.34E-08  &	3.95   & 1.34E-08  	& 3.95  & 1.34E-08  & 3.95   & 3.95E-08  & 3.88   \\ 
\bottomrule[1.5pt]
\end {tabular}}
\par
\captionof{table}{Experimental order of convergence (EOC) and $\mathrm{L}^{2}$ errors at time $T=5$ for the Euler manufactured solution test \eqref{MFS:Euler3D}. The moving mesh DGSEM is used with $N=3$ on a static mesh (top) and on a moving mesh (bottom) with the grid point distribution \eqref{GridpointDistribution3D}.}\label{tab:Euler3DEOCN=3}
\end{center}      
\begin{center}
{ \scriptsize 
\begin{tabular}{c|cccccccccc}\toprule[1.5pt]
   $K$ &
  ${L}^{2}\left(\rho\right)$ &  $EOC(\rho)$ &  ${L}^{2}\left(\rho u_{1} \right)$  &   $EOC(\rho u_{1})$  & ${L}^{2}\left(\rho u_{2} \right)$  &   $EOC(\rho u_{2})$  & ${L}^{2}\left(\rho u_{3} \right)$  &   $EOC(\rho u_{3})$  & ${L}^{2}\left(E\right)$   & $EOC(E)$  \\\midrule
$2^3$   & 6.99E-03       & -     &  6.64E-03     & - &  6.64E-03    & - & 6.64E-03  &  -   & 1.16E-02   & -    \\
$4^3$   &  4.02E-04     & 4.12   	&  3.97E-04   & 4.06   	&   3.97E-04  	& 4.06   & 3.97E-04  & 4.06  & 7.96E-04   &   3.87    \\   
$8^3$  &   4.50E-06   &    6.48  &    4.50E-06    &  6.47 	&   4.50E-06	 & 6.47   & 4.50E-06  &   6.47   &  1.16E-05  &  6.10       \\
$16^3$ &   1.37E-07  & 5.04     	&  1.38E-07    &  5.02	&  		1.38E-07	&    5.02 & 1.38E-07  &  5.02   & 3.66E-07  & 4.98      \\
$32^3$  &   4.33E-09  &   4.98  & 4.40E-09	&  4.97 &   4.40E-09	& 4.97  &   4.40E-09 & 4.97     &	1.16E-08 &  4.97   \\ 
$64^3$ &   1.36E-10    & 4.99   & 1.38E-10  &	4.99   &   1.38E-10	& 4.99  &  1.38E-10 & 4.99     &  3.66E-10 & 4.99    \\  
\midrule 
$2^3$   &  1.02E-02    & -     & 9.06E-03     & - & 9.06E-03     & - &  9.06E-03 &  -   & 1.45E-02  & -    \\
$4^3$   & 4.53E-04      &  4.50	& 4.13E-04  &  4.46 	&  4.13E-04  	& 4.46    & 4.13E-04    & 4.46 & 7.18E-04  & 4.33      \\   
$8^3$  & 1.10E-05     &  5.37    & 1.02E-05        & 5.35  	&   1.02E-05	 & 5.35  &  1.02E-05 &  5.35  & 1.86E-05   &  5.27      \\
$16^3$ & 1.91E-07   & 5.85     	& 1.72E-07    & 5.88  	&  	1.72E-07		&  5.88   &  1.72E-07 &  5.88 & 3.81E-07   &   5.61  \\
$32^3$  & 7.28E-09    & 4.71     & 6.33E-09	&  4.77 &  6.33E-09 	&  4.77 &  6.33E-09  &   4.77  &	1.38E-08 & 4.78   \\ 
$64^3$ & 2.79E-10     &  4.71 &  2.38E-10 &	4.74   & 2.38E-10  	& 4.74  & 2.38E-10  &  4.74  & 5.40E-10   & 4.68    \\ 
\bottomrule[1.5pt]
\end {tabular}}
\par
\captionof{table}{Experimental order of convergence (EOC) and $\mathrm{L}^{2}$ errors at time $T=5$ for the Euler manufactured solution test \eqref{MFS:Euler3D}. The moving mesh DGSEM is used with $N=4$ on a static mesh (top) and on a moving mesh (bottom) with the grid point distribution \eqref{GridpointDistribution3D}.}\label{tab:Euler3DEOCN=4}
\end{center}

\subsection{Entropy analysis validation}\label{sec:ECCheck3D}
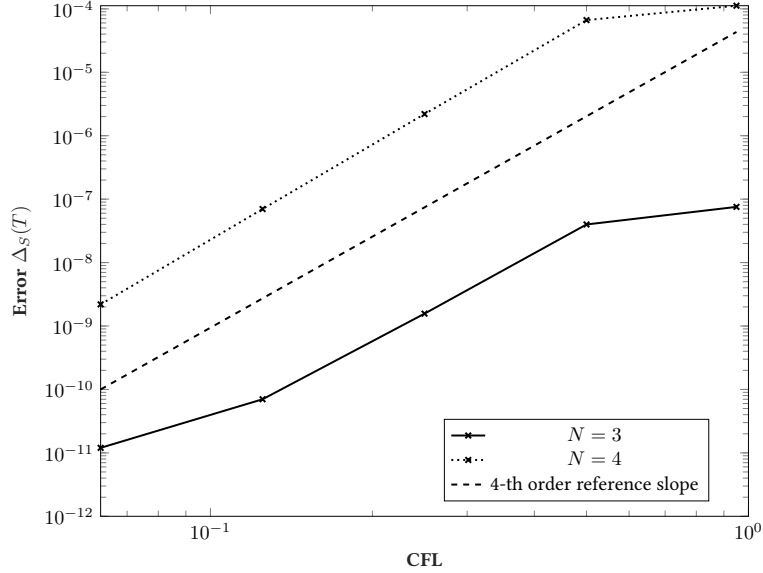
\begin{figure}[h]
\begin{center}
\begin{tikzpicture}[scale=0.75]
\begin{axis}[%
width=4.521in,
height=3.566in,
at={(0.758in,0.481in)},
scale only axis,
xmode=log,
xmin=0.0625,
xmax=1,
xminorticks=true,
xlabel style={font=\color{white!15!black}},
xlabel={\textbf{CFL}},
ymode=log,
ymin=1e-12,
ymax=0.000112133456999786,
yminorticks=true,
ylabel style={font=\color{white!15!black}},
ylabel={\textbf{Error} $\Delta_{S}(T)$ },
axis background/.style={fill=white},
legend style={at={(0.53,0.11)},anchor=west}]
\addplot [line width=1.000, color=black, mark=x, mark options={solid, black}]
  table[row sep=crcr]{%
0.0625	1.19992904501487e-11\\
0.125	6.99991176134063e-11\\
0.25	1.56699897502222e-09\\
0.5	3.98840001025746e-08\\
0.95	7.55079998526753e-08\\
};
\addlegendentry{$N=3$}

\addplot [line width=1.000, color=black,dotted , mark=x, mark options={solid, black}]
  table[row sep=crcr]{%
0.0625	2.19400142498216e-09\\
0.125	6.9891001430733e-08\\
0.25	2.19769400011671e-06\\
0.5	6.65075460002384e-05\\
0.95	0.000112133456999786\\
};
\addlegendentry{$N=4$}

\addplot [line width=1.000, color=black, dashed]
  table[row sep=crcr]{%
0.0625	1e-10\\
0.95	4.32267904852414e-05\\
};
\addlegendentry{4-th order reference slope}

\end{axis}
\end{tikzpicture}
\end{center}
\caption{Log-log plot of the entropy conservation errors $\Delta_{S}(T)$ for the Euler equations with initial data \eqref{TGVInitial}. The errors are given at time $T=13$ for polynomials with degree $N=3$ (solid line) and $N=4$ (dotted line) on a curved moving mesh with $K=16^3$ elements.}\label{fig:EntropyStability3Da}	
\end{figure}

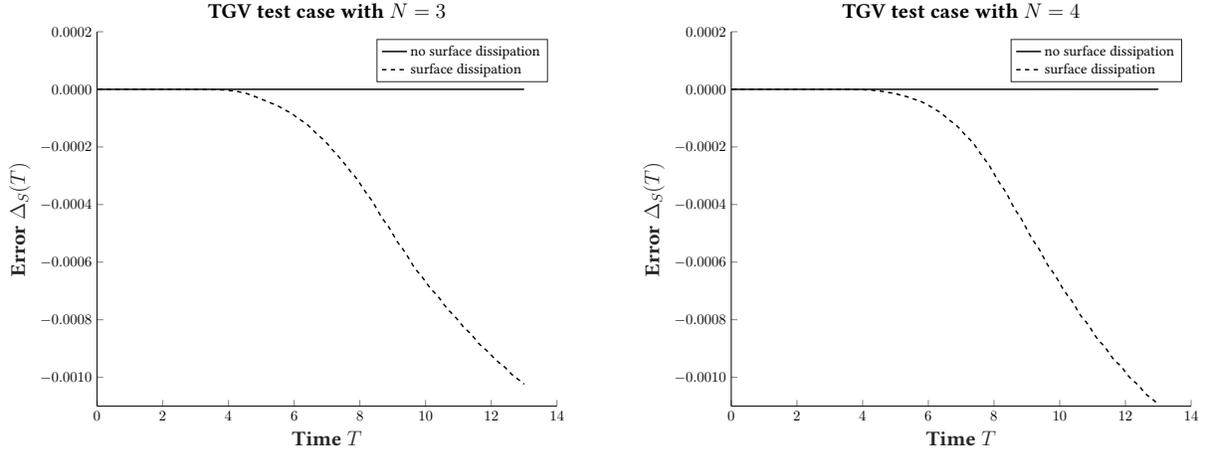
\begin{figure}[h]
\begin{center}
\begin{minipage}[t]{0.35\textwidth}
\begin{tikzpicture}[scale=0.55]
\begin{axis}[%
width=4.376in,
height=3.566in,
at={(0.903in,0.481in)},
scale only axis,
xmin=0,
xmax=14,
xlabel style={font=\color{white!15!black}},
xlabel={\Large \textbf{Time} $T$},
ymin=-0.0011,
ymax=0.0002,
yticklabel style={
            /pgf/number format/fixed,
            /pgf/number format/precision=4,
            /pgf/number format/fixed zerofill
        },
ylabel style={font=\color{white!15!black}},
ylabel={\Large \textbf{Error} $\Delta_{S}(T)$},
axis background/.style={fill=white},
title style={font=\bfseries},
title={\Large TGV test case with $N=3$},
legend style={legend cell align=left, align=left, draw=white!15!black},
    scaled y ticks=false,
    axis x line*=bottom,
    axis y line*=left
]
\addplot[black, line width=1.000,-] table {TVGRE1600N3CFL0_25NoDisEntropyError.txt};
\addlegendentry{no surface dissipation};
\addplot[black,dashed ,line width=1.000,-] table {TVGRE1600N3CFL0_25DisWintersEntropyError.txt};
\addlegendentry{surface dissipation}
\end{axis}
\end{tikzpicture}
\end{minipage}
\begin{minipage}[t]{0.15\textwidth}
\begin{tikzpicture}[scale=0.2]
\phantom{666}
\end{tikzpicture}
\end{minipage}
\begin{minipage}[t]{0.45\textwidth}
\begin{tikzpicture}[scale=0.55]
\begin{axis}[%
width=4.376in,
height=3.566in,
at={(0.903in,0.481in)},
scale only axis,
xmin=0,
xmax=14,
xlabel style={font=\color{white!15!black}},
xlabel={\Large \textbf{Time} $T$},
ymin=-0.0011,
ymax=0.0002,
yticklabel style={
            /pgf/number format/fixed,
            /pgf/number format/precision=4,
            /pgf/number format/fixed zerofill
        },
ylabel style={font=\color{white!15!black}},
ylabel={\Large \textbf{Error} $\Delta_{S}(T)$},
axis background/.style={fill=white},
title style={font=\bfseries},
title={\Large TGV test case with $N=4$},
legend style={legend cell align=left, align=left, draw=white!15!black},
    scaled y ticks=false,
    axis x line*=bottom,
    axis y line*=left
]
\addplot[black, line width=1.000,-] table {TVGRE1600N4CFL0_25NoDis_EntropyError.txt}; 
\addlegendentry{no surface dissipation};
\addplot[black,dashed ,line width=1.000,-] table {TVGRE1600N4CFL025DisWintersEntropyError.txt};
\addlegendentry{surface dissipation}
\end{axis}
\end{tikzpicture}
\end{minipage}
\end{center}		
\caption{Temporal evolution of the entropy conservation errors $\Delta_{S}(T)$ for the Euler equations with initial data \eqref{TGVInitial}. The flux in Appendix \ref{Chandrashekar_Flux} is used as surface flux without dissipation (solid line) and with the dissipation term in Appendix \ref{EulerDiss} (dashed line).}\label{fig:EntropyStability3Db}	
\end{figure}

The three dimensional Euler equations \eqref{eq:Euler} are equipped with the entropy/entropy flux pairs 
\begin{equation}\label{EulerEntropy}
s=-\frac{\rho\varsigma}{\gamma-1},\qquad f_{l}^{s}=-\frac{\rho\varsigma u_{l}}{\gamma-1},\qquad l=1,2,3,
\end{equation} 
where $\varsigma=\log\left(p\rho^{-\gamma}\right)$. We are interested in the behavior of the discrete entropy conservation error     
\begin{equation}\label{ErrorTortalEntropy3D}
\Delta_{S}(T)=\bar{S}\left(T\right)-\bar{S}\left(0\right),
\end{equation} 
where $\bar{S}\left(\cdot\right)$ is computed by \eqref{DGSEM:DisreteEntropyIntegral}. We investigate the inviscid Taylor-Green vortex (TGV) test case \cite{Shu2005} in the domain $\Omega=\left[0,2\pi\right]^{3}$. 
The inviscid TGV can be a challenging test case regarding the robustness of a numerical scheme, partly because the dynamics produce arbitrarily small scales. The flow field is thus by design under-resolved, which makes it a suitable test case to investigate the entropy conservation properties of the scheme.
The TGV evolves from the initial data    
\begin{align}\label{TGVInitial}
\begin{split}
\rho=&\quad 1, \\ 
\vec{u}=& \quad\left[\text{sin}\left(x_{1}\right)\text{cos}\left(x_{2}\right)\text{cos}\left(x_{3}\right),-\text{cos}\left(x_{1}\right)\text{sin}\left(x_{2}\right)\text{cos}\left(x_{3}\right),0\right]^{T}, \\
p=&\quad p_{0}+\frac{1}{16}\left(\text{cos}\left(2x_{1}\right)+\text{cos}\left(2x_{2}\right)\right)\left(\text{cos}\left(2x_{3}\right)+2\right).
\end{split}
\end{align}
To render the simulation close to incompressible, the Mach number $M_0=\frac{1}{\sqrt{\gamma p_{0}}}$ is set to $0.1$ by adjusting the pressure correspondingly. We run the simulation with $K=16^3$ elements and periodic boundary conditions. The final time is chosen to be $T=13$. Furthermore, we apply the flux function in Appendix \ref{Chandrashekar_Flux} to compute the derivative projection operator \eqref{SpatialDerivativeProjectionOperator2}. In Figure \ref{fig:EntropyStability3Da} we present a log-log plot of the entropy conservation error for $N=3,4$. We note that the flux in Appendix \ref{Chandrashekar_Flux} was used as surface flux without a dissipation term in these computations, rendering the semi-discrete discretization fully entropy conserving. As expected, we observe the reduction of the remaining entropy conservation error according to the order of the RK method for decreasing CFL numbers.  
In Figure \ref{fig:EntropyStability3Db} the temporal evolution of the entropy conservation errors $\Delta_{S}(T)$ is given. The CFL number is set to $C_{\text{CFL}}=0.125$ and polynomial degrees $N=3$ and $N=4$ are used. We observe that the entropy conservation error $\Delta_{S}(T)$ is constant in time (dashed line) when the flux in Appendix \ref{Chandrashekar_Flux} is used without a dissipation term. This indicates the entropy conservation in the TVG test case. On the other hand the entropy conservation error $\Delta_{S}(T)$ is decreasing in time (solid line) when the surface flux is stabilized by the dissipation term in Appendix \ref{EulerDiss}. Thus, the moving mesh DGSEM is an entropy stable scheme in this test case. These observations agree with the results in Theorem 3.3 and Corollary 3.5.   

\subsection{Free stream preservation validation}\label{sec:FreeStreamValidation}
We consider the Euler equations \eqref{eq:Euler} on the domain $\Omega=[0,2\pi]^3$ with the initial data
\begin{equation}\label{FreestreamSolution:Euler3D}
\mathbf{U}\left(x,t\right)=\begin{bmatrix}\rho\left(x,t\right)\\
\rho u_{1}\left(x,t\right)\\
\rho u_{2}\left(x,t\right)\\
\rho u_{3}\left(x,t\right)\\
E\left(x,t\right)
\end{bmatrix}=\begin{bmatrix}1 \\
0.3 \\
0 \\
0 \\
17 \\
\end{bmatrix}.
\end{equation}
The entropy stable DGSEM is applied with the flux function in Appendix \ref{Chandrashekar_Flux} as volume and surface flux as well as the  dissipation operator in Appendix \ref{EulerDiss} to stabilize the surface flux. We apply $K=16^3$ elements, the formula \eqref{GridpointDistribution3D} to describe the displacement of the mesh points and periodic boundary conditions are used in the simulation. Furthermore, the final time is set to $T=20$. In Table \ref{FreestreamSolution:Euler3D}, we present the $\mathrm{L}^{\infty}$ errors between the initial data \eqref{FreestreamSolution:Euler3D} and the numerical solution at time $T=20$ for polynomials of degree $N=3$ (top), $N=4$ (bottom) and different CFL numbers $C_{\text{CFL}}$. We observe that the errors are close to zero and vary slightly for the different CFL numbers. These results indicate the compliance of the free stream preservation property.      
\begin{center}
\begin{tabular}{c|ccccc}\toprule[1.5pt]
$C_{\text{CFL}}$ & 
$\mathrm{L}^{\infty}\left(\rho\right)$ & $\mathrm{L}^{\infty}\left(\rho u_{1}\right)$ & $\mathrm{L}^{\infty}\left(\rho u_{2}\right)$ & $\mathrm{L}^{\infty}\left(\rho u_{3}\right)$  & $\mathrm{L}^{\infty}\left(E\right)$ \\ \midrule
0.95 & 2.47E-14  & 1.40E-12 & 4.46E-12 &	4.48E-12 &	1.33E-12 	   \\
0.5 &  2.47E-14  &	1.40E-12 & 	4.46E-12 & 	4.48E-12 & 	1.33E-12 	 \\
0.25 & 2.70E-14  &	1.40E-12 &	4.46E-12 & 	4.48E-12 &	1.36E-12 \\
0.125 & 3.10E-14 &	1.40E-12 & 	4.46E-12 &	4.48E-12 &	1.43E-12 \\
0.0625 & 3.78E-14 & 1.40E-12 & 	4.46E-12 & 	4.48E-12 & 	1.56E-12 \\
 \midrule 
0.95 & 2.07E-14 & 	1.24E-12 & 	5.28E-12 & 	5.21E-12 & 	1.12E-12 \\
0.5 &  2.49E-14 & 	1.24E-12 & 	5.28E-12 & 	5.21E-12 & 	1.30E-12 \\
0.25 & 2.81E-14 &	1.24E-12 & 	5.28E-12 &	5.21E-12 & 	1.34E-12 \\
0.125 & 3.32E-14 & 	1.24E-12 & 	5.28E-12 & 	5.21E-12 & 	1.40E-12 \\
0.0625 & 4.24E-14 & 1.24E-12 & 	5.28E-12 &	5.21E-12 & 	1.59E-12  \\
\bottomrule[1.5pt]
\end{tabular}\par
\captionof{table}{Free stream preservation test for $N=3$ (top) and $N=4$ (bottom). The $\mathrm{L}^{\infty}$ errors measure the difference between the initial data \eqref{FreestreamSolution:Euler3D} and the numerical solution at time $T=20$ for different constants $C_{\text{CFL}}$.}\label{tab:FreestreamEuler3D}
\end{center}

\subsection{Robustness Test}\label{RobustnessTest}

As has been stated in Section \ref{sec:ECCheck3D} and noted in literature \cite{Moura2017,Winters2018}, the inviscid TGV is a notoriously challenging test case for the stability of a high order discretization. While for lower polynomial degrees calculations may be possible, high-order simulations are known to crash even if aliasing-reducing methods like polynomial dealiasing are used \cite{Moura2017}.
Thus, we use the TGV test case \eqref{TGVInitial} to demonstrate the increased robustness of the entropy stable moving mesh DGSEM. To do so, we run the simulation up to $T=13$ using a polynomial degree of $N=7$ on three different meshes employing $K_1=14^3$, $K_2=19^3$ and $K_3=26^3$ elements. These cases correspond to the most restrictive simulatons from \cite{Moura2017}.
Again, the point distribution given in \eqref{GridpointDistribution3D} is used. We use the flux function in Appendix \ref{Chandrashekar_Flux} as volume and surface flux and stabilize the surface flux by the dissipation operator in Appendix \ref{EulerDiss}.

Using the entropy stable moving DGSEM, we were able to run all simulations until final time. This shows that the consistent dissipation operators in combination with the entropy conservative volume fluxes can lead to superior stability properties. 

\section{Conclusions}\label{sec:conc}
In this work a moving mesh DGSEM to solve non-linear conservation laws has been constructed and analyzed. The semi-discrete method is provably entropy stable and the free stream preservation property is satisfied for each explicit $s$-stage Runge-Kutta method.      

We started the paper with a view on the ALE transformation for one dimensional conservation laws and constructed a moving mesh FV method. The condition \eqref{FV:EntropyCondition} to construct entropy conservative two-point fluxes in the ALE framework has been presented. This condition can be seen as an extension of Tadmor's discrete entropy condition \cite{Tadmor1987} in the ALE framework. It is remarkable that an entropy conservative moving mesh flux can be constructed directly from an entropy conservative static mesh flux, when a numerical state function with the property \eqref{FV:EntropyConditionState} is available. This approach has been used in the appendix to construct entropy conservative fluxes for the shallow water and Euler equations. 

Then, the moving mesh DGSEM has been presented for three dimensional conservation laws. The derivatives in space are approximated with high-order derivative matrices which are SBP operators. Furthermore, the  split form DG framework \cite{Gassner2016,Gassner2017} has been used to avoid aliasing in the discretization of the volume integrals. In addition, two-point flux functions with the generalized entropy condition \eqref{FV:EntropyCondition} are used in the split form DG framework. These modules in the spatial discretization are the basis to prove that the moving mesh DGSEM is an entropy stable scheme. It is worth to mention that the discrete entropy analysis requires merely the assumption that the time derivatives can be evaluated exactly. Operations like the integration-by-parts formula are mimicked    
by the SPB operators on the discrete level. 


The three dimensional Euler equations have been considered to verify the proven properties of the moving mesh DGSEM in our numerical experiments. We presented convergence tests for smooth test problems to verify that the split form DG framework provides also on a moving mesh a high-order accurate approximation. Furthermore, the numerical robustness tests in the Section \ref{RobustnessTest} emphasize the relevance of the entropy stable DGSEM, since the method was able to run the challenging inviscid TGV test case until final time.   

\section*{Acknowledgement} 
Gero Schn\"ucke and Gregor Gassner are supported by the European Research Council (ERC) under the European Union's Eights Framework Program Horizon 2020 with the research project Extreme, ERC grant agreement no. 714487. The authors gratefully acknowledge the support and the computing time on "Hazel Hen" provided by the HLRS through the project "hpcdg".


\begin{thebibliography}{99}

\bibitem{Barth2018} T. J. Barth. Numerical methods for gasdynamic systems on unstructured meshes. In: An introduction to recent developments in theory and numerics for conservation laws. Springer, Berlin, Heidelberg, 1999. S. 195-285.

\bibitem{Bellman} R. Bellman. Introduction to matrix analysis, volume 19 of Classics in Applied Mathematics. SIAM, Philadelphia, PA, 2nd Edition, 1987.

\bibitem{Bohm2018}  
M. Bohm, A. R. Winters, G. J. Gassner, D. Derigs, F. Hindenlang and J. Saur. An entropy stable nodal discontinuous Galerkin method for the resistive MHD equations. Part I: Theory and Numerical Verification. J. Comput. Phys. (2018),  \url{doi.org/10.1016/j.jcp.2018.06.027}.

\bibitem{Boscheri2017} W. Boscheri and M. Dumbser. Arbitrary-Lagrangian-Eulerian discontinuous Galerkin schemes with a posteriori subcell finite volume limiting on moving unstructured meshes. J. Comput. Phys. 346 (2017), 449-479.

\bibitem{Butcher1987} J. C. Butcher. The Numerical Analysis of Ordinary Differential Equations. Runge-Kutta and General Linear Methods. Wiley, Chichester, 1987.


\bibitem{chan2018} 
J. Chan. On discretely entropy conservative and entropy stable discontinuous Galerkin methods. J. Comput. Phys. 362 (2018), 346-374.

\bibitem{CHQZ:2006} 
C. Canuto, M. Y. Hussaini, A. Quarteroni and T. A. Zang. Spectral methods: Fundamentals in single domains. Springer Verlag, 2006.  

\bibitem{Chandrashekar2013} P. Chandrashekar. Kinetic energy preserving and entropy stable finite volume schemes for compressible Euler and Navier-Stokes equations. Commun. Comput. Phys. 14 (2013), 1252-1286.

\bibitem{Chandrashekar2017} P. Chandrashekar and M. Zenk. Well-balanced nodal discontinuous Galerkin method for Euler equations with gravity. J. Sci. Comput. 71 (2017), 1062-1093.



\bibitem{Chavent1989} G. Chavent and B. Cockburn. The Local Projection $P^0P^1$-Discontinuous Galerkin Finite Element Method for Scalar Conservation Laws. ESAIM-Math. Model. Num. ($M^2AN$) 23 (1989), 565-592.

\bibitem{Chen2017} T. Chen and  C.-W. Shu. Entropy stable high order discontinuous Galerkin methods with suitable quadrature rules for hyperbolic conservation laws. J. Comput. Phys. 345 (2017), 427-461.

\bibitem{Chiodaroli2015} 
E. Chiodaroli, C. De Lellis and O. Kreml. Global Ill-Posedness of the Isentropic System of Gas Dynamics. Commun. Pur. Appl. Math. 68 (2015), 1157-1190.

\bibitem{Cockburn2001} B. Cockburn and C.-W. Shu. Runge-Kutta discontinuous Galerkin methods for convection-dominated problems. J. Sci. Comput. 16 (2001), 173-261.

\bibitem{crean2018} J. Crean, J. E. Hicken, D. C. Del Rey Fern\'andez, D. W. Zingg and M. H. Carpenter. Entropy-stable summation-by-parts discretization of the Euler equations on general curved elements. J. Comput. Phys. 356 (2018), 410-438.

\bibitem{DiPerna1983} R. J. Di Perna. Convergence of approximate solutions to conservation laws. Arch. Ration. Mech. An. 82 (1983), 27-70.

\bibitem{DiPietro2011} D. A. Di Pietro and A. Ern. Mathematical aspects of discontinuous Galerkin methods. Vol. 69. Springer Science \& Business Media, 2011.

\bibitem{Donea2004} J. Donea, A. Huerta, J. P. Ponthot and A. Rodr\'iguez-Ferran. In Encyclopedia of Computational Mechanics, E. Stein, R. De Borst and T. J.R. Hughes (Eds.), Volume. 1: Fundamentals., Chapter 14: Arbitrary Lagrangian-Eulerian Methods, Wiley, 2004.

\bibitem{Farhat2001} C. Farhat, P. Geuzaine and C. Grandmont. The Discrete Geometric Conservation Law and the Nonlinear Stability of ALE Schemes for the Solution of Flow Problems on Moving Grids. J. Comput. Phys. 174 (2001), 669-694.

\bibitem{Fazio2003} R. Fazio and R. J. LeVeque. Moving-mesh methods for one-dimensional hyperbolic problems using CLAWPACK. Comput. Math. Appl. 45 (2003), 273-298.

\bibitem{Fernandez2014}
D. C. Del Rey Fern\'andez, J. E. Hicken, and D. W. Zingg. Review of summation-by-parts operators with simultaneous approximation terms for the numerical solution of partial differential equations. Comput. Fluids (2014), 171-196.

\bibitem{Fisher2013} T. C. Fisher and M. H. Carpenter. high order entropy stable finite difference schemes for nonlinear conservation laws: Finite domains. J. Comput. Phys. 252 (2013), 518-557.

\bibitem{Fjordholm2011} U. S. Fjordholm, S. Mishra and E. Tadmor. Well-balanced and energy stable schemes for the shallow water equations with discontinuous topography. J. Comput. Phys. 230 (2011), 5587-5609.

\bibitem{Friedrich2018} L. Friedrich, G. Schn\"ucke, A. R. Winters, D. C. Del Rey Fern\'andez, G. J. Gassner and M. H. Carpenter. Entropy Stable Space-Time Discontinuous Galerkin Schemes with Summation-by-Parts Property for Hyperbolic Conservation Laws. (2018) arXiv preprint arXiv:1808.08218.

\bibitem{Fu2018} P. Fu, G. Schn\"ucke, and Y. Xia. Arbitrary Lagrangian-Eulerian discontinuous Galerkin method for conservation laws on moving simplex meshes. Math. Comput. (2018). Available via \url{https://doi.org/10.1090/mcom/3417} (Cited \today) 

\bibitem{Gassner2013} G. J. Gassner. A skew-symmetric discontinuous Galerkin spectral element discretization and its relation to SBP-SAT finite difference methods. SIAM J. Sci. Comput. 35 (2013), A1233-A1253. 


\bibitem{Gassner2016} G. J. Gassner, A. R. Winters and D. A. Kopriva. Split form nodal discontinuous Galerkin schemes with summation-by-parts property for the compressible Euler equations. J. Comput. Phys. 327 (2016), 39-66.

\bibitem{Gassner2017} G. J. Gassner, A. R. Winters, F. J. Hindenlang and D. A. Kopriva. The BR1 scheme is stable for the compressible Navier-Stokes equations.  J Sci. Comput. 77 (2018), 1-47.

\bibitem{Gassner2016WB} 
G. J. Gassner, A. R. Winters and D. A. Kopriva. A well balanced and entropy conservative discontinuous Galerkin spectral element method for the shallow water equations. Appl. Math. Comput. 272 (2016), 291-308.


\bibitem{Godunov1961} S. K. Godunov. An interesting class of quasilinear systems. Dokl. Acad. Nauk SSSR. Vol. 139 (1961), 521-523. 


\bibitem{Farhat2000} H. Guillard and C. Farhat. On the significance of the geometric conservation law for flow computations on moving meshes, Comput. Method. Appl. M. 190 (2000), 1467-1482. 

\bibitem{Harten1983} A. Harten. On the symmetric form of systems of conservation laws with entropy. J. Comput. Phys. 49 (1983), 151-164. 


\bibitem{Hindenlang2015} 
F. J. Hindenlang, T. Bolemann and C.-D. Munz. Mesh curving techniques for high order discontinuous Galerkin simulations. In: IDIHOM: Industrialization of high-order methods-a top-down approach. Springer, Cham, 2015, 133-152.


\bibitem{Huang2010} W. Huang and R. D. Russell. Adaptive moving mesh methods. Vol. 174. Springer Science \& Business Media, 2010.

\bibitem{Ismail2009} F. Ismail and P. L. Roe. Affordable, entropy-consistent Euler flux functions II: Entropy production at shocks. J. Comput. Phys. 228 (2009), 5410-5436.




\bibitem{Kennedy2000} 
C. A. Kennedy, M. H. Carpenter and R. M. Lewis. Low-storage, explicit Runge-Kutta schemes for the compressible Navier-Stokes equations. Appl. Numer. Math. 35 (2000), 177-219.

\bibitem{Kopriva2006}
D. A. Kopriva. Metric identities and the discontinuous spectral element method on curvilinear meshes. J. Sci. Comout. 26 (2006), 301-327.

\bibitem{Kopriva2009} D. A. Kopriva. Implementing spectral methods for partial differential equations: Algorithms for scientists and engineers. Springer Science \& Business Media, 2009.


\bibitem{Kopriva2016} D. A. Kopriva, A. R. Winters, M. Bohm and G. J. Gassner. A provably stable discontinuous Galerkin spectral element approximation for moving hexahedral meshes. Comput. Fluids 139 (2016), 148-160. 


\bibitem{kreiss1}
H.-O. Kreiss and J. Oliger. Comparison of accurate methods for the integration of hyperbolic equations. Tellus 24  (1972), 199-215.

\bibitem{Kruzkov1970} S. N. Kru\v{z}kov. First order quasilinear equations in several independent variables. Math. USSR SB+ 10 (1970), 217-243.


\bibitem{Lesoinne1996}  M. Lesoinne and C. Farhat. Geometric conservation laws for flow problems with moving boundaries and deformable meshes, and their impact on aeroelastic computations. Comput. Method. Appl. M. 134 (1996), 71-90.

\bibitem{Lombard1979} C. K. Lombard and P. D. Thomas. Geometric conservation law and its application to flow computations on moving grids. AIAA J., 17 (1979), 1030-1037.

\bibitem{Lomtev1999JCP} I. Lomtev, R. M. Kirby and G. E. Karniadakis. A Discontinuous Galerkin ALE Method for Compressible Viscous Flows in Moving Domains. J. Comput. Phys. 155 (1999), 128-159.

\bibitem{Loubere2010} R. Loub\`ere, P. H. Maire, M. Shashkov, J. Breil and S. Galera. ReALE: A reconnection-based arbitrary-Lagrangian-Eulerian method. J. Comput. Phys. 229 (2010), 4724-4761.


\bibitem{Marinacci2013} F. Marinacci, R. Pakmor and V. Springel. The formation of disc galaxies in high-resolution moving-mesh cosmological simulations. Mon. on Not. R. Astron. Soc. 437 (2013), 1750-1775.

\bibitem{Mavriplis2006} D. J. Mavriplis and Z. Yang. Construction of the discrete geometric conservation law for high-order time-accurate simulations on dynamic meshes, J. Comput. Phys. 213 (2006), 557-573.

\bibitem{Minoli2011} C. A. A. Minoli and D. A. Kopriva. Discontinuous Galerkin spectral element approximations on moving meshes. J. Comput. Phys. 230 (2011), 1876-1902.

\bibitem{Mock1980} M. S. Mock. Systems of conservation laws of mixed type. J. Differ. Equations 37 (1980), 70-88.

\bibitem{Moura2017} R. C. Moura, G. Mengaldo, J. Peiró and S. J. Sherwin. On the eddy-resolving capability of high-order discontinuous Galerkin approaches to implicit LES / under-resolved DNS of Euler turbulence. J. Comput. Phys. 330 (2017) 615-623.

\bibitem{Nguyen2010JFS} V. T. Nguyen. An arbitrary Lagrangian–Eulerian discontinuous Galerkin method for simulations of flows over variable geometries. J. Fluid. Struct. 26 (2010), 312-329.


\bibitem{Persson2009} P. O. Persson, J. Bonet, J. Peraire. Discontinuous Galerkin Solution of the Navier-Stokes Equations on Deformable Domains. Comput. Method. Appl. M. 198 (2009), 1585-1595.

\bibitem{Ranocha2018} H. Ranocha. Generalised Summation-by-Parts Operators and Entropy Stability of Numerical Methods for Hyperbolic Balance Laws. Cuvillier Verlag, 2018.




\bibitem{Shu2005} 
C. W. Shu, W. S. Don, D. Gottlieb, O. Schilling and L. Jameson. Numerical convergence study of nearly incompressible, inviscid Taylor–Green vortex flow. J. Sci. Comput. 24 (2005), 1-27.

\bibitem{Shunn2007} L. Shunn and F. Ham. Method of manufactured solutions applied to variable-density flow solvers. Annual Research Briefs-2007, Center for Turbulence Research (2007), 155-168.

\bibitem{Springel2010} V. Springel. E pur si muove: Galilean-invariant cosmological hydrodynamical simulations on a moving mesh. Mon. on Not. R. Astron. Soc. 401 (2010), 791.851.

\bibitem{Tadmor1987} Tadmor, Eitan. The numerical viscosity of entropy stable schemes for systems of conservation laws. I. Math. Comput. 49 (1987), 91-103.

\bibitem{Tadmor2003} E. Tadmor. Entropy stability theory for difference approximations of nonlinear conservation laws and related time-dependent problems, Acta Numer. 12 (2003). 451-512. 


\bibitem{Toro2013} E. F. Toro. Riemann solvers and numerical methods for fluid dynamics: a practical introduction. Springer Science \& Business Media, 2013.

\bibitem{Persson2015}  L. Wang and P. O. Persson. High-order Discontinuous Galerkin Simulations on Moving Domains using ALE Formulations and Local Remeshing and Projections, 53rd AIAA Aerospace Sciences Meeting. AIAA SciTech Forum (AIAA 2015-0820).
Available via \url{https://dx.doi.org/10.2514/6.2015-0820} (Cited \today)

\bibitem{Wintermeyer2017} N. Wintermeyer, A. R. Winters, G. J. Gassner and D. A. Kopriva. An entropy stable nodal discontinuous Galerkin method for the two dimensional shallow water equations on unstructured curvilinear meshes with discontinuous bathymetry. J. Comput. Phys. 340 (2017), 200-242. 

\bibitem{Winters2018} A. R. Winters, R. C. Moura, G. Mengaldo, G. J. Gassner, S. Walch, J. Peiro, S. J. Sherwin. A comparative study on polynomial dealiasing and split form discontinuous Galerkin schemes for under-resolved turbulence computations. J. Comput. Phys. 372 (2018), 1-21.

\bibitem{Winters2017} 
A. R. Winters, D. Derigs, G. J. Gassner and S. Walch. A uniquely defined entropy stable matrix dissipation operator for high Mach number ideal MHD and compressible Euler simulations. J. Comput. Phys. 332 (2017), 274-289.


\bibitem{Winters2014} 
A. R. Winters and D. A. Kopriva. ALE-DGSEM approximation of wave reflection and transmission from a moving medium. J. Comput. Phys. 263 (2014), 233-267.  

\bibitem{Winters2014_A}
A. R. Winters. Discontinuous Galerkin spectral element approximations for the reflection and transmission of waves from moving material interfaces. Diss. The Florida State University, 2014.
\end{thebibliography}

\appendix

\renewcommand\thefigure{\thesection.\arabic{figure}}    
\setcounter{figure}{0}    

\section{Proof for Lemma 3.1}\label{sec:ProofGridVelocity}
In this section, we prove Lemma 3.1. For the sake of simplicity we present merely the proof in two dimensions. The three dimensional proof can be done by the same argumentation.     

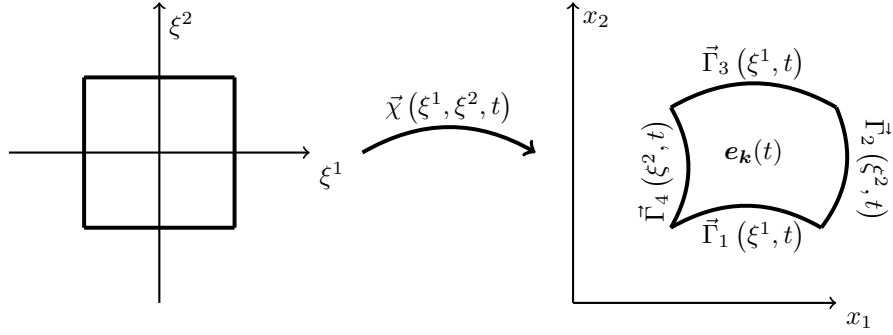
\begin{figure}[h]
\begin{center}
\begin{tikzpicture}[scale=2.0]

 \draw[thick,->] (1.75,0,0) -- (3.5,0,0) node[anchor = north west]{$x_{1}$};
 \draw[thick,->] (1.75,0,0) -- (1.75,2,0) node[anchor = north west]{$x_{2}$};

\node at (2.95,1)  {$\boldsymbol{e_{k}}(t)$};
\draw[line width=1.5pt,bend left,-] (2.4,0.5) to (3.4, 0.5);
\node at (2.95,0.45)  {$\vec{\Gamma}_{1}\left(\xi^{1},t\right)$};
\draw[line width=1.5pt,bend right,-] (3.4,0.5) to (3.5, 1.3);
\node[label=below:\rotatebox{-90}{$\vec{\Gamma}_{2}\left(\xi^{2},t\right)$}] at (3.75,1.35) {};
\draw[line width=1.5pt,bend right,-] (3.5, 1.3) to (2.4, 1.3);
\node at (2.95,1.6)  {$\vec{\Gamma}_{3}\left(\xi^{1},t\right)$};
\draw[line width=1.5pt, bend left,-] (2.4, 1.3) to (2.4,0.5);
\node[label=below:\rotatebox{90}{$\vec{\Gamma}_{4}\left(\xi^{2},t\right)$}] at (2.3,1.35) {};

\draw[line width=1.5pt,bend left,->] (0.35, 1.0) to (1.5, 1.0);
\node at (0.925,1.3)  {$\vec{\chi}\left(\xi^{1},\xi^{2},t\right)$};

\draw[thick,->] (-2,1,0) -- (0,1,0) node[anchor = north west]{$\xi^{1}$};
\draw[thick,->] (-1,0,0) -- (-1,2.0,0) node[anchor = north west]{$\xi^{2}$}; 

\draw[line width=1.5pt,left,-] (-1.5,0.5) to (-0.5, 0.5);
\draw[line width=1.5pt, right,-] (-0.5,0.5) to (-0.5, 1.5);
\draw[line width=1.5pt, right,-] (-0.5, 1.5) to (-1.5, 1.5);
\draw[line width=1.5pt,  left,-] (-1.5, 1.5) to (-1.5,0.5); 
\end{tikzpicture}
\caption{\label{ReferencePhysicalCell} Left the reference element $E=[-1,1]^2$ and on the right a general quadrilateral element $e_{k}(t)$ with the curved faces $\vec{\Gamma}_{1}\left(\xi^{1},t\right)$, $\vec{\Gamma}_{2}\left(\xi^{2},t\right)$, $\vec{\Gamma}_{3}\left(\xi^{1},t\right)$ and $\vec{\Gamma}_{4}\left(\xi^{2},t\right)$. The mapping $\vec{\chi}\left(\xi^{1},\xi^{2},t\right)$ connects $E$ and $e_{k}(t)$.}
\end{center}
\end{figure}

The two dimensional bijective isoparametric transfinite mapping can be computed as in Kopriva \cite[Chapter 6, equation (6.18)]{Kopriva2009}). The mapping is for all $\left[\xi^{1},\xi^{2}\right]^{T}\in E$  given by    
\begin{align}\label{IsoparametricTransformation}
\begin{split}
\vec{\chi}\left(\xi^{1},\xi^{2},t\right)=&\quad\frac{1}{2}\left[\left(1-\xi^{1}\right)\interpolation{N}{\left(\vec{\Gamma}_{4}\right)}\left(\xi^{2},t\right)+\left(1+\xi^{1}\right)\interpolation{N}{\left(\vec{\Gamma}_{2}\right)}\left(\xi^{2},t\right)\right]\\
&+\frac{1}{2}\left[\left(1-\xi^{2}\right)\interpolation{N}{\left(\vec{\Gamma}_{1}\right)}\left(\xi^{1},t\right)+\left(1+\xi^{2}\right)\interpolation{N}{\left(\vec{\Gamma}_{3}\right)}\left(\xi^{1},t\right)\right]\\
&-\frac{1}{4}\left[\left(1+\xi^{1}\right)\left\{ \left(1-\xi^{2}\right)\interpolation{N}{\left(\vec{\Gamma}_{1}\right)}\left(1,t\right)+\left(1+\xi^{2}\right)\interpolation{N}{\left(\vec{\Gamma}_{3}\right)}\left(1,t\right)\right\} \right]\\
&-\frac{1}{4}\left[\left(1-\xi^{1}\right)\left\{ \left(1-\xi^{2}\right)\interpolation{N}{\left(\vec{\Gamma}_{1}\right)}\left(-1,t\right)+\left(1+\xi^{2}\right)\interpolation{N}{\left(\vec{\Gamma}_{3}\right)}\left(-1,t\right)\right\} \right].
\end{split}
\end{align}
It is worth to mention that the mapping $\vec{\chi}\left(\xi^{1},\xi^{2},t\right)$ matches with the boundary faces in the interpolation points. The location of the curved faces $\vec{\Gamma}_{1}\left(\xi^{1},t\right)$, $\vec{\Gamma}_{2}\left(\xi^{2},t\right)$, $\vec{\Gamma}_{3}\left(\xi^{1},t\right)$ and $\vec{\Gamma}_{4}\left(\xi^{2},t\right)$ is sketched in Figure \ref{ReferencePhysicalCell}. 

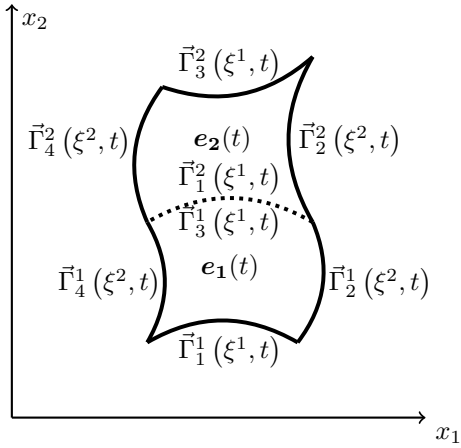
\begin{figure}[h]
\begin{center}
\begin{tikzpicture}[scale=2.0]
\node at (1.45,1)  {$\boldsymbol{e_{1}}(t)$};
 
\draw[line width=1.5pt,bend left,-] (0.9,0.5) to (1.9, 0.5);

\node at (1.45,0.45)  {$\vec{\Gamma}_{1}^{1}\left(\xi^{1},t\right)$};

\draw[line width=1.5pt,bend right,-] (1.9,0.5) to (2, 1.3);

\node at (2.45,0.9)  {$\vec{\Gamma}_{2}^{1}\left(\xi^{2},t\right)$};

\draw[line width=1.5pt,dotted,bend right,-] (2, 1.3) to (0.9, 1.3);

\node at (1.45,1.3)  {$\vec{\Gamma}_{3}^{1}\left(\xi^{1},t\right)$};

\draw[line width=1.5pt, bend left,-] (0.9, 1.3) to (0.9,0.5);

\node at (0.65,0.9)  {$\vec{\Gamma}_{4}^{1}\left(\xi^{2},t\right)$};


\node at (1.40,1.85)  {$\boldsymbol{e_{2}}(t)$};

\draw[line width=1.5pt,bend left,-] (2, 1.3) to (2, 2.4);

\node at (1.45,1.6)  {$\vec{\Gamma}_{1}^{2}\left(\xi^{1},t\right)$};

\draw[line width=1.5pt,bend left,-] (2, 2.4) to (1, 2.2);  

\node at (1.45,2.36)  {$\vec{\Gamma}_{3}^{2}\left(\xi^{1},t\right)$};

\draw[line width=1.5pt,bend right,-] (1, 2.2) to (0.9, 1.3);   

\node at (2.25,1.85)  {$\vec{\Gamma}_{2}^{2}\left(\xi^{2},t\right)$};

\node at (0.45,1.84)  {$\vec{\Gamma}_{4}^{2}\left(\xi^{2},t\right)$};

 \draw[thick,->] (0,0,0) -- (2.75,0,0) node[anchor = north west]{$x_{1}$};
 \draw[thick,->] (0,0,0) -- (0,2.75,0) node[anchor = north west]{$x_{2}$};
\end{tikzpicture}
\caption{\label{ElementsSharingTheSameBoundary} Two elements $\boldsymbol{e_{1}}(t)$ and $\boldsymbol{e_{2}}(t)$ of a conforming mesh sharing the same curved boundary (dotted curve).}
\end{center}
\end{figure}

In the following, $e_{1}(t)$ and $e_{2}(t)$ are two neighboring elements which share the same boundary face. Without lost of generality the elements share the face $\vec{\Gamma}_{3}^{1}=\vec{\Gamma}_{1}^{2}$ as it is illustrated in Figure \ref{ElementsSharingTheSameBoundary}.
Then, for the elements $e_{l}(t)$, $l=1,2$, the grid velocity field is given by  
\begin{align}\label{IsoparametricGridvelocity}
\begin{split}
\vec{\nu}^{l}\left(\xi^{1},\xi^{2},t\right)=&\quad\frac{1}{2}\left[\left(1-\xi^{1}\right)\interpolation{N}{\left(\frac{d}{dt}\vec{\Gamma}_{4}^{l}\right)}\left(\xi^{2},t\right)+\left(1+\xi^{1}\right)\interpolation{N}{\left(\frac{d}{dt}\vec{\Gamma}_{2}^{l}\right)}\left(\xi^{2},t\right)\right]\\
&+\frac{1}{2}\left[\left(1-\xi^{2}\right)\interpolation{N}{\left(\frac{d}{dt}\vec{\Gamma}_{1}^{l}\right)}\left(\xi^{1},t\right)+\left(1+\xi^{2}\right)\interpolation{N}{\left(\frac{d}{dt}\vec{\Gamma}_{3}^{l}\right)}\left(\xi^{1},t\right)\right]\\
&-\frac{1}{4}\left[\left(1+\xi^{1}\right)\left\{ \left(1-\xi^{2}\right)\interpolation{N}{\left(\frac{d}{dt}\vec{\Gamma}_{1}^{l}\right)}\left(1,t\right)+\left(1+\xi^{2}\right)\interpolation{N}{\left(\frac{d}{dt}\vec{\Gamma}_{3}^{l}\right)\left(1,t\right)}\right\} \right]\\
&-\frac{1}{4}\left[\left(1-\xi^{1}\right)\left\{ \left(1-\xi^{2}\right)\interpolation{N}{\left(\frac{d}{dt}\vec{\Gamma}_{1}^{l}\right)}\left(-1,t\right)+\left(1+\xi^{2}\right)\interpolation{N}{\left(\frac{d}{dt}\vec{\Gamma}_{3}^{l}\right)}\left(-1,t\right)\right\} \right],
\end{split}
\end{align}
since it holds the identity 
\begin{equation}
\frac{d}{dt}\interpolation{N}{\left(\vec{\Gamma}_{i}^{l}\right)}=\interpolation{N}{\left(\frac{d}{dt}\vec{\Gamma}_{i}^{l}\right)},\qquad l=1,2, \qquad \text{and} \qquad  i=1,2,3,4.
\end{equation} 
Furthermore, since for $l=1,2$, and $i=1,2,3,4$, the faces $\vec{\Gamma}_{i}^{l}\left(\cdot,\cdot,t\right)$ are continuously differentiable in the time interval $\left[0,T\right]$, it holds    
\begin{align}
\interpolation{N}{\left(\frac{d}{dt}\vec{\Gamma}_{4}^{1}\right)\left(1,t\right)}=&\interpolation{N}{\left(\frac{d}{dt}\vec{\Gamma}_{3}^{1}\right)}\left(-1,t\right),\qquad\interpolation{N}{\left(\frac{d}{dt}\vec{\Gamma}_{2}^{1}\right)}\left(1,t\right)=\interpolation{N}{\left(\frac{d}{dt}\vec{\Gamma}_{3}^{1}\right)}\left(1,t\right), \label{AppendixAssumption1} \\
\interpolation{N}{\left(\frac{d}{dt}\vec{\Gamma}_{4}^{2}\right)}\left(-1,t\right)=&\interpolation{N}{\left(\frac{d}{dt}\vec{\Gamma}_{1}^{2}\right)}\left(-1,t\right),\qquad\interpolation{N}{\left(\frac{d}{dt}\vec{\Gamma}_{2}^{2}\right)}\left(-1,t\right)=\interpolation{N}{\left(\frac{d}{dt}\vec{\Gamma}_{1}^{2}\right)}\left(1,t\right), \label{AppendixAssumption2}
\end{align}
\begin{equation}\label{AppendixAssumption3}
\frac{d}{dt}\vec{\Gamma}_{3}^{1}\left(\zeta_{j},t\right)=\frac{d}{dt}\vec{\Gamma}_{1}^{2}\left(\zeta_{j},t\right),\qquad j=0,\dots,N,
\end{equation}
where $\left\{ \zeta_{j}\right\} _{j=0}^{N}$ are interpolation points.

For the element $e_{1}(t)$ the points along the interface with $e_{2}(t)$ are mapped in the set $\left\{ \left(\xi,1\right):\ \xi\in\left[-1,1\right]\right\}$. Hence, the grid velocity becomes    
\begin{align}
\begin{split}
\vec{\nu}^{1}\left(\xi,1,t\right)=&
\interpolation{N}{\left(\frac{d}{dt}\vec{\Gamma}_{3}^{1}\right)}\left(\xi,t\right),\qquad \forall\xi\in\left[-1,1\right]  
\end{split}
\end{align}
by \eqref{AppendixAssumption1}. On the opposite, for the element $e_{2}(t)$ the points along the interface with $e_{1}(t)$ are mapped in the set $\left\{ \left(\xi,-1\right):\ \xi\in\left[-1,1\right]\right\}$ and we obtain
\begin{align}
\begin{split}
\vec{\nu}^{2}\left(\xi,-1,t\right)=&
\interpolation{N}{\left(\frac{d}{dt}\vec{\Gamma}_{1}^{2}\right)}\left(\xi,t\right),\qquad \forall\xi\in\left[-1,1\right], 
\end{split}
\end{align}
by \eqref{AppendixAssumption2}. Thus, we obtain 
$\vec{\nu}^{1}\left(\cdot,1,t\right)=\vec{\nu}^{2}\left(\cdot,-1,t\right)$ by \eqref{AppendixAssumption3}. This proves that the grid velocity is continuous in the interface points of the two neighboring elements. 

\section{Entropy stable moving mesh shallow water fluxes}\label{sec:App Shallow Water Flux}
In this section, we apply the methodology in Remark 2.2 to construct two entropy conservative moving mesh fluxes for the shallow water equations. The homogeneous shallow water equations (without bottom topography) are given by  
\begin{equation}\label{eq:Shallow}
\pderivative{\textbf{u}}{t} + \vec{\nabla} \cdot \blockvec{\textbf{f}}=\textbf{0}. 
\end{equation}
The state vector and the components of the block vector flux, $\blockvec{\textbf{f}}$, are defined as 
\begin{equation}
\mathbf{u}=\begin{bmatrix}h\\
hu_{1}\\
hu_{2}
\end{bmatrix},\qquad\textbf{f}_{1}=\begin{bmatrix}hu_{1}\\
hu_{1}^{2}+\frac{g}{2}h^{2}\\
hu_{1}u_{2}
\end{bmatrix},\qquad\textbf{f}_{2}=\begin{bmatrix}hu_{2}\\
hu_{1}u_{2}\\
hu_{2}^{2}+\frac{g}{2}h^{2}
\end{bmatrix},
\end{equation}
where the conserved states are the water height and the discharge $h\vec{u}=\left[hu_{1},hu_{2}\right]^{T}$.  The quantity $g$ is the gravitational constant. The shallow water equations are equipped with the entropy/entropy flux pairs 
\begin{equation}\label{ShallowEntropy}
s=\frac{1}{2}h\left|\vec{u}\right|^{2}+\frac{1}{2}gh^{2},\qquad f_{l}^{s}=\frac{1}{2}hu_{l}\left|\vec{u}\right|^{2}+gh^{2}u_{l},\qquad l=1,2.
\end{equation}
The entropy variables are 
\begin{equation}\label{ShallowEntropyVariables}
\mathbf{w}=\left[gh-\frac{1}{2}\left|\vec{u}\right|^{2},u_{1},u_{2}\right]^{T}. 
\end{equation}
The entropy functionals are given by 
\begin{equation}\label{ShallowEntropyFunctionals}
\phi=\mathbf{w}^{T}\mathbf{u}-s=\frac{1}{2}gh^{2},\qquad\psi_{l}=\mathbf{w}^{T}\mathbf{f}_{l}-f_{l}^{s}=\frac{1}{2}gh^{2}u_{l},\qquad l=1,2.
\end{equation}

\subsection{Entropy conservative shallow water flux based on the flux in \cite{Wintermeyer2017}}\label{WGWK_Flux}
Friedrich et al. \cite[Appendix A, Equation (A.4)]{Friedrich2018} constructed the following state function     
\begin{equation}\label{ShallowState}
\mathbf{U}^{\#}=\begin{bmatrix}\avg{h}\\
\avg{h}\avg{u_{1}}\\
\avg{h}\avg{u_{2}}
\end{bmatrix}.
\end{equation}
The state function \eqref{ShallowState} is consistent with the shallow water state, symmetric and satisfies 
 \begin{equation}\label{TadmorShallow1}
\jump{\mathbf{w}}^{T}\mathbf{U}^{\#}=\frac{1}{2}g\jump{h^{2}}=\jump{\phi}.
\end{equation}
In addition, an entropy conservative numerical flux for the shallow water equations has been constructed by Wintermeyer, Winters, Gassner, Kopriva (WGWK) \cite{Wintermeyer2017}. The flux is given by   
\begin{equation}\label{WGWK-Flux}
\mathbf{F}_{1}^{\text{EC\_WGWK}}=\begin{bmatrix}\avg{hu_{1}}\\
\avg{hu_{1}}\avg{u_{1}}+g\avg{h}^{2}-\frac{1}{2}\avg{h^{2}}\\
\avg{hu_{1}}\avg{u_{2}}
\end{bmatrix},\qquad 
\mathbf{F}_{2}^{\text{EC\_WGWK}}=\begin{bmatrix}\avg{hu_{2}}\\
\avg{hu_{2}}\avg{u_{1}}\\
\avg{hu_{2}}\avg{u_{2}}+g\avg{h}^{2}-\frac{1}{2}\avg{h^{2}}
\end{bmatrix}.
\end{equation}
The flux \eqref{WGWK-Flux} is consistent, symmetric and satisfies       
\begin{equation}\label{TadmorShallow2}
\jump{\mathbf{w}}^{T}\mathbf{F}_{l}^{\text{EC\_WGWK}}=\frac{1}{2}g\jump{h^{2}u_{l}}=\jump{\psi_{l}},\qquad l=1,2.
\end{equation}
The state function \eqref{ShallowState} and the flux \eqref{WGWK-Flux} are used to construct the flux  
\begin{subequations}\label{MMWGWK-Flux}
\begin{equation}\label{ShallowECFlux_X-Direction}
\mathbf{G}_{1}^{\text{EC\_WGWK}}=\mathbf{F}_{1}^{\text{EC\_WGWK}}-\avg{\nu_{1}}\mathbf{U}^{\#}=\begin{bmatrix}\avg{hu_{1}}-\avg{h}\avg{\nu_{1}}\\
\avg{hu_{1}}\avg{u_{1}}-\avg{h}\avg{u_{1}}\avg{\nu_{1}}+g\avg{h}^{2}-\frac{1}{2}\avg{h^{2}}\\
\avg{hu_{1}}\avg{u_{2}}-\avg{h}\avg{u_{2}}\avg{\nu_{1}}
\end{bmatrix},
\end{equation}
and
\begin{equation}\label{ShallowECFlux_Y-Direction}
\mathbf{G}_{2}^{\text{EC\_WGWK}}=\mathbf{F}_{2}^{\text{EC\_WGWK}}-\avg{\nu_{2}}\mathbf{U}^{\#}=\begin{bmatrix}\avg{hu_{2}}-\avg{h}\avg{\nu_{2}}\\
\avg{hu_{2}}\avg{u_{1}}-\avg{h}\avg{u_{1}}\avg{\nu_{2}}\\
\avg{hu_{2}}\avg{u_{2}}-\avg{h}\avg{u_{2}}\avg{\nu_{2}}+g\avg{h}^{2}-\frac{1}{2}\avg{h^{2}}
\end{bmatrix}.
\end{equation}
\end{subequations}
The flux functions \eqref{ShallowECFlux_X-Direction},  \eqref{ShallowECFlux_Y-Direction} are consistent with $\mathbf{g}_{1}$, $\mathbf{g}_{2}$, symmetric and it follows by \eqref{TadmorShallow1} and \eqref{TadmorShallow2}    
\begin{equation}\label{TadmorShallow3}
\jump{\mathbf{w}}^{T}\mathbf{G}_{l}^{\text{EC\_WGWK}}=\jump{\mathbf{w}}^{T}\mathbf{F}_{l}^{\text{EC\_WGWK}}-\avg{\nu_{l}}\jump{\mathbf{w}}^{T}\mathbf{U}^{\#}=\jump{\psi_{l}}-\avg{\nu_{l}}\jump{\phi},\qquad l=1,2.
\end{equation}

\subsection{Entropy conservative shallow water flux based on the flux in \cite{Fjordholm2011}}\label{Fjordholm}
Fjordholm, Mishra, Tadmor (FMT) constructed a further entropy conservative numerical flux for the shallow water equations in \cite{Fjordholm2011}. This flux is given by 
\begin{equation}\label{FMT-Flux}
\mathbf{F}_{1}^{\text{EC\_FMT}}=\begin{bmatrix}\avg{h}\avg{u_{1}}\\
\avg{h}\avg{u_{1}}^{2}+\frac{g}{2}\avg{h^{2}}\\
\avg{h}\avg{u_{1}}\avg{u_{2}}
\end{bmatrix}, 
\qquad 
\mathbf{F}_{2}^{\text{EC\_FMT}}=\begin{bmatrix}\avg{h}\avg{u_{2}}\\
\avg{h}\avg{u_{1}}\avg{u_{2}}\\
\avg{h}\avg{u_{2}}^{2}+\frac{g}{2}\avg{h^{2}}
\end{bmatrix}.
\end{equation}
The state function \eqref{ShallowState} and the flux \eqref{FMT-Flux} are used to construct the flux  
\begin{equation}\label{ShallowEC_FMT}
\mathbf{G}_{1}^{\text{EC\_FMT}}=\begin{bmatrix}\avg{h}\avg{u_{1}-\nu_{1}}\\
\avg{h}\avg{u_{1}-\nu_{1}}\avg{u_{1}}+\frac{1}{2}\avg{h^{2}}\\
\avg{h}\avg{u_{1}-\nu_{1}}\avg{u_{2}}
\end{bmatrix}, \qquad 
\mathbf{G}_{2}^{\text{EC\_FMT}}=\begin{bmatrix}\avg{h}\avg{u_{2}-\nu_{2}}\\
\avg{h}\avg{u_{2}-\nu_{2}}\avg{u_{1}}\\
\avg{h}\avg{u_{2}-\nu_{2}}\avg{u_{2}}+\frac{1}{2}\avg{h^{2}}
\end{bmatrix}.
\end{equation}
The flux functions \eqref{ShallowEC_FMT} are consistent with $\mathbf{g}_{1}$, $\mathbf{g}_{2}$, symmetric and satisfy \eqref{TadmorShallow3}. 

\subsection{Matrix dissipation term for the shallow water flux}\label{ShallowDiss}
In the Section \ref{sec:EntropyStability}, it has been shown that the matrix dissipation operators to stabilize the entropy conservative moving mesh flux are basically the same as on a static mesh. Thus, we use the shallow water fluxes from the previous Appendices \ref{WGWK_Flux} and \ref{Fjordholm} with the dissipation operators from Wintermeyer et al. \cite{Wintermeyer2017}.
In the following, the matrices to construct the entropy stable dissipation operators \eqref{EntropyBasedDissipation} are listed. The scaling matrix for the right eigenvalues is given by 
\begin{equation}
\matx{T}=\text{diag}\left(\frac{1}{\sqrt{2g}},\sqrt{h},\frac{1}{\sqrt{2g}}\right)
\end{equation}
in both directions. Therefore, we apply the following discrete scaling matrix 
\begin{equation}
\matx{T}^{\star}=\text{diag}\left(\frac{1}{\sqrt{2g}},\sqrt{\avg{h}},\frac{1}{\sqrt{2g}}\right).
\end{equation}
The average components of the other matrices to compute the dissipation term in the $x_{1}$-direction are  
\begin{equation}
\matx{R}_{1}^{\star}=\begin{bmatrix}1 & 0 & 1\\
\avg{u_{1}}-\sqrt{g\avg{h}} & 0 & \avg{u_{1}}+\sqrt{g\avg{h}}\\
\avg{u_{2}} & 1 & \avg{u_{2}}
\end{bmatrix},
\end{equation}
\begin{equation}
\Lambda_{1}=\text{diag}\left(\left|\avg{u_{1}-\nu_{1}}-\sqrt{g\avg{h}}\right|,\left|\avg{u_{1}-\nu_{1}}\right|,\left|\avg{u_{1}-\nu_{1}}+\sqrt{g\avg{h}}\right|\right).
\end{equation}
In the $x_{2}$-direction the components are given by
\begin{equation}
\matx{R}_{2}^{\star}=\begin{bmatrix}1 & 0 & 1\\
\avg{u_{1}} & 1 & \avg{u_{1}}\\
\avg{u_{2}}-\sqrt{g\avg{h}} & 0 & \avg{u_{2}}+\sqrt{g\avg{h}}
\end{bmatrix},
\end{equation}
\begin{equation}
\Lambda_{2}=\text{diag}\left(\left|\avg{u_{2}-\nu_{2}}-\sqrt{g\avg{h}}\right|,\left|\avg{u_{2}-\nu_{2}}\right|,\left|\avg{u_{2}-\nu_{2}}+\sqrt{g\avg{h}}\right|\right).
\end{equation}

\section{Entropy stable moving mesh Euler fluxes}\label{sec:App Euler}
We present entropy stable Cartesian fluxes  $\mathbf{G}_{l}^{\text{EC}}$, $l=1,2,3$, for the compressible Euler equations \eqref{eq:Euler} equipped with the entropy/entropy flux pairs \eqref{EulerEntropy}. Then the entropy variables are given by  
\begin{equation}\label{EulerEntropyVariables}
\mathbf{w}=\left[\frac{\gamma-\varsigma}{\gamma-1}-\beta\left|\vec{u}\right|^{2},2\beta u_{1},2\beta u_{2},2\beta u_{3},-2\beta\right]^{T},\qquad\text{with}\qquad\beta:=\frac{\rho}{2p}
\end{equation}
and the entropy functionals are given by 
\begin{equation}\label{EulerEntropyFunctionals}
\phi=\mathbf{w}^{T}\mathbf{u}-s=\rho,\qquad\psi_{l}=\mathbf{w}^{T}\mathbf{f}_{l}-f_{l}^{s}=\rho u_{l},\qquad l=1,2,3.
\end{equation}

\subsection{Entropy conservative Euler flux based on the flux in \cite{Chandrashekar2013}}\label{Chandrashekar_Flux}
In the following the logarithmic mean $\avg{\cdot}^{\text{log}}$ will be used.     For two positive states $a^{-}$ and $a^{+}$, the logarithmic mean is defined by 
\begin{equation}\label{LogMean}
\avg{a}^{\text{log}}:=\begin{cases}
\frac{\jump{a}}{\jump{\log\left(a\right)}}, & \text{if }a^{-}\neq a^{+},\\
a^{-}, & \text{if }a^{-}=a^{+}.
\end{cases}
\end{equation}
A numerically stable procedure to compute the logarithmic mean \eqref{LogMean} is provided by Ismail and Roe \cite[Appendix B]{Ismail2009}. Friedrich et al. \cite[Theorem 3]{Friedrich2018} constructed the following state function 
\begin{equation}\label{EulerState}
\mathbf{U}^{\#}=\begin{bmatrix}\avg{\rho}^{\text{log}}\\
\avg{\rho}^{\text{log}}\avg{u_{1}}\\
\avg{\rho}^{\text{log}}\avg{u_{2}}\\
\avg{\rho}^{\text{log}}\avg{u_{3}}\\
\frac{\avg{\rho}^{\text{log}}}{2\left(\gamma-1\right)\avg{\beta}^{\text{log}}}+\frac{1}{2}\avg{\rho}^{\text{log}}\overline{\left|\vec{u}\right|^{2}}
\end{bmatrix},
\end{equation}
where 
\begin{equation}
\overline{\left|\vec{u}\right|^{2}}=2\left(\avg{u_{1}}^{2}+\avg{u_{2}}^{2}+\avg{u_{3}}^{2}\right)-\left(\avg{u_{1}^{2}}+\avg{u_{2}^{2}}+\avg{u_{3}^{2}}\right).
\end{equation}
The state function \eqref{EulerState} is consistent, symmetric and it holds 
 \begin{equation}\label{TadmorEuler1}
\jump{\mathbf{w}}^{T}\mathbf{U}^{\#}=\jump{\rho}=\jump{\phi}.
\end{equation}
Furthermore, Chandrashekar constructed in \cite{Chandrashekar2013} a kinetic energy preserving and entropy conservative (KEPEC) numerical flux function for the compressible Euler equations. In the $x_{1}$-direction Chandrashekar's KEPEC flux is given by   
\begin{equation}\label{ChandrashekarX-Direction}
\mathbf{F}_{1}^{\text{EC\_CH}}=\begin{bmatrix}\avg{\rho}^{\text{log}}\avg{u_{1}}\\
\avg{\rho}^{\text{log}}\avg{u_{1}}^{2}+\frac{\avg{\rho}}{2\avg{\beta}}\\
\avg{\rho}^{\text{log}}\avg{u_{1}}\avg{u_{2}}\\
\avg{\rho}^{\text{log}}\avg{u_{1}}\avg{u_{3}}\\
\frac{\avg{\rho}^{\text{log}}\avg{u_{1}}}{2\left(\gamma-1\right)\avg{\beta}^{\text{log}}}+\frac{1}{2}\avg{\rho}^{\text{log}}\avg{u_{1}}\overline{\left|\vec{u}\right|^{2}}+\frac{\avg{\rho}\avg{u_{1}}}{2\avg{\beta}}
\end{bmatrix}.
\end{equation}
The flux \eqref{ChandrashekarX-Direction} is consistent with $\textbf{f}_1$ and symmetric. In particular, Chandrashekar proved that      
\begin{equation}\label{TadmorEuler}
\jump{\mathbf{w}}^{T}\mathbf{F}_{1}^{\text{EC\_CH}}=\jump{\rho u_{1}}=\jump{\psi_{1}}.
\end{equation}
The state function \eqref{EulerState} and the flux function \eqref{ChandrashekarX-Direction} are used to construct the flux 
\begin{equation}\label{EulerECFlux_X-Direction}
\mathbf{G}_{1}^{\text{EC\_CH}}=\mathbf{F}_{1}^{\text{EC\_CH}}-\avg{\nu_{1}}\mathbf{U}^{\#}=\begin{bmatrix}\avg{\rho}^{\text{log}}\avg{u_{1}-\nu_{1}}\\
\avg{\rho}^{\text{log}}\avg{u_{1}-\nu_{1}}\avg{u_{1}}+\frac{\avg{\rho}}{2\avg{\beta}}\\
\avg{\rho}^{\text{log}}\avg{u_{1}-\nu_{1}}\avg{u_{2}}\\
\avg{\rho}^{\text{log}}\avg{u_{1}-\nu_{1}}\avg{u_{3}}\\
\frac{\avg{\rho}^{\text{log}}\avg{u_{1}-\nu_{1}}}{2\left(\gamma-1\right)\avg{\beta}^{\text{log}}}+\frac{1}{2}\avg{\rho}^{\text{log}}\avg{u_{1}-\nu_{1}}\overline{\left|\vec{u}\right|^{2}}+\frac{\avg{\rho}\avg{u_{1}}}{2\avg{\beta}}
\end{bmatrix}.
\end{equation}
The flux \eqref{EulerECFlux_X-Direction} is consistent with $\mathbf{g}_{1}=\mathbf{f}_{1}-\nu_{1}\mathbf{u}$, symmetric and it follows 
\begin{equation}\label{TadmorEuler2}
\jump{\mathbf{w}}^{T}\mathbf{G}_{1}^{\text{EC\_CH}}=\jump{\mathbf{w}}^{T}\mathbf{F}_{1}^{\text{EC\_CH}}-\avg{\nu_{1}}\jump{\mathbf{w}}^{T}\mathbf{U}^{\#}=\jump{\psi_{1}}-\avg{\nu_{1}}\jump{\phi}
\end{equation}    
by \eqref{TadmorEuler1} and \eqref{TadmorEuler}. In the same way, the $x_{2}$-direction and the $x_{3}$-direction of Chandrashekar's KEPEC flux can be used to construct 
\begin{equation}\label{EulerECFlux_Y-Direction}
\mathbf{G}_{2}^{\text{EC\_CH}}=\begin{bmatrix}\avg{\rho}^{\text{log}}\avg{u_{2}-\nu_{2}}\\
\avg{\rho}^{\text{log}}\avg{u_{1}}\avg{u_{2}-\nu_{2}}\\
\avg{\rho}^{\text{log}}\avg{u_{2}}\avg{u_{2}-\nu_{2}}+\frac{\avg{\rho}}{2\avg{\beta}}\\
\avg{\rho}^{\text{log}}\avg{u_{2}-\nu_{2}}\avg{u_{3}}\\
\frac{\avg{\rho}^{\text{log}}\avg{u_{2}-\nu_{2}}}{2\left(\gamma-1\right)\avg{\beta}^{\text{log}}}+\frac{1}{2}\avg{\rho}^{\text{log}}\avg{u_{2}-\nu_{2}}\overline{\left|\vec{u}\right|^{2}}+\frac{\avg{\rho}\avg{u_{2}}}{2\avg{\beta}}
\end{bmatrix}
\end{equation}
and 
\begin{equation}\label{EulerECFlux_Z-Direction}
\mathbf{G}_{3}^{\text{EC\_CH}}=\begin{bmatrix}\avg{\rho}^{\text{log}}\avg{u_{3}-\nu_{3}}\\
\avg{\rho}^{\text{log}}\avg{u_{1}}\avg{u_{3}-\nu_{3}}\\
\avg{\rho}^{\text{log}}\avg{u_{2}}\avg{u_{3}-\nu_{3}}\\
\avg{\rho}^{\text{log}}\avg{u_{3}}\avg{u_{3}-\nu_{3}}+\frac{\avg{\rho}}{2\avg{\beta}}\\
\frac{\avg{\rho}^{\text{log}}\avg{u_{3}-\nu_{3}}}{2\left(\gamma-1\right)\avg{\beta}^{\text{log}}}+\frac{1}{2}\avg{\rho}^{\text{log}}\avg{u_{3}-\nu_{3}}\overline{\left|\vec{u}\right|^{2}}+\frac{\avg{\rho}\avg{u_{3}}}{2\avg{\beta}}
\end{bmatrix}.
\end{equation}
The fluxes \eqref{EulerECFlux_Y-Direction}, \eqref{EulerECFlux_Z-Direction}  are consistent with $\mathbf{g}_{2}=\mathbf{f}_{2}-\nu_{2}\mathbf{u}$, $\mathbf{g}_{3}=\mathbf{f}_{3}-\nu_{3}\mathbf{u}$, symmetric and satisfy 
\begin{align}\label{TadmorEuler3}
\jump{\mathbf{w}}^{T}\mathbf{G}_{2}^{\text{EC\_CH}}=&\jump{\mathbf{w}}^{T}\mathbf{F}_{2}^{\text{EC\_CH}}-\avg{\nu_{2}}\jump{\mathbf{w}}^{T}\mathbf{U}^{\#}=\jump{\psi_{2}}-\avg{\nu_{2}}\jump{\phi}, \\
\jump{\mathbf{w}}^{T}\mathbf{G}_{3}^{\text{EC\_CH}}=&\jump{\mathbf{w}}^{T}\mathbf{F}_{3}^{\text{EC\_CH}}-\avg{\nu_{3}}\jump{\mathbf{w}}^{T}\mathbf{U}^{\#}=\jump{\psi_{3}}-\avg{\nu_{3}}\jump{\phi}.
\end{align}

\subsection{Entropy conservative Euler flux based on the flux in \cite{Ranocha2018}}\label{Ranocha_Flux}
Ranocha constructed in his PHD thesis \cite{Ranocha2018} another KEPEC numerical flux for the compressible Euler equations. We proceed as in the Appendix  \ref{Chandrashekar_Flux} and use the state \eqref{EulerState} and Ranocha's flux to construct the following two-point flux functions    
\begin{subequations}\label{EulerECFlux2}
\begin{equation}\label{EulerECFlux2_X-Direction}
\mathbf{G}_{1}^{\text{EC\_R}}=\begin{bmatrix}\avg{\rho}^{\text{log}}\avg{u_{1}-\nu_{1}}\\
\avg{\rho}^{\text{log}}\avg{u_{1}-\nu_{1}}\avg{u_{1}}+\avg{p}\\
\avg{\rho}^{\text{log}}\avg{u_{1}-\nu_{1}}\avg{u_{2}}\\
\avg{\rho}^{\text{log}}\avg{u_{1}-\nu_{1}}\avg{u_{3}}\\
\left(\frac{\avg{\rho}^{\text{log}}}{2\left(\gamma-1\right)\avg{\beta}^{\text{log}}}+\frac{1}{2}\avg{\rho}^{\text{log}}\right)\avg{u_{1}-\nu_{1}}\overline{\left|\vec{u}\right|^{2}}+2\avg{p}\avg{u_{1}}-\avg{pu_{1}}
\end{bmatrix},
\end{equation}
\begin{equation}\label{EulerECFlux2_Y-Direction}
\mathbf{G}_{2}^{\text{EC\_R}}=\begin{bmatrix}\avg{\rho}^{\text{log}}\avg{u_{2}-\nu_{2}}\\
\avg{\rho}^{\text{log}}\avg{u_{1}}\avg{u_{2}-\nu_{2}}\\
\avg{\rho}^{\text{log}}\avg{u_{2}}\avg{u_{2}-\nu_{2}}+\avg{p}\\
\avg{\rho}^{\text{log}}\avg{u_{2}-\nu_{2}}\avg{u_{3}}\\
\left(\frac{\avg{\rho}^{\text{log}}}{2\left(\gamma-1\right)\avg{\beta}^{\text{log}}}+\frac{1}{2}\avg{\rho}^{\text{log}}\right)\avg{u_{2}-\nu_{2}}\overline{\left|\vec{u}\right|^{2}}+2\avg{p}\avg{u_{2}}-\avg{pu_{2}}
\end{bmatrix},
\end{equation}
and 
\begin{equation}\label{EulerECFlux3_Z-Direction}
\mathbf{G}_{3}^{\text{EC\_R}}=\begin{bmatrix}\avg{\rho}^{\text{log}}\avg{u_{3}-\nu_{3}}\\
\avg{\rho}^{\text{log}}\avg{u_{1}}\avg{u_{3}-\nu_{3}}\\
\avg{\rho}^{\text{log}}\avg{u_{2}}\avg{u_{3}-\nu_{3}}\\
\avg{\rho}^{\text{log}}\avg{u_{3}}\avg{u_{3}-\nu_{3}}+\avg{p}\\
\left(\frac{\avg{\rho}^{\text{log}}}{2\left(\gamma-1\right)\avg{\beta}^{\text{log}}}+\frac{1}{2}\avg{\rho}^{\text{log}}\right)\avg{u_{3}-\nu_{3}}\overline{\left|\vec{u}\right|^{2}}+2\avg{p}\avg{u_{3}}-\avg{pu_{3}}.
\end{bmatrix}.
\end{equation}
\end{subequations}
The flux functions \eqref{EulerECFlux2} are consistent with $\mathbf{g}_{1}$, $\mathbf{g}_{2}$, $\mathbf{g}_{3}$, symmetric and satisfy 
\begin{equation}\label{TadmorEuler4}
\jump{\mathbf{w}}^{T}\mathbf{G}_{l}^{\text{EC\_R}}=\jump{\psi_{l}}-\avg{\nu_{l}}\jump{\phi},\qquad l=1,2,3.
\end{equation}  

\subsection{Matrix dissipation term for the Euler flux}\label{EulerDiss}
We will use the Euler fluxes from the previous Appendices \ref{Chandrashekar_Flux} and \ref{Ranocha_Flux} with the dissipation operators form Winters et al. \cite{Winters2017}. In the following, the matrices to construct the entropy stable dissipation operators \eqref{EntropyBasedDissipation} are listed. The average components of the dissipation term in the $x_{1}$-direction are given by
\begin{equation}
\matx{R}_{1}^{\star}=\begin{bmatrix}1 & 1 & 0 & 0 & 1\\
\avg{u_{1}}-\bar{c} & \avg{u_{1}} & 0 & 0 & \avg{u_{1}}+\bar{c}\\
\avg{u_{2}} & \avg{u_{2}} & 1 & 0 & \avg{u_{2}}\\
\avg{u_{3}} & \avg{u_{3}} & 0 & 1 & \avg{u_{3}}\\
\bar{h}-\avg{u_{1}}\bar{c} & \frac{1}{2}\overline{\left|\vec{u}\right|^{2}} & \avg{u_{2}} & \avg{u_{3}} & \bar{h}+\avg{u_{1}}\bar{c}
\end{bmatrix},
\end{equation}
\begin{equation}
\matx{T}_{1}^{\star}=\text{diag}\left(\sqrt{\frac{\avg{\rho}^{\text{log}}}{2\gamma}},\sqrt{\frac{\left(\gamma-1\right)}{\gamma}\avg{\rho}^{\text{log}}},\sqrt{\frac{\avg{\rho}}{2\avg{\beta}}},\sqrt{\frac{\avg{\rho}}{2\avg{\beta}}},\sqrt{\frac{\avg{\rho}^{\text{log}}}{2\gamma}}\right),
\end{equation}
\begin{equation}
\Lambda_{1}=\text{diag}\left(\left|\avg{u_{1}-\nu_{1}}-\bar{c}\right|,\left|\avg{u_{1}-\nu_{1}}\right|,\left|\avg{u_{1}-\nu_{1}}\right|,\left|\avg{u_{1}-\nu_{1}}\right|,\left|\avg{u_{1}-\nu_{1}}+\bar{c}\right|\right),
\end{equation}
where 
\begin{equation}
\bar{c}:=\sqrt{\frac{\gamma\avg{\rho}}{2\avg{\rho}^{\text{log}}\avg{\beta}}}, \qquad 
\bar{h}:=\frac{\gamma}{2\left(\gamma-1\right)\avg{\beta}^{\text{log}}}+\frac{1}{2}\overline{\left|\vec{u}\right|^{2}}.
\end{equation}
In the $x_{2}$-direction the components are given by
\begin{equation}
\matx{R}_{2}^{\star}=\begin{bmatrix}1 & 0 & 1 & 0 & 1\\
\avg{u_{1}} & 1 & \avg{u_{1}} & 0 & \avg{u_{1}}\\
\avg{u_{2}}-\bar{c} & 0 & \avg{u_{2}} & 0 & \avg{u_{2}}+\bar{c}\\
\avg{u_{3}} & 0 & \avg{u_{3}} & 1 & \avg{u_{3}}\\
\bar{h}-\avg{u_{2}}\bar{c} & \avg{u_{1}} & \frac{1}{2}\overline{\left|\vec{u}\right|^{2}} & \avg{u_{3}} & \bar{h}+\avg{u_{2}}\bar{c}
\end{bmatrix},
\end{equation}
\begin{equation}
\matx{T}_{2}^{\star}=\text{diag}\left(\sqrt{\frac{\avg{\rho}^{\text{log}}}{2\gamma}},\sqrt{\frac{\avg{\rho}}{2\avg{\beta}}},\sqrt{\frac{\left(\gamma-1\right)}{\gamma}\avg{\rho}^{\text{log}}},\sqrt{\frac{\avg{\rho}}{2\avg{\beta}}},\sqrt{\frac{\avg{\rho}^{\text{log}}}{2\gamma}}\right),
\end{equation}
\begin{equation}
\Lambda_{2}=\text{diag}\left(\left|\avg{u_{2}-\nu_{2}}-\bar{c}\right|,\left|\avg{u_{2}-\nu_{2}}\right|,\left|\avg{u_{2}-\nu_{2}}\right|,\left|\avg{u_{2}-\nu_{2}}\right|,\left|\avg{u_{2}-\nu_{2}}+\bar{c}\right|\right),
\end{equation}
and in the $x_{3}$-direction the components are given by
\begin{equation}
\matx{R}_{3}^{\star}=\begin{bmatrix}1 & 0 & 0 & 1 & 1\\
\avg{u_{1}} & 1 & 0 & \avg{u_{1}} & \avg{u_{1}}\\
\avg{u_{2}} & 0 & 1 & \avg{u_{2}} & \avg{u_{2}}\\
\avg{u_{3}}-\bar{c} & 0 & 0 & \avg{u_{3}} & \avg{u_{3}}+\bar{c}\\
\bar{h}-\avg{u_{3}}\bar{c} & \avg{u_{1}} & \avg{u_{2}} & \frac{1}{2}\overline{\left|\vec{u}\right|^{2}} & \bar{h}+\avg{u_{3}}\bar{c}
\end{bmatrix},
\end{equation}
\begin{equation}
\matx{T}_{3}^{\star}=\text{diag}\left(\sqrt{\frac{\avg{\rho}^{\text{log}}}{2\gamma}},\sqrt{\frac{\avg{\rho}}{2\avg{\beta}}},\sqrt{\frac{\avg{\rho}}{2\avg{\beta}}},\sqrt{\frac{\left(\gamma-1\right)}{\gamma}\avg{\rho}^{\text{log}}},\sqrt{\frac{\avg{\rho}^{\text{log}}}{2\gamma}}\right),
\end{equation}
\begin{equation}
\Lambda_{3}=\text{diag}\left(\left|\avg{u_{3}-\nu_{3}}-\bar{c}\right|,\left|\avg{u_{3}-\nu_{3}}\right|,\left|\avg{u_{3}-\nu_{3}}\right|,\left|\avg{u_{3}-\nu_{3}}\right|,\left|\avg{u_{3}-\nu_{3}}+\bar{c}\right|\right).
\end{equation}


\section{Proofs of entropy conservation for advection terms}\label{sec:Proof:EC}
In this section, we apply the SBP formula or Abel transformation
\begin{equation}\label{ClassicSummationByParts}
\sum_{k=1}^{K}a_{k}\left(b_{k+1}-b_{k}\right)=a_{K}b_{K+1}-a_{1}b_{1}-\sum_{k=1}^{K-1}a_{k+1}\left(b_{k+1}-b_{k}\right)
\end{equation}  
for two sequences $\left\{ a_{k}\right\} _{k=1}^{K}$ and $\left\{ b_{k}\right\} _{k=1}^{K+1}$. Furthermore, we apply the following identities which result from the properties of the SBP operator $\matx{Q}$  
\begin{align}
\sum_{i,j=0}^{N}\matx{Q}_{ij}\jump{a}_{\left(i,j\right)}\avg{b}_{\left(i,j\right)}=& \quad
\sum_{i,j=0}^{N}\matx{Q}_{ij}a_{i}b_{j}-\left(a_{N}b_{N}-a_{0}b_{0}\right), \label{SBP1} \\
\sum_{i,j=0}^{N}\matx{Q}_{ij}\jump{a}_{\left(i,j\right)}\avg{b}_{\left(i,j\right)}\avg{c}_{\left(i,j\right)}
=& \quad\sum_{i,j=0}^{N}2\matx{Q}_{ij}a_{i}\avg{b}_{\left(i,j\right)}\avg{c}_{\left(i,j\right)}-\left(a_{N}b_{N}c_{N}-a_{0}b_{0}c_{0}\right), \label{SBP2} 
\end{align}
where $\left\{ a\right\} _{i=0}^{N}$, $\left\{ b\right\} _{i=0}^{N}$ and $\left\{ c\right\} _{i=0}^{N}$ are  generic nodal values. These identities can be proven in a similar way as the discrete split forms in Lemma 1 in \cite{Gassner2016}. Thus, we skip a proof in this paper.

\subsection{Proof for Theorem 2.1}\label{FV:EC}
We multiply the equation \eqref{FVCL} by the discrete entropy variables $\mathbf{w}_{k}$ and sum over all elements $I_{k}\left(t\right)$  
\begin{equation}\label{FV:TotalEntropyConservation1}
\sum_{k=1}^{K}\mathbf{w}_{k}^{T}\left(\frac{d \mathrm{J}_{k}\textbf{u}_{k}}{dt}\right)=-\sum_{k=1}^{K}\mathbf{w}_{k}^{T}\left(\mathbf{g}_{k+\frac{1}{2}}^{*}-\mathbf{g}_{k-\frac{1}{2}}^{*}\right).
\end{equation} 

The same arguments as in the computation of \eqref{Temporalpart1D} 
provide the identity   
\begin{align}\label{FV:TemporalEC1}
\begin{split}
\mathbf{w}_{k}^{T}\left(\frac{d \mathrm{J}_{k}\mathbf{u}_{k}}{dt}\right)
=&\quad\mathrm{J}_{k}\left(\frac{d}{dt}s_{k}\left(\mathbf{u}_{k}\right)\right)+\left(\frac{d\mathrm{J}_{k}}{dt}\right)\mathbf{w}_{k}^{T}\mathbf{u}_{k}\\
=&\quad\left(\frac{d}{dt}\mathrm{J}_{k}s_{k}\left(\mathbf{u}_{k}\right)\right)+\left(\frac{d\mathrm{J}_{k}}{dt}\right)\phi_{k} \\
=&\quad\frac{d}{dt}\left(\mathrm{J}_{k}s_{k}\left(\mathbf{u}_{k}\right)\right)+\left(\nu_{k+\frac{1}{2}}^{*}-\nu_{k-\frac{1}{2}}^{*}\right)\phi_{k},
\end{split}
\end{align}
where we used the D-GCL \eqref{FVGCL} in the last step. Next, we sum the equation \eqref{FV:TemporalEC1} over all elements $I_{k}\left(t\right)$ and obtain     
\begin{equation}\label{FV:TemporalEC2}
\sum_{k=1}^{K}\mathbf{w}_{k}^{T}\left(\frac{d \mathrm{J}_{k}\mathbf{u}_{k}}{dt}\right)=\frac{d}{dt}\sum_{k=1}^{K}\mathrm{J}_{k}s_{k}\left(\mathbf{u}_{k}\right)+\sum_{k=1}^{K}\left(\nu_{k+\frac{1}{2}}^{*}-\nu_{k-\frac{1}{2}}^{*}\right)\phi_{k}.
\end{equation}

By the entropy condition \eqref{FV:EntropyCondition} and the SBP formula \eqref{ClassicSummationByParts} follows    
\begin{align}\label{FV:SpatialEC1}
\begin{split}
\sum_{k=1}^{K}\mathbf{w}_{k}^{T}\left(\mathbf{g}_{k+\frac{1}{2}}^{*}-\mathbf{g}_{k-\frac{1}{2}}^{*}\right)
=& \quad \mathbf{w}_{K}^{T}\mathbf{g}_{K+\frac{1}{2}}^{*}-\mathbf{w}_{1}^{T}\mathbf{g}_{\frac{1}{2}}^{*} -\sum_{k=1}^{K-1}\left(\mathbf{w}_{k+1}-\mathbf{w}_{k}\right)^{T}\mathbf{g}_{k+\frac{1}{2}}^{*}\\
=& \quad  \mathbf{w}_{K}^{T}\mathbf{g}_{K+\frac{1}{2}}^{*}-\mathbf{w}_{1}^{T}\mathbf{g}_{\frac{1}{2}}^{*}-\sum_{k=1}^{K-1}\left(\psi_{k+1}-\psi_{k}\right)+\sum_{k=1}^{K-1}\nu_{k+\frac{1}{2}}^{*}\left(\phi_{k+1}-\phi_{k}\right).
\end{split}
\end{align}
The first sum in the last line is a telescope sum and it follows    
\begin{equation}\label{FV:SpatialEC2}
\sum_{k=1}^{K-1}\left(\psi_{k+1}-\psi_{k}\right)=\psi_{K}-\psi_{1}. 
\end{equation}
To evaluate the second sum we use again the SBP formula \eqref{ClassicSummationByParts} and obtain
\begin{align}\label{FV:SpatialEC3}
\begin{split}
\sum_{k=1}^{K-1}\nu_{k+\frac{1}{2}}^{*}\left(\phi_{k+1}-\phi_{k}\right) 
=& 
\sum_{k=1}^{K-1}\nu_{k+\frac{1}{2}}^{*}\left(\phi_{k+1}-\phi_{k}\right)+\sum_{k=1}^{K}\left(\nu_{k+\frac{1}{2}}^{*}-\nu_{k-\frac{1}{2}}^{*}\right)\phi_{k}-\sum_{k=1}^{K}\left(\nu_{k+\frac{1}{2}}^{*}-\nu_{k-\frac{1}{2}}^{*}\right)\phi_{k} \\
=& 
\nu_{K+\frac{1}{2}}^{*}\psi_{K}-\nu_{\frac{1}{2}}^{*}\psi_{1}-\sum_{k=1}^{K}\left(\nu_{k+\frac{1}{2}}^{*}-\nu_{k-\frac{1}{2}}^{*}\right)\phi_{k}.
\end{split}
\end{align}
Next, we plug the equations \eqref{FV:SpatialEC2}, \eqref{FV:SpatialEC3} in \eqref{FV:SpatialEC1} and obtain
\begin{align}\label{FV:SpatialEC6}
\begin{split}
\sum_{k=1}^{K}\mathbf{w}_{k}^{T}\left(\mathbf{g}_{k+\frac{1}{2}}^{*}-\mathbf{g}_{k-\frac{1}{2}}^{*}\right)
=&\quad -\sum_{k=1}^{K}\left(\nu_{k+\frac{1}{2}}^{*}-\nu_{k-\frac{1}{2}}^{*}\right)\phi_{k} \\
&\quad+\mathbf{w}_{K}^{T}\mathbf{g}_{K+\frac{1}{2}}^{*}-\psi_{K}+\nu_{K+\frac{1}{2}}^{*}\psi_{K}\\
&\quad-\left(\mathbf{w}_{1}^{T}\mathbf{g}_{\frac{1}{2}}^{*}-\psi_{1}+\nu_{\frac{1}{2}}^{*}\psi_{1}\right). 
\end{split}
\end{align}
Finally, we plug the equations \eqref{FV:TemporalEC2} and \eqref{FV:SpatialEC6} in \eqref{FV:TotalEntropyConservation1} and obtain 
\begin{equation}
\frac{d}{dt}\sum_{k=1}^{K}\mathrm{J}_{k}s\left(\mathbf{u}_{k}\right)=-\left(\mathbf{w}_{K}^{T}\mathbf{g}_{K+\frac{1}{2}}^{*}-\psi_{K}+\nu_{K+\frac{1}{2}}^{*}\psi_{K}\right)+\left(\mathbf{w}_{1}^{T}\mathbf{g}_{\frac{1}{2}}^{*}-\psi_{1}+\nu_{\frac{1}{2}}^{*}\psi_{1}\right).
\end{equation}

\subsection{Proof for equation \eqref{DGSEM:EntropyPreservation5}}\label{DGSEM_EC}
The flux $\blockvec{\mathbf{G}}{}^{\text{EC}}$ satisfies the symmetry property \eqref{FluxSymmetric} and the SBP property \eqref{SBP} provides 
\begin{equation}\label{DGSEM_VolumeCont1}
2\omega_{i}\mathcal{D}_{im}=2\mathcal{Q}_{im}=\mathcal{Q}_{im}-\mathcal{Q}_{mi}+\mathcal{B}_{im},\qquad i,m=0,\dots,N. 
\end{equation}
Thus, we obtain  
\begin{align}\label{DGSEM_VolumeCont2}
\begin{split}
\left\langle \Dprojection{N}\cdot\blockvec{\tilde{\mathbf{G}}}{}^{\text{EC}},\mathbf{W}\right\rangle_{N}  
=& \quad \sum_{j,k=0}^{N}\omega_{j}\omega_{k}\sum_{i,m,k=0}^{N}\mathcal{Q}_{im}\jump{\mathbf{W}}_{\left(i,m\right)j}^{T}\left(\blockvec{\mathbf{G}}{}^{\text{EC}}\left(\vec{\nu}_{ijk},\vec{\nu}_{mjk},\mathbf{U}_{ijk},\mathbf{U}_{mjk}\right)\cdot\avg{\mathrm{J}\vec{a}^{1}}_{\left(i,m\right)jk}\right) \\
&+\sum_{j,k=0}^{N}\omega_{j}\omega_{k}\sum_{i,m=0}^{N}\mathcal{B}_{im}\mathbf{W}_{ijk}^{T}\left(\blockvec{\mathbf{G}}{}^{\text{EC}}\left(\vec{\nu}_{ijk},\vec{\nu}_{mjk},\mathbf{U}_{ijk},\mathbf{U}_{mjk}\right)\cdot\avg{\mathrm{J}\vec{a}^{1}}_{\left(i,m\right)jk}\right) \\
&+\sum_{i,k=0}^{N}\omega_{i}\omega_{k}\sum_{j,m=0}^{N}\mathcal{Q}_{jm}\jump{\mathbf{W}}_{i\left(j,m\right),k}^{T}\left(\blockvec{\mathbf{G}}{}^{\text{EC}}\left(\vec{\nu}_{ijk},\vec{\nu}_{imk},\mathbf{U}_{ijk},\mathbf{U}_{imk}\right)\cdot\avg{\mathrm{J}\vec{a}^{2}}_{i\left(j,m\right)k}\right)\\
&+\sum_{i,k=0}^{N}\omega_{i}\omega_{k}\sum_{j,m=0}^{N}\mathcal{B}_{jm}\mathbf{W}_{ijk}^{T}\left(\blockvec{\mathbf{G}}{}^{\text{EC}}\left(\vec{\nu}_{ijk},\vec{\nu}_{imk},\mathbf{U}_{ijk},\mathbf{U}_{imk}\right)\cdot\avg{\mathrm{J}\vec{a}^{2}}_{i\left(j,m\right)k}\right) \\
&+\sum_{i,j=0}^{N}\omega_{i}\omega_{j}\sum_{j,m=0}^{N}\mathcal{Q}_{km}\jump{\mathbf{W}}_{ij\left(k,m\right)}^{T}\left(\blockvec{\mathbf{G}}{}^{\text{EC}}\left(\vec{\nu}_{ijk},\vec{\nu}_{ijm},\mathbf{U}_{ijk},\mathbf{U}_{ijm}\right)\cdot\avg{\mathrm{J}\vec{a}^{3}}_{ij\left(k,m\right)}\right)\\
&+\sum_{i,j=0}^{N}\omega_{i}\omega_{j}\sum_{j,m=0}^{N}\mathcal{B}_{km}\mathbf{W}_{ijk}^{T}\left(\blockvec{\mathbf{G}}{}^{\text{EC}}\left(\vec{\nu}_{ijk},\vec{\nu}_{ijm},\mathbf{U}_{ijk},\mathbf{U}_{ijm}\right)\cdot\avg{\mathrm{J}\vec{a}^{3}}_{ij\left(k,m\right)}\right)
\end{split}
\end{align}
by the same calculation as in \cite[Appendix C.1., Equations (C.4) and (C.5)]{Friedrich2018} or \cite[Appendix B.3., Equation (B.31)]{Gassner2017}. To evaluate the fluxes $\blockvec{\mathbf{G}}{}^{\text{EC}}$ at the element interfaces, we apply the consistence condition \eqref{ConsistentDGECFlux}, such that e.g. 
\begin{equation}\label{DGSEM_VolumeCont3}
\blockvec{\mathbf{G}}{}^{\text{EC}}\left(\vec{\nu}_{Njk},\vec{\nu}_{Njk},\mathbf{U}_{Njk},\mathbf{U}_{Njk}\right)=\blockvec{\mathbf{F}}\left(\mathbf{U}_{Njk}\right)-\avg{\vec{\nu}}_{\left(N,N\right)jk}\mathbf{U}_{Njk}=\blockvec{\mathbf{G}}_{Njk},\qquad j,k=0,\dots,N.
\end{equation}
This provides the identity    
\begin{align}\label{DGSEM_VolumeCont4}
\begin{split}
& \quad \sum_{j,k=0}^{N}\omega_{j}\omega_{k}\sum_{i,m=0}^{N}\mathcal{B}_{im}\mathbf{W}_{ijk}^{T}\left(\blockvec{\mathbf{G}}{}^{\text{EC}}\left(\vec{\nu}_{ijk},\vec{\nu}_{mjk},\mathbf{U}_{ijk},\mathbf{U}_{mjk}\right)\cdot\avg{\mathrm{J}\vec{a}^{1}}_{\left(i,m\right)jk}\right) \\
& +\sum_{i,k=0}^{N}\omega_{i}\omega_{k}\sum_{j,m=0}^{N}\mathcal{B}_{jm}\mathbf{W}_{ijk}^{T}\left(\blockvec{\mathbf{G}}{}^{\text{EC}}\left(\vec{\nu}_{ijk},\vec{\nu}_{imk},\mathbf{U}_{ijk},\mathbf{U}_{imk}\right)\cdot\avg{\mathrm{J}\vec{a}^{2}}_{i\left(j,m\right)k}\right) \\
& +\sum_{i,j=0}^{N}\omega_{i}\omega_{j}\sum_{k,m=0}^{N}\mathcal{B}_{km}\mathbf{W}_{ijk}^{T}\left(\blockvec{\mathbf{G}}{}^{\text{EC}}\left(\vec{\nu}_{ijk},\vec{\nu}_{ijm},\mathbf{U}_{ijk},\mathbf{U}_{ijm}\right)\cdot\avg{\mathrm{J}\vec{a}^{3}}_{ij\left(k,m\right)}\right)\\
=& \quad \sum_{j,k=0}^{N}\omega_{j}\omega_{k}\left[\mathbf{W}_{Njk}^{T}\left(\blockvec{\mathbf{G}}_{Njk}\cdot\left(\mathrm{J}\vec{a}^{1}\right)_{Njk}\right)-\mathbf{W}_{0jk}^{T}\left(\blockvec{\mathbf{G}}_{0jk}\cdot\left(\mathrm{J}\vec{a}^{1}\right)_{0jk}\right)\right] \\
& +\sum_{i,k=0}^{N}\omega_{i}\omega_{k}\left[\mathbf{W}_{iNk}^{T}\left(\blockvec{\mathbf{G}}_{iNk}\cdot\left(\mathrm{J}\vec{a}^{2}\right)_{iNk}\right)-\mathbf{W}_{i0k}^{T}\left(\blockvec{\mathbf{G}}_{i0k}\cdot\left(\mathrm{J}\vec{a}^{2}\right)_{i0k}\right)\right]\\ 
& +\sum_{i,j=0}^{N}\omega_{i}\omega_{j}\left[\mathbf{W}_{ijN}^{T}\left(\blockvec{\mathbf{G}}_{ijN}\cdot\left(\mathrm{J}\vec{a}^{3}\right)_{ijN}\right)-\mathbf{W}_{ij0}^{T}\left(\blockvec{\mathbf{G}}_{ij0}\cdot\left(\mathrm{J}\vec{a}^{3}\right)_{ij0}\right)\right] \\
=& \quad \int\limits _{\partial E,N}\mathbf{W}^{T}\left\{ \blockvec{\tilde{\mathbf{G}}}\cdot\hat{n}\right\} \,dS. 
\end{split}
\end{align}
Next, we investigate the first sum on the right hand side in the equation \eqref{DGSEM_VolumeCont2}. Since the fluxes $\mathbf{G}_{l}^{\text{EC}}$, $l=1,2,3$, satisfy the entropy condition \eqref{EntropyConditionDGECFlux}, follows   
\begin{align}\label{DGSEM_VolumeCont5}
\begin{split}
& \quad \sum_{j,k=0}^{N}\omega_{j}\omega_{k}\sum_{i,m=0}^{N}\mathcal{Q}_{im}\jump{\mathbf{W}}_{\left(i,m\right)jk}^{T}\left(\blockvec{\mathbf{G}}{}^{\text{EC}}\left(\vec{\nu}_{ijk},\vec{\nu}_{mjk},\mathbf{U}_{ijk},\mathbf{U}_{mjk}\right)\cdot\avg{\mathrm{J}\vec{a}^{1}}_{\left(i,m\right)jk}\right) \\
=& \quad \sum_{j,k=0}^{N}\omega_{j}\omega_{k}\sum_{i,m=0}^{N}\mathcal{Q}_{im}\jump{\mathbf{W}}_{\left(i,m\right)jk}^{T}\mathbf{G}_{1}^{\text{EC}}\left(\vec{\nu}_{ijk},\vec{\nu}_{mjk},\mathbf{U}_{ijk},\mathbf{U}_{mjk}\right)\avg{\mathrm{J}a_{1}^{1}}_{\left(i,m\right)jk} \\
& +\sum_{j,k=0}^{N}\omega_{j}\omega_{k}\sum_{i,m=0}^{N}\mathcal{Q}_{im}\jump{\mathbf{W}}_{\left(i,m\right)jk}^{T}\mathbf{G}_{2}^{\text{EC}}\left(\vec{\nu}_{ijk},\vec{\nu}_{mjk},\mathbf{U}_{ijk},\mathbf{U}_{mjk}\right)\avg{\mathrm{J}a_{2}^{1}}_{\left(i,m\right)jk} \\
& +\sum_{j,k=0}^{N}\omega_{j}\omega_{k}\sum_{i,m=0}^{N}\mathcal{Q}_{im}\jump{\mathbf{W}}_{\left(i,m\right)jk}^{T}\mathbf{G}_{3}^{\text{EC}}\left(\vec{\nu}_{ijk},\vec{\nu}_{mjk},\mathbf{U}_{ijk},\mathbf{U}_{mjk}\right)\avg{\mathrm{J}a_{3}^{1}}_{\left(i,m\right)jk}\\
=& \quad \sum_{j,k=0}^{N}\omega_{j}\omega_{k}\sum_{i,m=0}^{N}\mathcal{Q}_{im}\left(\jump{\Psi_{1}}_{\left(i,m\right)jk}-\avg{\nu_{1}}_{\left(i,m\right)jk}\jump{\Phi}_{\left(i,m\right)jk}\right)\avg{\mathrm{J}a_{1}^{1}}_{\left(i,m\right)jk} \\
&+\sum_{j,k=0}^{N}\omega_{j}\omega_{k}\sum_{i,m=0}^{N}\mathcal{Q}_{im}\left(\jump{\Psi_{2}}_{\left(i,m\right)jk}-\avg{\nu_{2}}_{\left(i,m\right)jk}\jump{\Phi}_{\left(i,m\right)jk}\right)\avg{\mathrm{J}a_{2}^{1}}_{\left(i,m\right)jk}  \\
&+\sum_{j,k=0}^{N}\omega_{j}\omega_{k}\sum_{i,m=0}^{N}\mathcal{Q}_{im}\left(\jump{\Psi_{3}}_{\left(i,m\right)jk}-\avg{\nu_{3}}_{\left(i,m\right)jk}\jump{\Phi}_{\left(i,m\right)jk}\right)\avg{\mathrm{J}a_{3}^{1}}_{\left(i,m\right)jk}.
\end{split}
\end{align}
For $l=1,2,3$, the SPB properties \eqref{SBP1} and \eqref{SBP2} provide 
\begin{align}\label{DGSEM_VolumeCont6}
\begin{split}
& \quad \sum_{j,k=0}^{N}\omega_{j}\omega_{k}\sum_{i,m=0}^{N}\mathcal{Q}_{im}\left(\jump{\Psi_{l}}_{\left(i,m\right)jk}-\avg{\nu_{l}}_{\left(i,m\right)jk}\jump{\Phi}_{\left(i,m\right)jk}\right)\avg{\mathrm{J}a_{l}^{1}}_{\left(i,m\right)jk} \\ 
=& -\sum_{j,k=0}^{N}\omega_{j}\omega_{k}\left[\left(\left(\Psi_{l}\right)_{Njk}-\left(\nu_{l}\right)_{Njk}\Phi_{Njk}\right)\left(\mathrm{J}a_{l}^{1}\right)_{Njk}-\left(\left(\Psi_{l}\right)_{0jk}-\left(\nu_{l}\right)_{0jk}\Phi_{0jk}\right)\left(\mathrm{J}a_{l}^{1}\right)_{0jk}\right] \\ 
&+\sum_{i,j,k=0}^{N}\omega_{i}\omega_{j}\omega_{k}\left(\Psi_{l}\right)_{ijk}\sum_{m=0}^{N}\mathcal{D}_{im}\left(\mathrm{J}a_{l}^{1}\right)_{mjk}-\sum_{i,j,k=0}^{N}\omega_{i}\omega_{j}\omega_{k}\Phi_{ijk}\sum_{m=0}^{N}2\,\mathcal{D}_{im}\avg{\nu_{l}}_{\left(i,m\right)jk}\avg{\mathrm{J}a_{l}^{1}}_{\left(i,m\right)jk}.
\end{split}
\end{align}
Hence, we obtain the identity 
\begin{align}\label{DGSEM_VolumeCont7}
\begin{split}
& \quad \sum_{j,k=0}^{N}\omega_{j}\omega_{k}\sum_{i,m=0}^{N}\mathcal{Q}_{im}\jump{\mathbf{W}}_{\left(i,m\right)jk}^{T}\left(\blockvec{\mathbf{G}}{}^{\text{EC}}\left(\vec{\nu}_{ijk},\vec{\nu}_{mjk},\mathbf{U}_{ijk},\mathbf{U}_{mjk}\right)\cdot\avg{\mathrm{J}\vec{a}^{1}}_{\left(i,m\right)jk}\right) \\ 
=&-\sum_{j,k=0}^{N}\omega_{j}\omega_{k}\left[\left(\left(\vec{\Psi}\right)_{Njk}-\left(\vec{\nu}\right)_{Njk}\Phi_{Njk}\right)\cdot\left(\mathrm{J}\vec{a}^{1}\right)_{Njk}-\left(\left(\vec{\Psi}\right)_{0jk}-\left(\vec{\nu}\right)_{0jk}\Phi_{0jk}\right)\left(\mathrm{J}\vec{a}^{1}\right)_{0jk}\right] \\
& +\sum_{i,j,k=0}^{N}\omega_{i}\omega_{j}\omega_{k}\left(\vec{\Psi}\right)_{ijk}\cdot\left(\sum_{m=0}^{N}\mathcal{D}_{im}\left(\mathrm{J}\vec{a}^{1}\right)_{mjk}\right) \\
&-\sum_{i,j,k=0}^{N}\omega_{i}\omega_{j}\omega_{k}\Phi_{ijk}\sum_{m=0}^{N}2\,\mathcal{D}_{im}\avg{\vec{\nu}}_{\left(i,m\right)jk}\cdot\avg{\mathrm{J}\vec{a}^{1}}_{\left(i,m\right)jk}.
\end{split}
\end{align}
By the same computation, the third sum on the right hand side in the equation \eqref{DGSEM_VolumeCont2}  becomes 
\begin{align}\label{DGSEM_VolumeCont8}
\begin{split}
& \quad\sum_{i,k=0}^{N}\omega_{i}\omega_{k}\sum_{j,m=0}^{N}\mathcal{Q}_{jm}\jump{\mathbf{W}}_{i\left(j,m\right)k}^{T}\left(\blockvec{\mathbf{G}}{}^{\text{EC}}\left(\vec{\nu}_{ijk},\vec{\nu}_{imk},\mathbf{U}_{ijk},\mathbf{U}_{imk}\right)\cdot\avg{\mathrm{J}\vec{a}^{2}}_{i\left(j,m\right)k}\right) \\
=& -\sum_{i,k=0}^{N}\omega_{i}\omega_{k}\left[\left(\left(\vec{\Psi}\right)_{iNk}-\left(\vec{\nu}\right)_{iNk}\Phi_{iNk}\right)\cdot\left(\mathrm{J}\vec{a}^{2}\right)_{iNk}-\left(\left(\vec{\Psi}\right)_{i0k}-\left(\vec{\nu}\right)_{i0k}\Phi_{i0k}\right)\left(\mathrm{J}\vec{a}^{2}\right)_{i0k}\right] \\
&+\sum_{i,j,k=0}^{N}\omega_{i}\omega_{j}\omega_{k}\left(\vec{\Psi}\right)_{ijk}\cdot\left(\sum_{m=0}^{N}\mathcal{D}_{jm}\left(\mathrm{J}\vec{a}^{2}\right)_{imk}\right) \\
&-\sum_{i,j,k=0}^{N}\omega_{i}\omega_{j}\omega_{k}\Phi_{ijk}\sum_{m=0}^{N}2\,\mathcal{D}_{jm}\avg{\vec{\nu}}_{i\left(j,m\right)k}\cdot\avg{\mathrm{J}\vec{a}^{2}}_{i\left(j,m\right)k}
\end{split}
\end{align}
and the sum next-to-last on the right hand side in the equation \eqref{DGSEM_VolumeCont2} becomes 
\begin{align}\label{DGSEM_VolumeCont8a}
\begin{split}
& \quad\sum_{i,j=0}^{N}\omega_{i}\omega_{j}\sum_{k,m=0}^{N}\mathcal{Q}_{km}\jump{\mathbf{W}}_{ij\left(k,m\right)}^{T}\left(\blockvec{\mathbf{G}}{}^{\text{EC}}\left(\vec{\nu}_{ijk},\vec{\nu}_{ijm},\mathbf{U}_{ijk},\mathbf{U}_{ijm}\right)\cdot\avg{\mathrm{J}\vec{a}^{3}}_{ij\left(k,m\right)}\right) \\
=& -\sum_{i,j=0}^{N}\omega_{i}\omega_{j}\left[\left(\left(\vec{\Psi}\right)_{ijN}-\left(\vec{\nu}\right)_{ijN}\Phi_{ijN}\right)\cdot\left(\mathrm{J}\vec{a}^{3}\right)_{ijN}-\left(\left(\vec{\Psi}\right)_{ij0}-\left(\vec{\nu}\right)_{ij0}\Phi_{ij0}\right)\left(\mathrm{J}\vec{a}^{3}\right)_{ij0}\right] \\
&+\sum_{i,j,k=0}^{N}\omega_{i}\omega_{j}\omega_{k}\left(\vec{\Psi}\right)_{ijk}\cdot\left(\sum_{m=0}^{N}\mathcal{D}_{km}\left(\mathrm{J}\vec{a}^{3}\right)_{ijm}\right) \\
&-\sum_{i,j,k=0}^{N}\omega_{i}\omega_{j}\omega_{k}\Phi_{ijk}\sum_{m=0}^{N}2\,\mathcal{D}_{km}\avg{\vec{\nu}}_{ij\left(k,m\right)}\cdot\avg{\mathrm{J}\vec{a}^{3}}_{ij\left(k,m\right)}.
\end{split}
\end{align}
The definition of the derivative projection operator \eqref{SpatialDerivativeProjectionOperator1} in the D-GCL \eqref{DGSEM:D-GLC} provides  
\begin{align}\label{DGSEM_VolumeCont9}
\begin{split}
&-\sum_{i,j,k=0}^{N}\omega_{i}\omega_{j}\omega_{k}\Phi_{ijk}\sum_{m=0}^{N}2\,\mathcal{D}_{im}\avg{\vec{\nu}}_{\left(i,m\right)jk}\cdot\avg{\mathrm{J}\vec{a}^{1}}_{\left(i,m\right)jk} \\
&-\sum_{i,j,k=0}^{N}\omega_{i}\omega_{j}\omega_{k}\Phi_{ijk}\sum_{m=0}^{N}2\,\mathcal{D}_{jm}\avg{\vec{\nu}}_{i\left(j,m\right)k}\cdot\avg{\mathrm{J}\vec{a}^{2}}_{i\left(j,m\right)k}\\
&-\sum_{i,j,k=0}^{N}\omega_{i}\omega_{j}\omega_{k}\Phi_{ijk}\sum_{m=0}^{N}2\,\mathcal{D}_{km}\avg{\vec{\nu}}_{ij\left(k,m\right)}\cdot\avg{\mathrm{J}\vec{a}^{3}}_{ij\left(k,m\right)}\\
=&-\sum_{i,j,k=0}^{N}\omega_{i}\omega_{j}\omega_{k}\Phi_{ijk}\left(\Dprojection{N}\cdot\vec{\tilde{\nu}}_{ijk}\right)=-\left\langle \Dprojection{N}\cdot\vec{\tilde{\nu}},\Phi\right\rangle _{N}.
\end{split}
\end{align}
Next, we plug the equations \eqref{DGSEM_VolumeCont4}, \eqref{DGSEM_VolumeCont7}, \eqref{DGSEM_VolumeCont8}, \eqref{DGSEM_VolumeCont8a} in the equation \eqref{DGSEM_VolumeCont2} and apply the identity \eqref{DGSEM_VolumeCont9}. This results in the identity  
\begin{align}\label{DGSEM_VolumeCont10}
\begin{split}
\left\langle \Dprojection{N}\cdot\blockvec{\tilde{\mathbf{G}}}{}^{\text{EC}},\mathbf{W}\right\rangle_{N}  
=& \quad \int\limits _{\partial E,N}\left[\mathbf{W}^{T}\left\{ \blockvec{\tilde{\mathbf{G}}}\cdot\hat{n}\right\} -\left(\vec{\tilde{\Psi}}-\vec{\tilde{\nu}}\Phi\right)\cdot\hat{n}\right]\,dS -\left\langle \Dprojection{N}\cdot\vec{\tilde{\nu}},\Phi\right\rangle _{N} \\
&+\sum_{i,j,k=0}^{N}\omega_{i}\omega_{j}\omega_{k}\left(\vec{\Psi}\right)_{ijk}\cdot\left(\sum_{m=0}^{N}\left(\mathcal{D}_{im}\left(\mathrm{J}\vec{a}^{1}\right)_{mjk}+\mathcal{D}_{jm}\left(\mathrm{J}\vec{a}^{2}\right)_{imk}+\mathcal{D}_{km}\left(\mathrm{J}\vec{a}^{3}\right)_{ijm}\right)\right). 
\end{split}
\end{align}
The definition of the contravariant vector flux functions \eqref{ContravariantFunction} and contravariant block vector flux functions \eqref{ContravariantBlockVector} provide the equality   
\begin{align}\label{DGSEM_VolumeCont12}
\begin{split}
\mathbf{W}^{T}\left\{ \blockvec{\tilde{\mathbf{G}}}\cdot\hat{n}\right\} -\left(\vec{\tilde{\Psi}}-\vec{\tilde{\nu}}\Phi\right)\cdot\hat{n}  
=& \quad \sum_{l,r=1}^{3}\mathrm{J}a_{l}^{r}\left(\mathbf{W}^{T}\mathbf{G}_{l}-\mathbf{\Psi}_{l}+\nu_{l}\Phi\right)\hat{n}^{r} \\
=&\quad\sum_{l,r=1}^{3}\mathrm{J}a_{l}^{r}\left(F_{l}^{s}-\nu_{l}S\right)\hat{n}^{r} \\
=& \quad \left(\hat{s}\vec{n}\right)\cdot\left(\vec{F}^{s}-\vec{\nu}S\right)=\tilde{F}_{\hat{n}}^{s}-\tilde{\nu}_{\hat{n}}S.
\end{split}
\end{align}
Therefore, the equation \eqref{DGSEM_VolumeCont10} simplifies to
\begin{equation}\label{DGSEM_VolumeCont13}
\left\langle \Dprojection{N}\cdot\blockvec{\tilde{\mathbf{G}}}{}^{\text{EC}},\mathbf{W}\right\rangle_{N} =\quad\int\limits _{\partial E,N}\left(\tilde{F}_{\hat{n}}^{s}-\tilde{\nu}_{\hat{n}}S\right)\,dS-\left\langle \Dprojection{N}\cdot\vec{\tilde{\nu}},\Phi\right\rangle _{N}, 
\end{equation}
since the discrete volume weighted contravariant vectors  $\mathrm{J}\vec{a}^{\alpha}$, $\alpha=1,2,3$, are computed by the conservative curl form \eqref{DiscreteContravariantVectors} and the discrete metric identities \eqref{DiscreteMetricIdentities} are satisfied in the LGL points. 
\end{document}